\documentclass{amsart}

\usepackage{amsfonts}
\usepackage{amssymb}
\usepackage{psfig}

\newcommand\B{{{\mathcal B}}}

\newcommand\R{{{\mathbf R}}}

\newcommand\Sys{{{\mathcal S}}}

\newcommand\Z{{{\mathbf Z}}}
\newcommand\Y{{{\mathbf Y}}}
\newcommand\X{{{\mathbf X}}}
\newcommand\D{{{\mathcal D}}}
\newcommand\E{{{\mathbf E}}}
\newcommand\W{{{\mathbf W}}}
\newcommand\bfS{{{\mathbf \bfS}}}
\newcommand\x{{{\mathbf x}}}

\newcommand\m{{{\mathbf m}}}

\newcommand\h{{{\mathbf h}}}

\renewcommand\P{{{\mathcal P}}}
\newcommand\primes{{{\mathbf P}}}

\newcommand\lcm{{{\operatorname{lcm}}}}
\newcommand\mes{{{\operatorname{mes}}}}
\renewcommand\mod{{{\operatorname{mod}}}}
\newcommand\EXP{{{\operatorname{Exp}}}}
\newcommand\Res{{{\operatorname{Res}}}}
\renewcommand\Re{{{\operatorname{Re}}}}
\renewcommand\Im{{{\operatorname{Im}}}}
\newcommand\eps{\varepsilon}
\def \endprf{\hfill  {\vrule height6pt width6pt depth0pt}\medskip}

\theoremstyle{plain}
  \newtheorem{theorem}[subsection]{Theorem}

  \newtheorem{proposition}[subsection]{Proposition}
  
  \newtheorem{lemma}[subsection]{Lemma}
  \newtheorem{corollary}[subsection]{Corollary}

\theoremstyle{remark}
  \newtheorem{remark}[subsection]{Remark}
  \newtheorem{remarks}[subsection]{Remarks}
  \newtheorem{example}[subsection]{Example}
  \newtheorem{examples}[subsection]{Examples}

\theoremstyle{definition}
  \newtheorem{definition}[subsection]{Definition}

\begin{document}

\title[Polynomial progressions in primes]{The primes contain arbitrarily long polynomial progressions}

\author{Terence Tao}
\address{Department of Mathematics, UCLA, Los Angeles CA 90095-1555}
\email{tao@@math.ucla.edu}

\author{Tamar Ziegler}
\address{Department of Mathematics, The University of Michigan, Ann Arbor MI 48109}
\email{tamarz@@umich.edu}

\thanks{The second author was partially supported by NSF grant  No. DMS-0111298.
This work was initiated at a workshop held at the CRM in Montreal. The authors would like to
thank the CRM for their hospitality.}
\begin{abstract}  
We establish the existence of infinitely many \emph{polynomial} progressions in the primes; more precisely, given any integer-valued polynomials $P_1, \ldots, P_k \in \Z[\m]$ in one unknown $\m$ with $P_1(0) = \ldots = P_k(0) = 0$ and any $\eps > 0$, we show that there are infinitely many integers $x,m$ with $1 \leq m \leq x^\eps$ such that $x+P_1(m), \ldots, x+P_k(m)$ are simultaneously prime.  The arguments are based on those in \cite{gt-primes}, which treated the linear case $P_i = (i-1)\m$ and $\eps=1$; the main new features are a localization of the shift parameters (and the attendant Gowers norm objects) to both coarse and fine scales, the use of PET induction to linearize the polynomial averaging, and some elementary estimates for the number of points over finite fields in certain algebraic varieties.
\end{abstract}

\maketitle
\setcounter{tocdepth}{1}
\tableofcontents

\section{Introduction}

In 1975, Szemer\'edi \cite{szemeredi} proved that any subset $A$ of integers of positive upper density $\limsup_{N \rightarrow \infty} \frac{|A \cap [N]|}{|[N]|} > 0$ contains arbitrarily long arithmetic progressions. Throughout this paper $[N]$ denotes the discrete interval $[N] := \{1,\ldots,N\}$, and $|X|$ denotes the cardinality of a finite set $X$.
Shortly afterwards, Furstenberg \cite{furstenberg} gave an ergodic-theory proof of Szemer\'edi's
theorem. Furstenberg observed that questions about configurations in subsets of positive density in the integers correspond to recurrence questions for sets of positive measure in a probability measure
preserving system.  This observation is now known as the \emph{Furstenberg correspondence principle}.

In 1978, S\'ark\"ozy \cite{sarkozy}\footnote{S\'ark\"ozy actually proved a stronger theorem for the polynomial $P = \m^2$
providing an upper bound for density of a set $A$ for which $A-A$ does not contain a perfect 
square. His estimate was later improved by Pintz, Steiger, and Szemer\'edi in \cite{PSS}, and then generalized in \cite{BPSS} for $P = \m^k$ and then \cite{slijepcevic} for arbitrary $P$ with $P(0)=0$.} (using the Hardy-Littlewood circle method) and Furstenberg \cite{furstenberg-book} (using the correspondence principle, and ergodic theoretic methods) proved independently that for any polynomial\footnote{We use $\Z[\m]$ to denote the space of polynomials of one variable $\m$ with integer-valued coefficients; see Section \ref{notation-sec} for further notation along these lines.} $P \in \Z[\m]$ with $P(0)=0$,
any  set $A \subset \Z$ of positive density contains a pair of points $x,y$ with difference $y-x=P(m)$ for some positive integer $m \ge 1$.
In $1996$  Bergelson and Leibman \cite{bl} proved, by purely ergodic theoretic means\footnote{Unlike Szemer\'edi's theorem or S\'ark\"ozy's theorem, no non-ergodic proof of the Bergelson-Leibman theorem in its full generality is currently known. However, in this direction Green \cite{green} has shown by Fourier-analytic methods that any set of integers of positive density contains a triple $\{x,x+n,x+2n\}$ where $n$ is a non-zero sum of two squares.}, a vast generalization of the 
Fustenberg-S\'ark\"ozy  theorem -  establishing the existence of arbitrarily long polynomial progressions in sets of positive density in the integers.

\begin{theorem}[Polynomial Szemer\'edi theorem]\label{pSZ}\cite{bl} 
Let $A \subset \Z$ be a set of positive upper density, i.e.
$\limsup_{N \rightarrow \infty} \frac{|A \cap [N]|}{|[N]|} >0$.  Then given any integer-valued polynomials $P_1, \ldots, P_k \in \Z[\m]$ in one unknown $\m$ with $P_1(0) = \ldots = P_k(0) = 0$, $A$ contains infinitely many progressions 
of the form $x+P_1(m), \ldots, x+P_k(m)$ with $m > 0$.
\end{theorem} 

\begin{remark} By shifting $x$ appropriately, one may assume without loss of generality that one of the polynomials $P_i$ vanishes, e.g. $P_1 = 0$.  We shall rely on this ability to normalize one polynomial of our choosing to be zero at several points in the proof, most notably in the ``PET induction'' step in Section \ref{pet-sec}.  The arguments in \cite{bl} also establish a generalization of this theorem to higher dimensions, which will be important to us to obtain a certain uniformly quantitative version of this theorem later (see Theorem \ref{pSZ-quant} and Appendix \ref{uniform-sec}).
\end{remark}

The ergodic theoretic methods, to this day, have the limitation of only being able to 
handle sets of positive density in the integers, although this density is allowed to be arbitrarily small. However in $2004$, Green and Tao \cite{gt-primes} discovered a \emph{transference principle} which allowed one (at least in principle) to reduce questions about 
configurations in special sets of zero density (such as the primes $\P := \{2,3,5,7,\ldots\}$) to questions about 
sets of positive density in the integers. This opened the door to transferring the Szemer\'edi type results which are known for sets of positive upper density in the integers to the prime numbers.  Applying this transference principle to Szemer\'edi's theorem, Green and Tao showed that there are arbitrarily long arithmetic progressions in the prime numbers\footnote{Shortly afterwards, the transference principle was also combined in \cite{tao:gaussian} with the multidimensional Szemer\'edi theorem \cite{fk0} (or more precisely a hypergraph lemma related to this theorem, see \cite{tao:hyper}) to establish arbitrarily shaped constellations in the \emph{Gaussian} primes.  A much simpler transference principle is also available for dense subsets of \emph{genuinely random} sparse sets; see \cite{tao-transference}.}.

In this paper we prove a transference principle for polynomial configurations, which then allows us to use (a uniformly quantitative version of) the Bergelson-Leibman theorem to prove the existence of arbitrarily long \emph{polynomial} progressions in the primes, or more generally in large subsets of the primes.   More precisely, the main result of this paper is the following.

\begin{theorem}[Polynomial Szemer\'edi theorem for the primes]\label{mainthm} Let $A \subset \P$ be a set of primes of positive relative upper density in the primes, i.e.
$\limsup_{N \rightarrow \infty} \frac{|A \cap [N]|}{|\P \cap [N]|} > 0$.  Then given any $P_1, \ldots, P_k \in \Z[\m]$ with $P_1(0) = \ldots = P_k(0) = 0$, $A$ contains infinitely many progressions of the form $x+P_1(m), \ldots, x+P_k(m)$ with $m > 0$.
\end{theorem}

\begin{remarks} The main result of \cite{gt-primes} corresponds to Theorem \ref{mainthm} in the linear case $P_i := (i-1)\m$.  The case $k=2$ of this theorem follows from the results of \cite{PSS}, \cite{BPSS}, \cite{slijepcevic}, which in fact address arbitrary sets of integers with logarithmic type sparsity, and whose proof is more direct, proceeding via the Hardy-Littlewood circle method and not via the transference principle.  As a by-product of our proof, we shall also be able to impose the bound $m \leq x^\eps$ for any fixed $\eps > 0$, and thus (by diagonalization) that $m = x^{o(1)}$; see Remark \ref{short-rem}.  Our results for the case $A = \P$ are consistent with what is predicted by the Bateman-Horn conjecture \cite{bh}, which remains totally open in general (though see \cite{linprimes} for some partial progress in the linear case).
\end{remarks}

\begin{remark} In view of the generalization of Theorem \ref{pSZ} to higher dimensions in \cite{bl} it is reasonable to conjecture that an analogous result to Theorem \ref{mainthm} also holds in higher dimensions, thus any subset of $\P^d$ of positive relative upper density should contain infinitely many polynomial constellations, for any choice of polynomials which vanish at the origin.  This is however still open even in the linear case, the key difficulty being that the tensor product of pseudorandom measures is not pseudorandom.  In view of \cite{tao:gaussian} however, it should be possible (though time-consuming) to obtain a counterpart to Theorem \ref{mainthm} for the Gaussian primes.
\end{remark}

\begin{remark} The arguments in this paper are mostly quantitative and finitary, and in particular avoid the use of the axiom of choice.  However, our proof relies crucially on the Bergelson-Leibman theorem, Theorem \ref{pSZ}, and more precisely on a certain multidimensional generalization of that theorem \cite[Theorem $A_0$]{bl}.  At present, the only known proof of that theorem (in \cite{bl}) requires Zorn's lemma and thus our results here are also currently dependent on the axiom of choice. However, it is expected that the Bergelson-Leibman theorem will eventually be proven by other means which do not require the axiom of choice; for instance, the one-dimensional version of this theorem (i.e. Theorem \ref{pSZ}) can be established via the machinery of characteristic factors and Gowers-Host-Kra seminorms, by modifying the arguments in \cite{nikos-kra}, \cite{host-kra}, and this does not require the axiom of choice; this already allows us to establish Theorem \ref{mainthm} without the axiom of choice in the homogeneous case when $P_i(\m) = c_i \m^d$ for all $i=1,\ldots,k$ and some constants $c_1,\ldots,c_k, d$ (since in this case the $W$ factor in Theorem \ref{pSZ-quant} can be easily eliminated without introducing additional dimensions).  In a similar spirit, our arguments do not currently provide any effective bound for the first appearance of a pattern $x+P_1(m), \ldots, x+P_k(m)$ in the set $A$, but one expects that the Bergelson-Leibman theorem will eventually be proven with an effective bound (e.g. by extending the arguments in \cite{gowers}), in which case Theorem \ref{mainthm} will automatically come with an effective bound also.
\end{remark}

The philosophy of the proof is similar to the one in \cite{gt-primes}. The first key idea is to think of the primes (or any large subset thereof) as a set of positive relative density in the set of \emph{almost primes}, which (after some application of sieve theory, as in the work of Goldston and Y{\i}ld{\i}r{\i}m \cite{goldston-yildirim-old1}) can be shown to exhibit a somewhat pseudorandom behavior. Actually,
for technical reasons it is more convenient to work not with the sets of primes and almost primes, but rather with certain normalized weight functions $0 \leq f \leq \nu$ which are\footnote{This is an oversimplification, ignoring the ``$W$-trick'' necessary to eliminate local obstructions to uniformity; see Section \ref{notation-sec} for full details.} supported (or concentrated) on the primes and almost primes respectively, with $\nu$ obeying certain \emph{pseudorandom measure}\footnote{The term measure is a bit misleading. It is better to think of $\nu$ as the Radon-Nikodym derivative of a measure. Still, we stick to this terminology so as not to confuse the reader who is familiar with \cite{gt-primes}} properties.  The functions $f, \nu$ are unbounded, but have bounded expectation (mean).  A major step in the argument is then the establishment of a \emph{Koopman-von Neumann type structure theorem} which decomposes $f$ (except for a small error) as a sum $f = f_{U^\perp} + f_U$, where $f_{U^\perp}$ is a non-negative \emph{bounded} function with large expectation, and $f_U$ is an error which is unbounded but is so \emph{uniform} (in a Gowers-type sense) that it has a negligible impact on the (weighted) count of polynomial progressions.  The remaining component $f_{U^\perp}$ of $f$, being bounded, non-negative, and of large mean, can then be handled by (a quantitative version of) the Bergelson-Leibman theorem.

\begin{remark} As remarked in \cite{gt-primes}, the above transference arguments can be categorized as a kind of ``finitary'' ergodic theory. In the language of traditional (infinitary) ergodic theory, $f_{U^\perp}$ is analogous to a conditional expectation of $f$ relative to a suitable \emph{characteristic factor} for the polynomial average being considered.  Based on this analogy, and on the description of this characteristic factor in terms of nilsystems (see \cite{host-kra}, \cite{leibman}), one would hope that $f_{U^\perp}$ could be constructed out of nilsequences. In the case of linear averages this correspondence has already some roots in reality; see \cite{inverseU3}. In the special case $A=\P$ one can then hope to use analytic number theory methods to show that $f_{U^\perp}$ is essentially constant, which would lead to a more precise version of Theorem \ref{mainthm} in which one obtains a precise \emph{asymptotic} for the number of polynomial progressions in the primes, with $x$ and $n$ confined to various ranges.
In the case of progressions of length $4$ (or for more general linear patterns, assuming certain unproven conjectures), such an asymptotic was already established in
\cite{linprimes}.  While we expect similar asymptotics to hold for polynomial progressions, we do not pursue this interesting question here\footnote{One fundamental new difficulty that arises in the polynomial case is that it seems that one needs to control \emph{short} correlation sums between primes and nilsequences, such as on intervals of the form $[x, x + x^\eps]$, instead of the long correlation sums (such as on $[x,2x]$) which appear in the linear theory.  Even assuming strong conjectures such as GRH, it is not clear how to obtain such control.}.
\end{remark}

As we have already mentioned, the proof of Theorem \ref{mainthm} closely follows the arguments in \cite{gt-primes}.  However, some significant new difficulties arise when adapting those arguments\footnote{If the measure $\nu$ for the almost primes enjoyed infinitely many pseudorandomness conditions then one could adapt the arguments in \cite{tao-transference} to obtain Theorem \ref{mainthm} rather quickly.  Unfortunately, in order for $f$ to have non-zero mean, one needs to select a moderately large sieve level $R = N^{\eta_2}$ for the measure $\nu$, which means that one can only impose finitely many (though arbitarily large) such pseudorandomness conditions on $\nu$.  This necessitates the use of the (lengthier) arguments in \cite{gt-primes} rather than \cite{tao-transference}.} to the polynomial setting.  The most fundamental such difficulty arises in one of the very first steps of the argument in \cite{gt-primes}, in which one localizes the pattern $x + P_1(m), \ldots, x + P_k(m)$ to a finite interval $[N] = \{1,\ldots,N\}$.  In the linear case $P_i = (i-1)\m$ this localization restricts both $x$ and $m$ to be of size $O(N)$.  However, in the polynomial case, while the base point $x$ is still restricted to size $O(N)$, the shift parameter $m$ is now restricted to a much smaller range $O(M)$, where $M := N^{\eta_0}$ and $0 < \eta_0 < 1$ is a small constant depending on $P_1,\ldots,P_k$ (one can take for instance $\eta_0 := 1/2d_*$, where $d_*$ is the largest degree of the polynomials $P_1,\ldots,P_k$).  This eventually forces us to deal with localized averages of the form\footnote{This is an oversimplification as we are ignoring the need to first invoke the ``$W$-trick'' to eliminate local obstructions from small moduli and thus ensure that the almost primes behave pseudorandomly.  See Theorem \ref{thmquant} for the precise claim we need.}
\begin{equation}\label{pn-gowers}
 \E_{m \in [M]} \int_X T^{P_1(m)} f \ldots T^{P_k(m)} f
\end{equation}
where $X := \Z_N := \Z/N\Z$ is the cyclic group with $N$ elements,
$f: X \to \R^+$ is a weight function associated to the set $A$, and $T g(x) := g(x-1)$ is the shift operator on $X$.  Here we use the ergodic theory-like\footnote{Traditional ergodic theory would deal with the case where the underlying measure space $\Z_N$ is infinite and the shift range $M$ is going to infinity, thus informally $N=\infty$ and $M \to \infty$.  Unraveling the Furstenberg correspondence principle, this is equivalent to the setting where $N$ is finite (but going to infinity) and $M = \omega(N)$ is a very slowly growing function of $N$.    In \cite{gt-primes} one is instead working in the regime where $M=N$ are going to infinity at the \emph{same} rate.  The situation here is thus an intermediate regime where $M = N^{\eta_0}$ goes to infinity at a polynomially slower rate than $N$.  In the linear setting, all of these regimes can be equated using the random dilation trick of Varnavides \cite{var}, but this trick is only available in the polynomial setting if one moves to higher dimensions, see Appendix \ref{uniform-sec}.} notation
$$ \E_{n \in Y} F(n) := \frac{1}{|Y|} \sum_{n \in Y} F(n)$$
for any finite non-empty set $Y$, and
\begin{equation}\label{x-def}
 \int_X f := \E_{x \in X} f(x) = \frac{1}{N} \sum_{x \in X} f(x).
\end{equation}
We shall normalize $f$ to have mean $\int_X f = \eta_3$  and will also have the pointwise bound $0 \leq f \leq \nu$ for some ``pseudorandom measure'' $\nu$ associated to the almost primes at a sieve level $R := N^{\eta_2}$ for some\footnote{The ``missing'' values of $\eta$, such as $\eta_1$, will be described more fully in Section \ref{notation-sec}.} $0 < \eta_3 \ll \eta_2 \ll \eta_0$ (so $M$ is asymptotically larger than any fixed power of $R$).  The functions $f$, $\nu$ will be defined formally in \eqref{fdef} and \eqref{nudef} respectively, but for now let us simply remark that we will have the bound $\int_X \nu = 1 + o(1)$, together with many higher order correlation estimates on $\nu$.  

Let us defer the (sieve-theoretic) discussion of the pseudorandomness of $\nu$ for the moment, and focus instead on the (finitary) ergodic theory components of the argument.
If we were in the linear regime $M=N$ used in \cite{gt-primes} (with $N$ assumed prime for simplicity), then repeated applications of the Cauchy-Schwarz inequality (using the PET induction method) would eventually let us control the average \eqref{pn-gowers} in terms of \emph{Gowers uniformity norms} such as
$$ \|f\|_{U^d(\Z_N)} := \left( \E_{\vec m^{(0)}, \vec m^{(1)} \in [N]^d} \int_X \prod_{\omega \in \{0,1\}^d} 
T^{m^{(\omega_1)}_1 + \ldots + m^{(\omega_d)}_d} f\right)^{1/2^d}$$
for some sufficiently large $d$ (depending on $P_1,\ldots,P_k$; eventually they will be of size $O(1/\eta_1)$ for some $\eta_2 \ll \eta_1 \ll \eta_0$), where $\omega = (\omega_1,\ldots,\omega_d)$
and $\vec m^{(i)} = (m^{(i)}_1, \ldots, m^{(i)}_d)$ for $i=0,1$.  If instead we were
in the pseudo-infinitary regime $M = M(N)$ for some slowly growing function $M: \Z^+ \to \Z^+$, repeated applications of the van der Corput lemma and PET induction would allow one to control these averages by the \emph{Gowers-Host-Kra seminorms} $\|f\|_d$ from \cite{host-kra}, which in our finitary setting would be something like
$$ \|f\|_d := \left( \E_{\vec m^{(0)}, \vec m^{(1)} \in [M_1] \times \ldots \times [M_d]} \int_X \prod_{\omega \in \{0,1\}^d} 
T^{m^{(\omega_1)}_1 + \ldots + m^{(\omega_d)}_d} f\right)^{1/2^d}$$
where $M_1,\ldots,M_d$ are slowly growing functions of $N$ which we shall deliberately keep unspecified\footnote{In the traditional ergodic setting $N=\infty, M \to \infty$, one would take multiple limit superiors as $M_1,\ldots,M_d \to \infty$, choosing the order in which these parameters go to infinity carefully; see \cite{host-kra}.}.  In our intermediate setting $M = N^{\eta_0}$, however, neither of these two quantities seem to be exactly appropriate.  Instead, after applying the van der Corput lemma and PET induction one ends up considering an averaged localized Gowers norm of the form\footnote{Again, this is a slight oversimplification as we are ignoring the effects of the ``$W$-trick''.}
$$ \| f\|_{U^{\vec Q([H]^t)}_{\sqrt{M}}} := \left( \E_{\vec h \in [H]^t} 
\E_{\vec m^{(0)}, \vec m^{(1)} \in [\sqrt{M}]^d} \int_X \prod_{\omega \in \{0,1\}^d} 
T^{Q_1(\vec h) m^{(\omega_1)}_1 + \ldots + Q_d(\vec h) m^{(\omega_d)}_d} f\right)^{1/2^d}
$$
where $H = N^{\eta_7}$ is a small power of $N$ (much smaller than $M$ or $R$), $t$ is a natural number depending only on $P_1,\ldots,P_d$ (and of size $O(1/\eta_1)$), and $Q_1,\ldots,Q_d \in \Z[\h_1,\ldots,\h_t]$ are certain polynomials (of $t$ variables $\h_1,\ldots,\h_t$) which depend on $P_1,\ldots,P_d$.  Indeed, we will eventually be able (see Theorem \ref{gvn}) to establish a polynomial analogue of the \emph{generalized von Neumann theorem} in \cite{gt-primes}, which roughly speaking will assert that (if $\nu$ is sufficiently pseudorandom) any component of $f$ which is ``locally Gowers-uniform'' in the sense that the above norm is small, and which is bounded pointwise by $O(\nu+1)$, will have a negligible impact on the average \eqref{pn-gowers}.  To exploit this fact, we shall essentially repeat the arguments in \cite{gt-primes} (with some notational changes to deal with the presence of the polynomials $Q_i$ and the short shift ranges) to establish (assuming $\nu$ is sufficiently pseudorandom) an analogue of the \emph{Koopman-von Neumann-type structure theorem} in that paper, namely a decomposition $f = f_{U^\perp} + f_U$ (modulo a small error) where $f_{U^\perp}$ is bounded by $O(1)$, is non-negative and has mean roughly $\delta$, and $f_U$ is locally Gowers-uniform and thus has a negligible impact on \eqref{pn-gowers}.  Combining this with a suitable quantitative version (Theorem \ref{pSZ-quant}) of the Bergelson-Leibman theorem one can then conclude Theorem \ref{mainthm}.

We have not yet discussed how one constructs the measure $\nu$ and establishes the required pseudorandomness properties.
We shall construct $\nu$ as a truncated divisor sum at level $R = N^{\eta_2}$, although instead of using the Goldston-Y{\i}ld{\i}r{\i}m divisor sum as in \cite{goldston-yildirim-old1}, \cite{gt-primes} we shall use a smoother truncation (as in 
\cite{tao-gy-notes}
\cite{host-survey}, \cite{linprimes}) as it is slightly easier to estimate\footnote{However, in contrast to the arguments in \cite{host-survey}, \cite{linprimes} we will not be able to completely localize the estimations on the Riemann-zeta function $\zeta(s)$ to a neighborhood of the pole $s=1$, for rather minor technical reasons, and so will continue to need the classical estimates \eqref{zeta-bounds} on $\zeta(s)$ near the line $s = 1+it$.}.  The pseudorandomness conditions then reduce, after standard sieve theory manipulations, to the entirely local problem of understanding the pseudorandomness of the functions $\Lambda_p: F_p \to \R^+$ on finite fields $F_p$, defined for all primes $p$ by $\Lambda_p(x) := \frac{p}{p-1}$ when $x \neq 0\ \mod p$ and $\Lambda_p(x) = 0$ otherwise.  Our pseudorandomness conditions shall involve polynomials, and so one is soon faced with the standard arithmetic problem of counting the number of points over $F_p$ of an algebraic variety.  Fortunately, the polynomials that we shall encounter will be \emph{linear} in one or more of the variables of interest, which allows us to obtain a satisfactory count of these points without requiring deeper tools from arithmetic such as class field theory or the Weil conjectures.

\subsection{Acknowledgements}

The authors thank Brian Conrad for valuable discussions concerning algebraic varieties, Peter Sarnak for encouragement, Vitaly Bergelson and Ben Green for help with the references, and Elon Lindenstrauss, Akshay Venkatesh and Lior Silberman  for helpful conversations.  We also thank the anonymous referee and Tim Gowers for useful suggestions and corrections, and Le Thai Hoang and Julia Wolf for pointing out an error in the published version of this paper.

\section{Notation and initial preparation}\label{notation-sec}

In this section we shall fix some important notation, conventions, and assumptions which will then be used throughout
the proof of Theorem \ref{mainthm}.  Indeed, all of the sub-theorems and lemmas used to prove Theorem \ref{mainthm} will be understood to use the conventions and assumptions in this section.  We thus recommend that the reader go through this section carefully before moving on to the other sections of the paper.

Throughout this paper we fix the set $A \subset \P$ and the polynomials $P_1,\ldots,P_k \in \Z[\m]$ appearing in Theorem \ref{mainthm}.  Henceforth we shall assume that the polynomials are all distinct, since duplicate polynomials clearly have no impact on the conclusion of Theorem \ref{mainthm}.  Since we are also assuming $P_i(0) = 0$ for all $i$, we conclude that
\begin{equation}\label{pij-nonconstant}
P_i - P_{i'} \hbox{ is non-constant for all } 1 \leq i < i' \leq k.
\end{equation}

By hypothesis, the upper density
$$ \delta_0 := \limsup_{N' \to \infty}  \frac{|A \cap [N']|}{|\P \cap [N']|}$$
is strictly positive.  We shall allow all implied constants to depend on the quantities $\delta_0,P_1,\ldots,P_k$.

By the prime number theorem
\begin{equation}\label{pnt}
 |\P \cap [N']| = (1 + o(1)) \frac{N'}{\log N'}
 \end{equation}
we can find an infinite sequence of integers $N'$ going to infinity such that
\begin{equation}\label{wn-eq}
 |A \cap [N']| > \frac{1}{2} \delta_0 \frac{N'}{\log N'}.
\end{equation}
Henceforth the parameter $N'$ is always understood to obey \eqref{wn-eq}.

Throughout this paper, we let $d_*$ denote the largest degree of the polymials $P_1,\ldots,P_k$.  In addition to this quantity, we shall also need eight (!) small quantities
$$ 0 < \eta_7 \ll \eta_6 \ll \eta_5 \ll \eta_4 \ll \eta_3 \ll \eta_2 \ll \eta_1 \ll \eta_0 \ll 1$$
which depend on $\delta_0$ and on $P_1,\ldots,P_k$.  All of the assertions in this paper shall be made under the implicit assumption that $\eta_0$ is sufficiently small depending on $\delta_0$, $P_1$, $\ldots$, $P_k$; that $\eta_1$ is sufficiently small depending on $\delta_0$, $P_1$, $\ldots$, $P_k$, $\eta_0$; and so forth down to $\eta_7$, which is assumed sufficiently small depending on $\delta_0$, $P_1$, $\ldots$, $P_k$, $\eta_0$, $\ldots$, $\eta_6$ and should thus be viewed as being extremely microscopic in size.  For the convenience of the reader we briefly and informally describe the purpose of each of the $\eta_i$, their approximate size, and the importance of being that size, as follows.

\begin{itemize}

\item The parameter $\eta_0$ controls the coarse scale $M := N^{\eta_0}$.  It can be set equal to $1/2d_*$.  If one desires the quantity $m$ in Theorem \ref{mainthm} to be smaller than $x^\eps$ then one can achieve this by choosing $\eta_0$ to be less than $\eps$.
The smallness of $\eta_0$ is necessary in order to deduce Theorem \ref{mainthm} from Theorem \ref{thmquant} below.

\item  The parameter $\eta_1$ (or more precisely its reciprocal $1/\eta_1$) controls the degree of pseudorandomness needed on a certain measure $\nu$ to appear later.  Due to the highly recursive nature of the ``PET induction'' step (Section \ref{pet-sec}), it will need to be rather small; it is essentially the reciprocal of an Ackermann function of the maximum degree $d_*$ and the number of polynomials $k$.  The smallness of $\eta_1$ is needed in order to estimate all the correlations of $\nu$ which arise in the proofs of Theorem \ref{gvn} and Theorem \ref{struct-thm}.

\item The parameter $\eta_2$ controls the sieve level $R := N^{\eta_2}$.  It can be taken to be $c \eta_1 / d_*$ for some small absolute constant $c > 0$.  It needs to be small relative to $\eta_1$ in order that the inradius bound of Proposition \ref{corrprop} is satisfied.

\item The parameter $\eta_3$ measures the density of the function $f$.  It is basically of the form $c \delta_0 \eta_2$ for some small absolute constant $c > 0$.  It needs to be small relative to $\eta_2$ in order to establish the mean bound \eqref{f-mean}.

\item The parameter $\eta_4$ measures the degree of accuracy required in the Koopman-von Neumann type structure theorem
(Theorem \ref{struct-thm}).  It needs to be substantially smaller than $\eta_3$ to make the proof of Theorem \ref{pSZ-rel} in Section \ref{finale-sec} work.  The exact dependence on $\eta_3$ involves the quantitative bounds arising from the Bergelson-Leibman theorem (see Theorem \ref{pSZ-quant}).  In particular, as the only known proof of this theorem is infinitary, no explicit bounds for $\eta_4$ in terms of $\eta_3$ are currently available.

\item The parameter $\eta_5$ controls the permissible error allowed when approximating indicator functions by a smoother object, such as a polynomial; it needs to be small relative to $\eta_4$ in order to make the proof of the abstract structure theorem (Theorem \ref{abstract}) work correctly.  It can probably be taken to be roughly of the form $\exp(- C / \eta_4^C )$ for some absolute constant $C > 0$, though we do not attempt to make $\eta_5$ explicit here.

\item The parameter $\eta_6$ (or more precisely $1/\eta_6$) controls the maximum degree, dimension, and number of the polynomials that are encountered in the argument.  It needs to be small relative to $\eta_5$ in order for the polynomials arising in the proof of Proposition \ref{nu-uniform} to obey the orthogonality hypothesis \eqref{feeble} of Theorem \ref{abstract}.  It can in principle be computed in terms of $\eta_5$ by using a sufficiently quantitative version of the Weierstrass approximation theorem, though we do not do so here.

\item The parameter $\eta_7$ controls the fine scale $H := N^{\eta_7}$, which arises during the ``van der Corput'' stage of the proof in Section \ref{pet-sec}.  It needs to be small relative to $\eta_6$ in order that the ``clearing denominators'' step in the proof of Proposition \ref{nuorthog} works correctly. It is probably safe to take $\eta_7$ to be $\eta_6^{100}$ although we shall not explicitly do this.  On the other hand, $\eta_7$ cannot vanish entirely, due to the need to average out the influence of ``bad primes'' in Corollary \ref{pfc2} and Theorem \ref{pcc}.

\end{itemize}

It is crucial to the argument that the parameters are ordered in exactly the above way.  In particular, the fine scale $H = N^{\eta_7}$ needs to be much smaller than the coarse scale $M = N^{\eta_0}$.

We use the following asymptotic notation:

\begin{itemize}
\item We use $X = O(Y)$, $X \ll Y$, or $Y \gg X$ to denote the estimate $|X| \leq C Y$ for some quantity $0 < C < \infty$ which can depend on $\delta_0, P_1, \ldots, P_k$.  If we need $C$ to also depend on additional parameters we denote this by subscripts, e.g. $X = O_K(Y)$ means that $|X| \leq C_K Y$ for some $C_K$ depending on $\delta_0, P_1, \ldots, P_k, K$.
\item We use $X = o(Y)$ to denote the estimate $|X| \leq c(N') Y$, where $c$ is a quantity depending on $\delta_0,P_1,\ldots,P_k,\eta_0,\ldots, \eta_7, N'$ which goes to zero as $N' \to \infty$ for each fixed choice of
$\delta_0, P_1, \ldots, P_k, \eta_0, \ldots, \eta_7$.  If we need $c(N')$ to depend on additional parameters we denote this by subscripts, e.g. $X = o_K(Y)$ denotes the estimate $|X| \leq c_K(N') Y$, where $c_K(N')$ is a quantity which goes to zero as $N' \to \infty$ for each fixed choice of $\delta_0, P_1, \ldots, P_k, \eta_0, \ldots, \eta_7, K$.
\end{itemize}

We shall implicitly assume throughout that $N'$ is sufficiently large depending on $\delta_0, P_1, \ldots, P_k, \eta_0, \ldots, \eta_7$; in particular, all quantities of the form $o(1)$ will be small. 

Next, we perform the ``$W$-trick'' from \cite{gt-primes} to eliminate obstructions to uniformity arising from small moduli.  We shall need a slowly growing function $w = w(N')$ of $N'$.  For sake of concreteness\footnote{Actually, the arguments here work for \emph{any} choice of function $w: \Z^+ \to \Z^+$ which is bounded by $\frac{1}{10} \log \log \log N'$ and which goes to infinity as $N' \to \infty$.  This is important if one wants an explicit lower bound on the number of polynomial progressions in a certain range.} we shall set
$$ w := \frac{1}{10} \log \log \log N'.$$
We then define the quantity $W$ by
$$ W := \prod_{p < w} p$$
and the integer\footnote{Unlike previous work such as \cite{gt-primes}, we will not need to assume that $N$ is prime (which is the finitary equivalent of the underlying space $X$ being totally ergodic), although it would not be hard to ensure that this were the case if desired.  This is ultimately because we shall clear denominators as soon as they threaten to occur, and so there will be no need to perform division in $X = \Z_N$.  On the other hand, this clearing of denominators will mean that many (fine) multiplicative factors such as $Q(\vec h)$ shall attach themselves to the (coarse-scale) shifts one is averaging over.  In any case, the ``$W$-trick'' of passing from the integers $\Z$ to a residue class $W \cdot \Z + b$ can already be viewed as a kind of reduction to the totally ergodic setting, as it eliminates the effects of small periods.} $N$ by
\begin{equation}\label{ndef}
 N := \lfloor \frac{N'}{2W} \rfloor.
 \end{equation}
Here and in the sequel, all products over $p$ are understood to range over primes, and $\lfloor x \rfloor$ is the greatest integer less than or equal to $x$.  The quantity $W$ will be used to eliminate the local obstructions to pseudorandomness arising from small prime moduli; one can think of $W$ (or more precisely the cyclic group $\Z_W$) as the finitary counterpart of the ``profinite factor'' generated by the periodic functions in infinitary ergodic theory.  From the prime number theorem \eqref{pnt} one sees that
\begin{equation}\label{wsmall}
 W \ll \log \log N 
 \end{equation}
and
\begin{equation}\label{ncomp}
 N' = N^{1+o(1)}.
\end{equation}
In particular the asymptotic limit $N' \to \infty$ is equivalent to the asymptotic limit $N \to \infty$ for the purposes of the $o()$ notation, and so we shall now treat $N$ as the underlying asymptotic parameter instead of $N'$.

From \eqref{wn-eq}, \eqref{wsmall}, \eqref{ndef} we have
$$ |A \cap [\frac{1}{2} WN] \backslash [w]| \gg W \frac{N}{\log N}$$
(recall that implied constants can depend on $\delta_0$).
On the other hand, since $A$ consists entirely of primes, all the elements in $A \backslash [w]$ are coprime to $W$.
By the pigeonhole principle\footnote{In the case $A = \P$, we may use the prime number theorem in arithmetic progressions (or the Siegel-Walfisz theorem) to choose $b$, for instance to set $b=1$.  However, we will not need to exploit this ability to fix $b$ here.}, we may thus find $b = b(N) \in [W]$ coprime to $W$ such that
\begin{equation}\label{b-def}
 |\{ x \in [\frac{1}{2} N]: Wx + b \in A \}| \gg \frac{W}{\phi(W)} \frac{N}{\log N}
\end{equation}
where $\phi(W) = \prod_{p < w} (p-1)$ is the Euler totient function of $W$, i.e. the number of elements of $[W]$ which are coprime to $W$.

Let us fix this $b$.  We introduce the underlying measure space $X := \Z_N = \Z/N\Z$, with the uniform probability measure given by \eqref{x-def}.  We also introduce the \emph{coarse scale} $M := N^{\eta_0}$, the \emph{sieve level} $R := N^{\eta_2}$, and the \emph{fine scale} $H := N^{\eta_7}$.  It will be important to observe the following size hierarchy:
\begin{equation}\label{hierarchy}
1 \ll W \ll W^{1/\eta_6} \ll H \ll H^{1/\eta_6} \ll R \ll R^{1/\eta_1} \ll M \ll N = |X|.
\end{equation}
Indeed each quantity on this hierarchy is larger than any bounded power of the preceding quantity, for suitable choices of the $\eta$ parameters, for instance $R^{O(1/\eta_1)} \leq M^{1/4}$.

\begin{remark}\label{range-remark}  In the linear case \cite{gt-primes} we have $M = N$, while the parameter $H$ is not present (or can be thought of as $O(1)$).  We shall informally refer to parameters of size\footnote{Later on we shall also encounter some parameters of size $O(\sqrt{M})$ or $O(M^{1/4})$, which we shall also consider to be coarse-scale.} $O(M)$ as \emph{coarse-scale parameters}, and parameters of size $H$ as \emph{fine-scale parameters}; we shall use the symbol $m$ to denote coarse-scale parameters and $h$ for fine-scale parameters (reserving $x$ for elements of $X$).  
Note that because the sieve level $R$ is intermediate between these parameters, we will be able to easily average the pseudorandom measure $\nu$ over coarse-scale parameters, but not over fine parameters.  Fortunately, our averages will always involve at least one coarse-scale parameter, and after performing the coarse-scale averages first we will have enough control on main terms and error terms to then perform the fine averages.  The need to keep the fine parameters short arises because at one key ``Weierstrass approximation'' stage to the argument, we shall need to control the product of an extremely large number (about $O(1/\eta_6)$, in fact) of averages (or more precisely ``dual functions''), and this will cause many fine parameters to be multiplied together in order to clear denominators.  This is still tolerable because $H$ remains smaller than $R, M, N$ even after being raised to a power $O(1/\eta_6)$.  Note it is key here that the number of powers $O(1/\eta_6)$ does not depend on $\eta_7$.  It will therefore be important to keep large parts of our argument uniform in the choice $\eta_7$, although we can and will allow $\eta_7$ to influence $o(1)$ error terms.  The quantity $H$ (and thus $\eta_7$) will not actually make an impact on the argument until Section \ref{rel-sec}, when the local Gowers norms are introduced.
\end{remark}

We define the standard shift operator $T: X \to X$ on $X$ by $Tx := x+1$, with the associated action  on functions $g: X \to \R$ by $T g := g \circ T^{-1}$, thus $T^n g(x) = g(x-n)$ for any $n \in \Z$.  We introduce the normalized counting function $f: X \to \R^+$ by setting 
\begin{equation}\label{fdef}
 f(x) := \frac{\phi(W)}{W} \log R \hbox{ whenever } x \in [\frac{1}{2} N] \hbox{ and } Wx + b \in A 
\end{equation}
and $f(x) = 0$ otherwise, where we identify $[\frac{1}{2} N]$ with a subset of $\Z_N$ in the usual manner.  The
use of $\log R$ instead of $\log N$ as a normalizing factor is necessary in order to bound $f$ pointwise by the pseudorandom measure $\nu$ which we shall encounter in later sections; the ratio $\eta_2$ between $\log R$ and $\log N$ represents the relative density between the primes and the almost primes.
Observe from \eqref{b-def} that $f$ has relatively large mean:
$$
\int_X f \gg \eta_2.
$$
In particular we have
\begin{equation}\label{f-mean}
\int_X f \geq \eta_3.
\end{equation}

\begin{remark} We will eventually need to take $\eta_2$ (and hence $\eta_3$) to be quite small, in order to ensure that the measure $\nu$ obeys all the required pseudorandomness properties (this is controlled by the parameter $\eta_1$, which has not yet made a formal appearance).  Fortunately, the Bergelson-Leibman theorem (Theorem \ref{pSZ}, or more precisely Theorem \ref{pSZ-quant} below) works for sets of arbitrarily small positive density, or equivalently for (bounded) functions of arbitrarily small positive mean\footnote{As in \cite{gt-primes}, the exact quantitative bound provided by this theorem (or more precisely Theorem \ref{pSZ-quant}) will not be relevant for qualitative results such as Theorem \ref{mainthm}.  Of course, such bounds would be important if one wanted to know how soon the first polynomial progression in the primes (or a dense subset thereof) occurs; for instance such bounds influence how small $\eta_4$ and thus all subsequent $\eta$s need to be, which in turn influences the exact size of the final $o(1)$ error in Theorem \ref{thmquant}.  Unfortunately, as the only known proof of Theorem \ref{pSZ} proceeds via infinitary ergodic theory, no explicit bounds are currently known, however it is reasonable to expect (in view of results such as \cite{gowers}, \cite{tao:ergodic}) that effective bounds will eventually become available.}.  This allows us to rely on fairly crude constructions for $\nu$ which will be easier to estimate.  This is in contrast to the recent work of Goldston, Y{\i}ld{\i}r{\i}m, and Pintz \cite{pintz} on prime gaps, in which it was vitally important that the density of the prime counting function relative to the almost prime counting function be as high as possible, which in turn required a near-optimal (and thus highly delicate) construction of the almost prime counting function.
\end{remark}

To prove Theorem \ref{mainthm} it will suffice to prove the following quantitative estimate:

\begin{theorem}[Polynomial Szemer\'edi theorem in the primes, quantitative version]\label{thmquant} Let the notation and assumptions be as above.  Then we have
\begin{equation}\label{main-est}
 \E_{m \in [M]} \int_X T^{P_1(Wm)/W} f \ldots T^{P_k(Wm)/W} f \geq \frac{1}{2} c(\frac{\eta_3}{2}) - o(1)
\end{equation}
where the function $c()$ is the one appearing in Theorem \ref{pSZ-quant}.  (Observe that since $P_i(0) = 0$ for all $1 \leq i \leq k$, the polynomial $P_k(W\m)/W$ has integer coefficients.)
\end{theorem}

Indeed, suppose that the estimate \eqref{main-est} held.  Then by expanding all the averages and using \eqref{fdef} we conclude that
\begin{align*}
&|\{ (x,m) \in X \times [M]: x + P_i(Wm)/W \in [\frac{1}{2} N] \hbox{ and } W(x+P_i(Wm)/W)+b
\in A \}| \\
&\quad \geq (\frac{1}{2} c(\frac{\eta_3}{2}) - o(1)) M N \left(\frac{W}{\phi(W) \log N}\right)^k.
\end{align*}
Here we are using the fact (from \eqref{hierarchy}) that $P_i(Wm)$ is much less than $N/2$ for $m \in [M]$, and so one cannot ``wrap around'' the cyclic group $\Z_N$. 
Observe that each element in the set on the left-hand side yields a different pair $(x',m') := (Wx+b,Wm)$ with the property that $x' + P_1(m'), \ldots, x' + P_k(m') \in A$.  On the other hand, as $N \to \infty$, the right-hand side goes to infinity.  The claim follows.

\begin{remark}\label{short-rem} The above argument in fact proves slightly more than is stated by Theorem \ref{mainthm}.  Indeed, it establishes a large number of pairs $(x',m')$ with $x' + P_1(m'), \ldots, x' + P_k(m') \in A$, $x' \in [N]$, and $m' \in [M]$; more precisely, there are at least\footnote{To obtain such a bound it is important to remember that we can take $w$ and hence $W$ to be as slowly growing as one pleases; see \cite{gt-primes} for further discussion.  Note that if $A = \P$ is the full set of primes then the Bateman-Horn conjecture \cite{bh} predicts an asymptotic of the form $(\gamma + o(1)) NM / \log^k N$ for an explicitly computable $\gamma$; we do not come close to verifying this conjecture here.} $c NM / \log^k N$ such pairs for some $c$ depending on $\delta_0$, $P_1$, $\ldots$, $P_k$, $\eta_0$, $\ldots$, $\eta_7$\footnote{The arguments in this paper can be easily generalized to give  a lower bound 
of  $c NM_1\ldots M_r / \log^k N$ on the number of tuples $(x',m'_1,\ldots,m'_r)$ with $x' + P_1(m'_1,\ldots,m'_r), \ldots, x' + P_k(m'_1,\ldots,m'_r) \in A$, $x' \in [N]$,  $m'_i \in [M_i]$ $1\le i\le r$, and $P_j \in \Z[\m_1,\ldots,\m_r]$ for $1\le j\le k$. To obtain this one would only need to slightly modify the arguments in section \ref{gvn-proof-sec} (see Remark \ref{multivar-rm}), whereas the rest of the proof remains the same.}.  By throwing away the contribution of those $x'$ of size $\ll N$ (which can be done either by modifying $f$ in the obvious manner, or by using a standard upper bound sieve to estimate this component) one can in fact assume that $x'$ is comparable to $N$.  Similarly one may assume $m'$ to be comparable to $N^{\eta_0}$.  The upshot of this is that for any given $\eta_0 > 0$ one in fact obtains infinitely many ``short'' polynomial progressions $x' + P_1(m'), \ldots, x' + P_k(m')$ with $m'$ comparable to $(x')^{\eta_0}$.  One can take smaller and smaller values of $\eta_0$ and diagonalize to obtain the same statement with the bound $m' = (x')^{o(1)}$.  This stronger version of Theorem \ref{mainthm} is already new in the linear case $P_i = (i-1)\m$, although it is not too hard to modify the arguments in \cite{gt-primes} to establish it.  Note that an inspection of the Furstenberg correspondence principle reveals that the Bergelson-Liebman theorem (Theorem \ref{pSZ}) has an even stronger statement in this direction, namely that if $A$ has positive upper density in the \emph{integers} $\Z$ rather than the primes $\P$, then there exists a \emph{fixed} $m \neq 0$ for which the set  
$\{ x: x+P_i(m) \in A \hbox{ for all } 1 \leq i \leq k \}$ is infinite (in fact, it can be chosen to have 
positive upper density).  Such a statement might possibly be true for primes (or dense subsets of the primes) but is well beyond the technology of this paper.  For instance, to establish such a statement even in the simple case $P_1=0$, $P_2 = \m$ is tantamount to asserting that the primes have bounded gaps arbitrarily often, which is still not known unconditionally even after the recent breakthroughs in \cite{pintz}.  On the other hand it may be possible to establish such a result with a logarithmic dependence between $x'$ and $m'$, e.g. $m' \ll \log^{O(1)} x'$.  We do not pursue this issue here.
\end{remark} 

\begin{figure}[tb]
\centerline{\psfig{figure=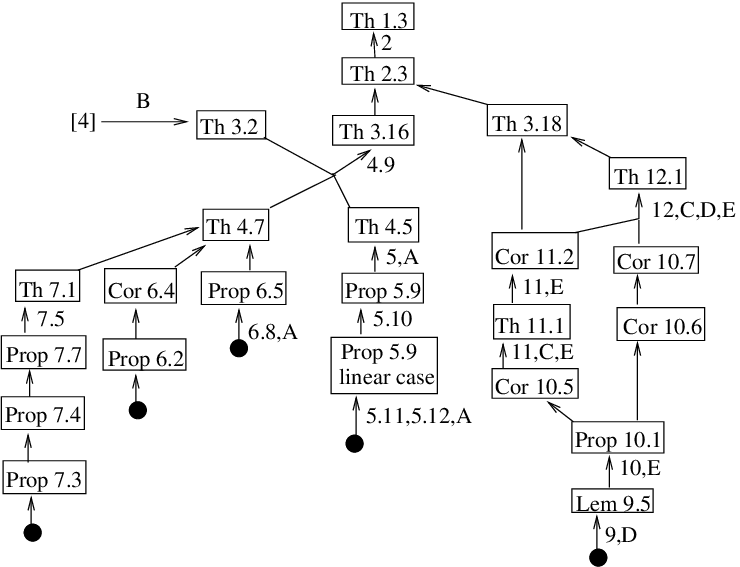}}
\caption{The main theorems in this paper and their logical dependencies.  The numbers and letters next to the arrows indicate the section(s) where the implication is proven, and which appendices are used; if no section is indicated, the result is proven immediately after it is stated.  Self-contained arguments are indicated using a filled-in circle.}
\label{fig1}
\end{figure}

It remains to prove Theorem \ref{thmquant}.  This shall occupy the remainder of the paper.  The proof is lengthy, but splits into many non-interacting parts; see Figure \ref{fig1} for a diagram of the logical dependencies of this paper.

\subsection{Miscellaneous notation}

To conclude this section we record some additional notation which will be used heavily throughout this paper.

We have already used the notation $\Z[\m]$ to denote the ring of integer-coefficient polynomials in one indeterminate\footnote{We shall use boldface letters to denote abstract indeterminates, reserving the non-boldface letters for concrete realizations of these indeterminates, which in this paper will always be in the ring of integers $\Z$.} $\m$.  More generally one can consider $\Z[\x_1,\ldots,\x_d]$, the ring of integer-coefficient polynomials in $d$ indeterminates $\x_1,\ldots,\x_d$.  More generally still we have $\Z[\x_1,\ldots,\x_d]^D$, the space of $D$-tuples of polynomials in $\Z[\x_1,\ldots,\x_d]$; note that each element of this space defines a polynomial map from $\Z^d$ to $\Z^D$.  Thus we shall think of elements of $\Z[\x_1,\ldots,\x_d]^D$ as $D$-dimensional-valued integer-coefficient polynomials over $d$ variables.  The \emph{degree} of a monomial $\x_1^{n_1} \ldots \x_d^{n_d}$ is $n_1+ \ldots + n_d$; the degree of a polynomial in $\Z[\x_1,\ldots,\x_d]^D$ is the highest degree of any monomial which appears in any component of the polynomial; we adopt the convention that the zero polynomial has degree $-\infty$.  We say that two $D$-dimensional-valued polynomials $\vec P, \vec Q \in \Z[\x_1,\ldots,\x_d]^D$ are \emph{parallel} if we have $n \vec P = m \vec Q$ for some non-zero integers $n,m$.

If $\vec n = (n_1,\ldots,n_D)$ and $\vec m = (m_1,\ldots,m_D)$ are two vectors in $\Z^d$, we use
$\vec n \cdot \vec m := n_1 m_1 + \ldots n_D m_D \in \Z$ to denote their dot product.

If $f: X \to \R$ and $g: X \to \R$ are two functions, we say that $f$ is \emph{pointwise bounded by} $g$, and write $f \leq g$, if we have $f(x) \leq g(x)$ for all $x \in X$.  Similarly, if $g: X \to \R^+$ is non-negative, we write $f = O(g)$ if we have $f(x) = O(g(x))$ uniformly for all $x \in X$.
If $A \subset X$, we use $1_A: X \to \{0,1\}$ to denote the indicator function of $A$; thus $1_A(x) = 1$ when $x \in A$ and $1_A(x) = 0$ when $x \not \in A$.  Given any statement $P$, we use $1_P$ to denote $1$ when $P$ is true and $0$ when $P$ is false, thus for instance $1_A(x) = 1_{x \in A}$.

We define a \emph{convex body} to be an open bounded convex subset of a Euclidean space $\R^d$.  We define the \emph{inradius} of a convex body to be the radius of the largest ball that is contained inside the body; this will be a convenient measure of how ``large'' a body is\footnote{In our paper there will only be essentially two types of convex bodies: ``coarse-scale'' convex bodies of inradius at least $M^{1/4}$, and ``fine-scale'' convex bodies, of inradius at least $\gg H$.  In almost all cases, the convex bodies will in fact simply be rectangular boxes.}.

\section{Three pillars of the proof}

As in \cite{gt-primes}, our proof of Theorem \ref{thmquant} rests on three independent pillars - a quantitative Szemer\'edi-type theorem (proven by traditional ergodic theory), a transference principle (proven by finitary ergodic theory), and the construction of a pseudorandom majorant $\nu$ for $f$ (with the pseudorandomness proven by sieve theory).  In this section we describe each these pillars separately, and state where they are proven.  

\subsection{The quantitative Szemer\'edi-type theorem}

Theorem \ref{thmquant} concerns a multiple polynomial average of an unbounded function $f$.  To control such an object, we first need to establish an estimate for \emph{bounded} functions $g$.  This is achieved as follows (cf. \cite[Proposition 2.3]{gt-primes}):

\begin{theorem}[Polynomial Szemer\'edi theorem, quantitative version]\label{pSZ-quant} Let the notation and assumptions be as in the previous section.  Let $\delta > 0$, and let $g: X \to \R$ be any function obeying the pointwise bound $0 \leq g \leq 1 + o(1)$ and the mean bound
$\int_X g \geq \delta - o(1)$.  Then we have
\begin{equation}\label{main-est-g}
 \E_{m \in [M]} \int_X T^{P_1(Wm)/W} g \ldots T^{P_k(Wm)/W} g \geq c(\delta) - o(1)
\end{equation}
for some $c(\delta) > 0$ depending on $\delta,P_1,\ldots,P_k$, but independent of $N$ or $W$.  
\end{theorem}

It is not hard to see that this theorem implies Theorem \ref{pSZ}.  The converse is not immediately obvious (the key point being, of course, that the bound $c(\delta)$ in \eqref{main-est-g} is uniform in both $N$ and $W$); however, it is not hard to deduce Theorem \ref{pSZ-quant} from (a multidimensional version of) Theorem \ref{pSZ} and the Furstenberg correspondence principle; one can also use the uniform version of the Bergelson-Leibman theorem proved in \cite{bmht}. As the arguments here are fairly standard, and are unrelated to those in the remainder of the paper, we defer the proof of Theorem \ref{pSZ-quant} to Appendix \ref{uniform-sec}.

\subsection{Pseudorandom measures}

To describe the other two pillars of the argument it is necessary for the measure $\nu$ to make its appearance.  (The precise properties of $\nu$, however, will not actually be used until Sections \ref{gvn-proof-sec} and \ref{polydual-sec}.)

\begin{definition}[Measure]\cite{gt-primes} 
A \emph{measure} is a non-negative function $\nu: X \to \R^+$ with the total mass estimate
\begin{equation}\label{measure} \int_X \nu = 1 + o(1)
\end{equation}
and the crude pointwise bound
\begin{equation}\label{nu-point}
\nu \ll_\eps N^\eps
\end{equation}
for any $\eps > 0$.
\end{definition}

\begin{remark} As remarked in \cite{gt-primes}, it is really $\nu \mu_X$ which is a measure rather than $\nu$, where $\mu_X$ is the uniform probability measure on $X$; $\nu$ should be more accurately referred to as a ``probability density'' or ``weight function''.  However we retain the terminology ``measure'' for compatibility with \cite{gt-primes}.  
The condition \eqref{nu-point} is needed here to discard certain error terms arising from the boundary effects of shift ranges (such as those arising from the van der Corput lemma).  This condition does not prominently feature in \cite{gt-primes}, as the shifts range over all of $\Z_N$, which has no boundary.  Fortunately, \eqref{nu-point} is very easy to establish for the majorant that we shall end up using.  We note though that while the right-hand side of \eqref{nu-point} does not look too large, we cannot possibly afford to allow factors such as $N^\eps$ to multiply into error terms such as $o(1)$ as these terms will almost certainly cease to be small at that point.  Hence we can only really use \eqref{nu-point} in situations where we already have a polynomial gain in $N$, which can for instance arise by exploiting the gaps in \eqref{hierarchy}.
\end{remark}

The simplest example of a measure is the constant measure $\nu \equiv 1$.  Another model example worth keeping in mind is the random measure where $\nu(x) = \log R$ with independent probability $1/\log R$ for each $x \in X$, and $\nu(x) = 0$ otherwise.
The following definitions attempt to capture certain aspects of this random measure, which will eventually be satisfied by a certain truncated divisor sum concentrated on almost primes.  These definitions are rather technical, and their precise form is only needed in later sections of the paper.  They are somewhat artificial in nature, being a compromise between the type of control needed to establish the relative polynomial Szemer\'edi theorem (Theorem \ref{pSZ-rel}) and the type of control that can be easily verified for truncated divisor sums (Theorem \ref{major}).  It may well be that a simpler notion of pseudorandomness can be given.

\begin{definition}[Polynomial forms condition]\label{polyform-def} Let $\nu: X \to \R^+$ be a measure.  We say that $\nu$ obeys the \emph{polynomial forms condition} if, given any $0 \leq J, d \leq 1/\eta_1$, any polynomials $Q_1,\ldots,Q_J \in \Z[\m_1,\ldots,\m_d]$ of $d$ unknowns of total degree at most $d_*$ and coefficients at most $W^{1/\eta_1}$, with $Q_j - Q_{j'}$ non-constant for every distinct $j,j' \in [J]$, for every $\eps > 0$, and for every convex body $\Omega \subset \R^d$ of inradius at least $N^\eps$, and contained in the ball $B(0,M^2)$, we have the bound
\begin{equation}\label{polyform}
 \E_{\vec h \in \Omega \cap \Z^d} \int_X \prod_{j \in [J]} T^{Q_j(\vec h)} \nu = 1 + o_{\eps}(1).
\end{equation}
\end{definition}

Note the first appearance of the parameter $\eta_1$, which is controlling the degree of the pseudorandomness here.  Note also that the bound is uniform in the coefficients of the polynomials $Q_1,\ldots,Q_J$.

\begin{examples}  The mean bound \eqref{measure} is a special case of \eqref{polyform}; another simple example is
$$ \E_{h \in [H]} \int_X \nu T^h \nu T^{Wh^2} \nu = 1 + o(1).$$
Observe that the smaller one makes $\eta_1$, the stronger the polynomial forms condition becomes.
\end{examples}

\begin{remark} Definition \ref{polyform-def} is a partial analogue of the ``linear forms condition'' in \cite{gt-primes}.  The parameter $\eta_1$ is playing multiple roles, controlling the degree, dimension, number and size of the polynomials in question.  It would be more natural to split this parameter into four parameters to control each of these attributes separately, but we have chosen to artificially unify these four parameters in order to simplify the notation slightly.  The parameter $\eps$ will eventually be set to be essentially $\eta_7$, but we leave it arbitrary here to emphasize that the definition of pseudorandomness does not depend on the choice of $\eta_7$ (or $H$).  This will be important later, basically because we need to select $\nu$ (or more precisely $\eta_2$ (or $R$), which is involved in the construction of $\nu$) before we are allowed to choose $\eta_7$.
\end{remark}

The next condition is in a similar spirit, but considerably more complicated; it allows for arbitrarily many factors in the average, as long as they have a partly linear structure, and that they are organized into relatively small groups, with a separate coarse-scale averaging applied to each of the groups.

\begin{definition}[Polynomial correlation condition]\label{polycor-def} Let $\nu: X \to \R^+$ be a measure.  We say that $\nu$ obeys the \emph{polynomial correlation condition} if, given any $0 \leq D, J, L \leq 1/\eta_1$, any integers $D', D'', K > 0$, and any $\eps > 0$, and given any vector-valued polynomials
\begin{align*}
\vec P_{j} &\in \Z[ \h_1, \ldots, \h_{D''} ]^{D}\\
\vec Q_{j,k} &\in  \Z[ \h_1, \ldots, \h_{D''} ]^{D'}\\
\vec S_l &\in \Z[ \h_1, \ldots, \h_{D''} ]^{D'}
\end{align*}
of degree at most $1/\eta_1$ for $j \in [J]$, $k \in [K]$, $l \in [L]$ obeying the non-degeneracy conditions
\begin{itemize}
	\item For any distinct $j,j' \in [J]$ and any $k \in [K]$, the $D+D'$-dimensional-valued polynomials $(\vec P_j, \vec Q_{j,k})$ and $(\vec P_{j'}, \vec Q_{j',k})$ are not parallel.
	\item The coefficients of $\vec P_j$ and $\vec S_l$ are bounded in magnitude by $W^{1/\eta_1}$.
	\item The $D'$-dimensional-valued polynomials $\vec S_l$ are distinct as $l$ varies in $[L]$.
\end{itemize}
and given any convex body $\Omega \subset \R^D$ of inradius at least $M^{1/4}$ and convex bodies $\Omega' \subset \R^{D'}$, $\Omega'' \subset \R^{D''}$ of inradius at least $N^\eps$, with all convex bodies contained in $B(0,M^2)$, then
\begin{equation}\label{polycor-eq}
\begin{split}
 &\E_{\vec n \in \Omega' \cap \Z^{D'}} \E_{\vec h \in \Omega'' \cap \Z^{D''}} \int_X  \\
&\left[ \prod_{k \in [K]} \E_{\vec m \in \Omega \cap \Z^D} \prod_{j \in [J]} T^{\vec P_j(\vec h) \cdot \vec m + \vec Q_{j,k}(\vec h) \cdot \vec n} \nu \right ] \prod_{l \in [L]} T^{\vec S_l(\vec h) \cdot \vec n} \nu \\
&\quad = 1 + o_{D',D'',K,\eps}(1)
\end{split}
\end{equation}
\end{definition}

\begin{remark} It will be essential here that $D', D'', K$ can be arbitrarily large\footnote{An analogous phenomenon occurs in the correlation condition in \cite{gt-primes}, where it was essential that the exponent $q$ appearing in that condition (which is roughly analogous to $K$ here) could be arbitrarily large.}; otherwise, this condition becomes essentially a special case of the polynomial forms condition.  Indeed in our argument, these quantities will get as large as $O(1/\eta_6)$, which is far larger than $1/\eta_1$.  As in the preceding definition, $\eps$ will eventually be set to equal essentially $\eta_7$, but we refrain from doing so here to keep the definition of pseudorandomness independent of $\eta_7$, to avoid the appearance of circularity in the argument.
\end{remark}

\begin{remark} The correlation condition \eqref{polycor-eq} would follow from the polynomial forms condition \eqref{polyform} if we had the pointwise bounds
\begin{equation}\label{aas} \E_{\vec m \in \Omega \cap \Z^D} \prod_{j \in [J]} T^{\vec P_j(\vec h) \cdot \vec m + \vec Q_{j,k}(\vec h) \cdot \vec n} \nu = 1 + o(1)
\end{equation}
for each $k \in [K]$ and all $\vec h, \vec n$.  Unfortunately, such a bound is too optimistic to be true: for instance, if $\vec P_j(\vec h) = Q_{j,k}(\vec h) = 0$ then the left-hand side is an average of $\nu^J$, which is almost certainly much larger than $1$.  In the number-theoretic applications in which $\nu$ is supposed to concentrate on almost primes, one also has similar problems when
$\vec P_j(\vec h), Q_{j,k}(\vec h)$ are non-zero but very smooth (i.e. they have many small prime factors slightly larger than $w$).  In \cite{gt-primes} these smooth cases were modeled by a weight function $\tau$, which obeyed arbitrarily large moment conditions which led to integral estimates analogous to \eqref{polycor-eq}.  In this paper we have found it more convenient to not explicitly create the weight function, instead placing the integral estimate \eqref{polycor-eq} in the correlation condition hypothesis directly.  In fact one can view \eqref{polycor-eq} as an assertion that \eqref{aas} holds ``asymptotically almost everywhere'' (cf. Proposition \ref{essbound} below).
\end{remark}

\begin{remark} One could generalize \eqref{polycor-eq} slightly by allowing the number of terms
$J$ in the $j$ product to depend on $k$, but we will not need this strengthening and in any event it follows automatically from \eqref{polycor-eq} by a H\"older inequality argument similar to that used
in Lemma \ref{avglem} below.
\end{remark}

\begin{definition}[Pseudorandom measure]  A \emph{pseudorandom measure} is any measure $\nu$ which obeys both the polynomial forms condition and the correlation condition.
\end{definition}

The following lemma (cf. \cite[Lemma 3.4]{gt-primes}) is useful:

\begin{lemma}\label{avglem} If $\nu$ is a pseudorandom measure, then so is $\nu_{1/2} := (1 + \nu)/2$ (possibly with slightly different decay rates for the $o(1)$ error terms).
\end{lemma}

\begin{proof}  It is clear that $\nu_{1/2}$ satisfies \eqref{measure} and \eqref{nu-point}.  Because
$\nu$ obeys the polynomial forms condition \eqref{polyform}, one can easily verify using the binomial
formula that $\nu_{1/2}$ does also.  Now we turn to the polynomial correlation condition, which requires a little more care.  Setting $\vec Q_{j,k}$ to be independent of $k$, we obtain that
\begin{align*}
 &\E_{\vec n \in \Omega' \cap \Z^{D'}} \E_{\vec h \in \Omega'' \cap \Z^{D''}} \int_X  \\
&\left[ \E_{\vec m \in \Omega \cap \Z^D} \prod_{j \in [J]} T^{\vec P_j(\vec h) \cdot \vec m + \vec Q_{j}(\vec h) \cdot \vec n} \nu \right ]^K \prod_{l \in [L]} T^{\vec S_l(\vec h) \cdot \vec n} \nu \\
&\quad = 1 + o_{D',D'',K,\eps}(1)
\end{align*}
for all $K \geq 0$ and $\vec P_j$, $\vec Q_j$, $\vec S_l$ obeying the hypotheses of the correlation condition.  By the binomial formula this implies that
\begin{align*}
 &\E_{\vec n \in \Omega' \cap \Z^{D'}} \E_{\vec h \in \Omega'' \cap \Z^{D''}} \int_X  \\
&\left[ \E_{\vec m \in \Omega \cap \Z^D} \prod_{j \in [J]} T^{\vec P_j(\vec h) \cdot \vec m + \vec Q_{j}(\vec h) \cdot \vec n} \nu - 1 \right ]^K \prod_{l \in [L]} 
T^{\vec S_l(\vec h) \cdot \vec n} \nu \\
&\quad = 0^K + o_{D',D'',K,\eps}(1).
\end{align*}
(Recall of course that $0^0=1$.)
Let us take $K$ to be a large \emph{even} integer.
Another application of the binomial formula allows one to replace the final $\nu$ by $\nu_{1/2}$.  
By the triangle inequality in a weighted Lebesgue norm $l^K$, we may then replace the other occurrences $\nu$
by $\nu_{1/2}$ also:
\begin{align*}
 &\E_{\vec n \in \Omega' \cap \Z^{D'}} \E_{\vec h \in \Omega'' \cap \Z^{D''}} \int_X  \\
&\left[ \E_{\vec m \in \Omega \cap \Z^D} \prod_{j \in [J]} T^{\vec P_j(\vec h) \cdot \vec m + \vec Q_{j}(\vec h) \cdot \vec n} \nu_{1/2} - 1 \right ]^K \prod_{l \in [L]} 
T^{\vec S_l(\vec h) \cdot \vec n} \nu_{1/2} \\
&\quad = 0^K + o_{D',D'',K,\eps}(1).
\end{align*}
This was only proven for even $K$, but follows also for odd $K$ by the Cauchy-Schwarz inequality \eqref{cauchy-schwarz}.  By H\"older's inequality we obtain a similar statement when the $Q_j$ are now allowed to vary in $k$:
\begin{align*}
 &\E_{\vec n \in \Omega' \cap \Z^{D'}} \E_{\vec h \in \Omega'' \cap \Z^{D''}} \int_X  \\
&\left[ \prod_{k \in [K]} (\E_{\vec m \in \Omega \cap \Z^D} \prod_{j \in [J]} T^{\vec P_j(\vec h) \cdot \vec m + \vec Q_{j,k}(\vec h) \cdot \vec n} \nu_{1/2}  - 1 ) \right ] \prod_{l \in [L]} 
T^{\vec S_l(\vec h) \cdot \vec n} \nu_{1/2} \\
&\quad = 0^K + o_{D',D'',K,\eps}(1).
\end{align*}
Applying the binomial formula again we see that $\nu_{1/2}$ obeys \eqref{polycor-eq} as desired.
\end{proof}

\subsection{Transference principle}

We can now state the second pillar of our argument (cf. \cite[Theorem 3.5]{gt-primes}).

\begin{theorem}[Relative polynomial Szemer\'edi theorem]\label{pSZ-rel}
Let the notation and assumptions be as in Section \ref{notation-sec}.  Then given any pseudorandom measure $\nu$ and any $g: X \to \R$ obeying the pointwise bound $0 \leq g \leq \nu$ and the mean bound 
\begin{equation}\label{mean}
\int_X g \geq \eta_3,
\end{equation}
we have  
\begin{equation}\label{main-est-g-rel}
 \E_{m \in [M]} \int_X T^{P_1(Wm)/W} g \ldots T^{P_k(Wm)/W} g \geq \frac{1}{2} c(\frac{\eta_3}{2}) - o(1)
\end{equation}
where $c()$ is the function appearing in Theorem \ref{pSZ-quant}.
\end{theorem}

Apart from inessential factors of $2$ (and the substantially worse decay rates concealed within the $o(1)$ notation), this theorem is significantly stronger than Theorem \ref{pSZ-quant}, which is essentially the special case $\nu=1$.  In fact we shall derive Theorem \ref{pSZ-rel} from Theorem \ref{pSZ-quant} using the transference principle technology from \cite{gt-primes}.  The argument is lengthy and will occupy Sections \ref{rel-sec}-\ref{struct-sec}.  

\subsection{Construction of the majorant}

To conclude Theorem \ref{thmquant} from Theorem \ref{pSZ-rel} and \eqref{f-mean} it clearly suffices to show (cf. \cite[Proposition 9.1]{gt-primes})

\begin{theorem}[Existence of pseudorandom majorant]\label{major}  
Let the notation and assumptions be as in Section \ref{notation-sec}.  Then there exists a pseudorandom measure $\nu$ such that the function $f$ defined in \eqref{fdef} enjoys the the pointwise bound $0 \leq f \leq \nu$.
\end{theorem}

This is the third pillar of the argument.  The majorant $\nu$ acts as an ``enveloping sieve'' for the primes (or more precisely, for the primes equal to $b$ modulo $W$), in the sense of \cite{Ramare}, \cite{Ramare-Ruzsa}.  It is defined explicitly in Section \ref{sec8}.
However, for the purposes of the proof of the other pillars of the argument
(Theorem \ref{pSZ-quant} and Theorem \ref{pSZ-rel}) it will not be necessary to know the precise definition of $\nu$, only that $\nu$ majorizes $f$ and is pseudorandom.  In order to establish this
pseudorandomness it is necessary that $\eta_2$ is small compared to $\eta_1$.  On the other hand,
observe that $\nu$ does not depend on $H$ and thus is insensitive to the choice of $\eta_7$.

The proof of Theorem \ref{major} follows similar lines to those in \cite{gt-primes}, \cite{linprimes}, except that the ``local'' or ``singular series'' calculation is more complicated, as one is now forced to count solutions to one or more \emph{polynomial} equations over $F_p$ rather than linear equations.  Fortunately, it turns out that the polynomials involved happen to be \emph{linear} in at least one ``coarse-scale'' variable, and so the number of solutions can be counted relatively easily, without recourse to any deep arithmetic facts (such as the Weil conjectures).  We establish Theorem \ref{major} in Sections \ref{sec8}-\ref{pcc-sec},
using some basic facts about convex bodies, solutions to polynomial equations in $F_p$, and distribution of prime numbers which are recalled in Appendices \S \ref{convex-sec}-\ref{prime-sec}
respectively.

\section{Overview of proof of transference principle}\label{rel-sec}

We now begin the proof of the relative polynomial Szemer\'edi theorem (Theorem \ref{pSZ-rel}).  As in \cite{gt-primes}, this theorem will follow quickly from three simpler components.  The first is the uniformly quantitative version of the ordinary polynomial Szemer\'edi theorem, Theorem \ref{pSZ-quant}, which will be proven in Appendix \ref{uniform-sec}.  The second is a ``polynomial generalized von Neumann theorem'' (Theorem \ref{gvn}) which allows us to 
neglect the contribution of sufficiently ``locally Gowers-uniform'' contributions to \eqref{main-est-g-rel}.  The third is a ``local Koopman-von Neumann structure theorem'' (Theorem \ref{struct-thm}) which decomposes a function $0 \leq f \leq \nu$ (outside of a negligible set) into a bounded positive component $f_{U^\perp}$ and a locally Gowers-uniform error $f_U$.  The purpose
of this section is to formally state the latter two components and show how they imply Theorem \ref{pSZ-rel}; the proofs of these components will then occupy subsequent sections of the paper.

The pseudorandom measure $\nu$ plays no role in the ordinary polynomial Szemer\'edi theorem, Theorem \ref{pSZ-quant}.  In the von Neumann theorem, Theorem \ref{gvn}, the pseudorandomness of $\nu$ is exploited via the polynomial forms condition (Definition \ref{polyform-def}).  In the structure theorem, Theorem \ref{struct-thm}, it is instead the polynomial correlation condition (Definition \ref{polycor-def}) which delivers the benefits of pseudorandomness.

\subsection{Local Gowers norms}

As mentioned in the introduction, a key ingredient in the proof of Theorem \ref{pSZ-rel} will be the introduction of a norm $\| \|_{U^{\vec Q([H]^t,W)}_M}$ which controls averages such as those in \eqref{main-est-g-rel}.  It is here that the parameter $\eta_7$ first makes an appearance, via the shift range $H$.  The purpose of this subsection is to define these norms formally.

Let $f: X \to \R$ be a function.  For any $d \geq 1$, recall that the (global) \emph{Gowers uniformity norm} $\|f\|_{U^d}$ of $f$ is defined by the formula
$$ \|f\|_{U^d}^{2^d} := \E_{m_1,\ldots,m_d \in \Z_N} \int_X \prod_{(\omega_1,\ldots,\omega_d) \in \{0,1\}^d}
T^{\omega_1 m_1 + \ldots + \omega_d m_d} f.$$
An equivalent definition is
$$ \|f\|_{U^d}^{2^d} := \E_{m^{(0)}_1,\ldots,m^{(0)}_d,m^{(1)}_1,\ldots,m^{(1)}_d \in \Z_N} 
\int_X \prod_{(\omega_1,\ldots,\omega_d) \in \{0,1\}^d}
T^{m^{(\omega_1)}_1 + \ldots + m^{(\omega_d)}_d} f$$
as can be seen by making the substitutions $m^{(1)}_i := m^{(0)}_i + m_i$ and shifting the integral by $m^{(0)}_1 + \ldots + m^{(0)}_d$.

We will not directly use the global Gowers norms in this paper, because the range of the shifts $m$ in those norms is too large for our applications.  Instead, we shall need local versions of this norm.  For any steps $a_1,\ldots,a_d \in \Z$, we define the \emph{local Gowers uniformity norm} $U^{a_1,\ldots,a_d}_{\sqrt{M}}$ by\footnote{We will need to pass from shifts of size $O(M)$ to shifts of size $O(\sqrt{M})$ to avoid dealing with certain boundary terms (similar to those arising in the van der Corput lemma).}
\begin{equation}\label{local-gowers}
 \|f\|_{U^{a_1,\ldots,a_d}_{\sqrt{M}}}^{2^d} := \E_{m^{(0)}_1,\ldots,m^{(0)}_d,m^{(1)}_1,\ldots,m^{(1)}_d \in [\sqrt{M}]} 
\int_X \prod_{(\omega_1,\ldots,\omega_d) \in \{0,1\}^d}
T^{m^{(\omega_1)}_1 a_1 + \ldots + m^{(\omega_d)}_d a_d} f.
\end{equation}
Thus for instance, when $\sqrt{M}=N$ and $a_1,\ldots,a_d$ are invertible in $\Z_N^\times$, then the $U^{a_1,\ldots,a_d}_{\sqrt{M}}$ norm is the same as the $U^d$ norm.  When $\sqrt{M}$ is much smaller than $N$, however, there appears to be no obvious comparison between these two norms.  It is not immediately obvious that the local Gowers norm is indeed a norm, but we shall show this in Appendix \ref{gowers-sec} where basic properties of these norms are established.  In practice we shall take $a_1,\ldots,a_d$ to be rather small compared to $R$ or $M$, indeed these steps will have size $O(H^{O(1)})$.

\begin{remark} One can generalize this norm to complex valued functions $f$ by conjugating those factors of $f$ for which $\omega_1 + \ldots + \omega_d$ is odd.  If we then set $f = e(\phi) = e^{2\pi i \phi}$ for some phase function $\phi: X \to \R/\Z$, then the local Gowers
$\|f\|_{U^{a_1,\ldots,a_d}_{\sqrt{M}}}$ norm is informally measuring the extent to which the $d$-fold difference
$$ \sum_{\omega_1, \ldots, \omega_d \in \{0,1\}} (-1)^{\omega_1 + \ldots + \omega_d} \phi( x + m^{(\omega_1)}_1 a_1 + \ldots + m^{(\omega_d)}_d a_d )$$
is close to zero, where $x$ ranges over $X$ and $m^{(0)}_i, m^{(1)}_i$ range over $[M]$ for $i \in [d]$.  Even more informally, these norms are measuring the extent to which $\phi$ ``behaves like'' a polynomial of degree less than $d$ on arithmetic progressions of the form
$$ \{ x + m_1 a_1 + \ldots + m_d a_d: m_1, \ldots, m_d \in [\sqrt{M}] \}$$
where $x \in X$ is arbitrary.  The global Gowers norm $U^d$, in contrast, measures similar behavior over the entire space $X$.
\end{remark}

We shall estimate the Gowers-uniform contributions to \eqref{main-est-g-rel} via repeated application of the van der Corput lemma using the standard \emph{polynomial exhaustion theorem} (PET) induction scheme.  This will eventually allow us to control these contributions, not by a single local Gowers-uniform norm, but rather by an \emph{average} of such norms, in which the shifts $h_1,\ldots,h_d$ are fine and parameterized by a certain polynomial.  More precisely, given
any $t \geq 0$ and any $d$-tuple $\vec Q = (Q_1,\ldots,Q_d) \in \Z[\h_1,\ldots,\h_t,\W]^d$ of polynomials, we define the \emph{averaged local Gowers uniformity norm} $U^{\vec Q([H]^t,W)}_{\sqrt{M}}$ by the formula
\begin{equation}\label{avgdef}
 \| f \|_{U^{\vec Q([H]^t,W)}_{\sqrt{M}}}^{2^d} := \E_{\vec h \in [H]^t} 
\|f\|_{U^{Q_1(\vec h,W),\ldots,Q_d(\vec h,W)}_{\sqrt{M}}}^{2^d}.
\end{equation}
Inserting \eqref{local-gowers} we thus have
\begin{equation}\label{avgdef-2}
\begin{split}
 \| f \|_{U^{\vec Q([H]^t,W)}_{\sqrt{M}}}^{2^d} := &\E_{\vec h \in [H]^t} 
\E_{m^{(0)}_1,\ldots,m^{(0)}_d,m^{(1)}_1,\ldots,m^{(1)}_d \in [\sqrt{M}]} 
\int_X \\
&\quad \prod_{(\omega_1,\ldots,\omega_d) \in \{0,1\}^d}
T^{m^{(\omega_1)}_1 Q_1(\vec h, W) + \ldots + m^{(\omega_d)}_d Q_d(\vec h, W)} f.
\end{split}
\end{equation}
In Appendix \ref{gowers-sec} we show that the local Gowers uniformity norms are indeed norms; by the triangle inequality in $l^{2^d}$, this implies that the averaged local Gowers uniformity norms are also norms.  To avoid degeneracies we will assume that none of the polynomials $Q_1,\ldots,Q_d$ vanish.

\begin{remark} The distinction between local Gowers uniform norms and their averaged counterparts is a necessary feature of our ``quantitative'' setting.  In the ``qualitative'' setting of traditional (infinitary) ergodic theory (where $X$ is infinite), there is no need for this sort of distinction; if the local Gowers uniformity norms go to zero as $M \to \infty$ for the shifts $h_1=\ldots=h_d = 1$, then it is not hard (using various forms of the Cauchy-Schwarz-Gowers inequality, such as those in Appendix \ref{gowers-sec}) to show the same is true for any other fixed choice of shifts $h_1,\ldots,h_d$, and hence the averaged norms will also go to zero as $M \to \infty$ for any fixed choice of $\vec Q$ and $H$.  The converse implications are also easy to establish.  Thus one can use a single local Gowers uniformity norm, $U^{1,\ldots,1}_{\sqrt{M}}$, to control everything in the limit $M \to \infty$ with $H$ bounded; this then corresponds to the Gowers-Host-Kra seminorms used in \cite{host-kra}, \cite{leibman} to control polynomial averages.  However in our more quantitative setting, where $H$ is allowed to grow like a (very small) power of $N$, we cannot afford to use the above equivalences (as they will amplify the $o(1)$ errors in our arguments to be unacceptably large), and so must turn instead to the more complicated-seeming averaged local Gowers uniformity norms.  
\end{remark}

\subsection{The polynomial generalized von Neumann theorem}

We are now ready to state the second main component of the proof of Theorem \ref{pSZ-rel} (the first component being Theorem \ref{pSZ-quant}).

\begin{theorem}[Polynomial generalized von Neumann theorem]\label{gvn}
Let the notation and assumptions be as in Section \ref{notation-sec}.  Then there exists a $d \geq 2$, $t \geq 0$ of size $O(1)$ and $d$-tuple $\vec Q \in \Z[\h_1,\ldots,\h_t,\W]^d$ of degree at most $d_*$ and coefficients $O(1)$, with none of the components of $\vec Q$ vanishing, as well as a constant $c > 0$ depending only on $P_1,\ldots,P_k$,
such that one has the inequality
\begin{equation}\label{gvn-est}
|\E_{m \in [M]} \int_X T^{P_1(Wm)/W} g_1 \ldots T^{P_k(Wm)/W} g_k| \ll \min_{1 \leq i \leq k} \| g_i \|_{U^{\vec Q([H]^t,W)}_{\sqrt{M}}}^c + o(1)
\end{equation}
for any functions $g_1,\ldots,g_k: X \to \R$ obeying the pointwise bound $|g_i| \leq 1 + \nu$ for all $1 \leq i \leq k$ and $x \in X$, and some pseudorandom measure $\nu$.
\end{theorem}

This theorem is a local polynomial analogue of \cite[Proposition 5.3]{gt-primes}.  It will be proven by a vast number
of applications of the van der Corput lemma and the Cauchy-Schwarz inequality following the standard PET induction scheme; the idea is to first apply the van der Corput lemma repeatedly to linearize the polynomials $P_1,\ldots,P_k$, and then apply Cauchy-Schwarz repeatedly to estimate the linearized averages by local Gowers norms.  The presence of the measure $\nu$ will cause a large number of shifts of $\nu$ to appear as weights, but these will ultimately be controllable via the polynomial forms condition (Definition \ref{polyform-def}).
The final values of $d$ and $t$ obtained will be very large (indeed, they exhibit Ackermann-type behavior in the maximal degree of $P_1,\ldots,P_k$) but can be chosen to be small compared to $1/\eta_1$, which controls the pseudorandomness of $\nu$.

The proof of Theorem \ref{gvn} is elementary but rather lengthy (and notation-intensive), and shall occupy all of Section \ref{gvn-proof-sec}.
The $\nu=1$ case of this theorem is a finitary version of a similar result in \cite{host-kra}, while the linear case of this theorem (when the $P_i-P_j$ are all linear) is essentially in \cite{gt-primes}.  Indeed the proof of this theorem will use a combination of the techniques from both of these papers.

\subsection{The local Koopman-von Neumann theorem}

The third major component of the proof of Theorem \ref{pSZ-rel} is the following structure theorem.

\begin{theorem}[Structure theorem]\label{struct-thm}
Let the notation and assumptions be as in Section \ref{notation-sec}.  Let $t \geq 0, d \geq 2$ be of size $O(1)$, and let $\vec Q \in \Z[\h_1,\ldots,\h_t,\W]^d$ be polynomials of degree at most $d_*$ and coefficients $O(1)$ (with none of the components of $\vec Q$ vanishing).  Then given any pseudorandom measure $\nu$ and any $g: X \to \R^+$ with the pointwise bound $0 \leq g \leq \nu$, there exist functions $g_{U^\perp}, g_U: X \to \R$ with the pointwise bound
\begin{equation}\label{gdecomp}
0 \leq g_{U^\perp}(x) + g_U(x) \leq g(x)
\end{equation}
of $g$ obeying the following estimates:
\begin{itemize}
  \item(Boundedness of structured component) We have the pointwise bound
  \begin{equation}\label{struct-bound}
    0 \leq g_{U^\perp}(x) \leq 1.
  \end{equation}
	\item($g_{U^\perp}$ captures most of the mass)  We have
	\begin{equation}\label{most-mass}
	  \int_X g_{U^\perp} \geq \int_X g - O(\eta_4) - o(1).
  \end{equation} 
  \item(Uniformity of unstructured component) We have the bound
  \begin{equation}\label{unstruct-bound}
   \| g_U \|_{U^{\vec Q([H]^t,W)}_{\sqrt{M}}} \leq \eta_4^{1/2^d} + o(1).
  \end{equation}
\end{itemize}
\end{theorem}

\begin{remark} Note the first apperance of the parameter $\eta_4$, which is controlling the accuracy of this structure theorem.  One can make this accuracy as strong as desired, but at the cost of pushing $\eta_7$ (and thus $H$) down, which will ultimately worsen many of the $o(1)$ errors appearing here and elsewhere.
\end{remark}

Theorem \ref{struct-thm} is the most technical and difficult component of the entire paper, and is proven in Sections \ref{polydual-sec}-\ref{struct-sec}.  It is a ``finitary ergodic theory'' argument which relies on iterating a certain ``dichotomy between structure and randomness''.  Here, the randomness is measured using the local Gowers uniformity norm $U^{\vec Q([H]^t,W)}_{\sqrt{M}}$.  To measure the structured component we need the machinery of \emph{dual functions}, as in \cite{gt-primes}, together with an energy incrementation argument which we formalize abstractly in Theorem \ref{abstract}.  A key point will be that $\nu-1$ is ``orthogonal'' to these dual functions in a rather strong sense (see Proposition \ref{nuorthog}), which will be the key to approximating functions bounded by $\nu$ with functions bounded by $1$.  This will be accomplished by a rather tricky series of applications of the Cauchy-Schwarz inequality and will rely heavily on the polynomial correlation condition (Definition \ref{polycor-def}).

\subsection{Proof of Theorem \ref{pSZ-rel}}\label{finale-sec}

Using Theorem \ref{gvn} and Theorem \ref{struct-thm} we can now quickly prove Theorem \ref{pSZ-rel} (and hence Theorem \ref{mainthm}, assuming Theorem \ref{major}), following the same argument as in \cite{gt-primes}.  

Let the notation and assumptions be as in Section \ref{notation-sec}.  Let $\nu$ be a pseudorandom measure, and let $g: X \to \R$ obey the pointwise bound $0 \leq g \leq \nu$ and \eqref{mean}.

Let $d > 0$, $t > 0$ and $\vec Q$ be as in Theorem \ref{gvn}; these expressions depend only on $P_1,\ldots,P_k$ and so we do not need to explicitly track their influence on the $O()$ and $o()$ notation.  Applying Theorem \ref{struct-thm}, we thus obtain functions $g_U$, $g_{U^\perp}$ obeying the properties claimed in that theorem.  From \eqref{gdecomp} we have
\begin{align*}
& \E_{m \in [M]} \int_X T^{P_1(Wm)/W} g \ldots T^{P_k(Wm)/W} g \geq \\
&\quad \E_{m \in [M]} \int_X T^{P_1(Wm)/W} (g_{U^\perp}+g_U) \ldots T^{P_k(Wm)/W} (g_{U^\perp}+g_U).
\end{align*}
We expand the right-hand side into $2^k = O(1)$ terms.  Consider any of the $2^k-1$ of these terms which involves at least
one factor of $g_U$.  From \eqref{gdecomp}, \eqref{struct-bound} we know that $g_U$ and $g_{U^\perp}$ are both
bounded pointwise in magnitude by $\nu + 1 + o(1)$, which is $O(\nu + 1)$ when $N$ is large enough.  Thus by Theorem \ref{gvn} and \eqref{unstruct-bound}, the contribution of all of these terms can be bounded in magnitude by 
$$ \ll \| g_U \|_{U^{\vec Q([H]^t,W)}_{\sqrt{M}}}^c + o(1) \ll \eta_4^{c/2^d} + o(1)$$
for some $c > 0$ depending only on $P_1,\ldots,P_k$.
On the other hand, from \eqref{most-mass}, \eqref{mean} and the choice of parameters we have
$$ \int_X g_{U^\perp} \geq \frac{1}{2} \eta_3.$$
Applying this, \eqref{struct-bound}, and Theorem \ref{pSZ-quant} obtain
$$
 \E_{m \in [M]} \int_X T^{P_1(Wm)/W} g_{U^\perp} \ldots T^{P_k(Wm)/W} g_{U^\perp} \geq c(\eta_3/2) > 0.$$
Putting all this together we conclude
$$
 \E_{m \in [M]} \int_X T^{P_1(Wm)/W} g \ldots T^{P_k(Wm)/W} g \geq c(\eta_3/2) - O(\eta_4^{c/2^d}) - o(1).$$
As $\eta_4$ is chosen small compared to $\eta_3$, Theorem \ref{pSZ-rel} follows.

\endprf

\section{Proof of generalized von Neumann theorem}\label{gvn-proof-sec}

In this section we prove Theorem \ref{gvn}.  In a nutshell, our argument here will be a rigorous implementation of the following scheme:
\begin{align*}
\hbox{polynomial } &\hbox{average} \ll \hbox{weighted linear average} + o(1) \quad\quad \hbox{(van der Corput)} \\
&\ll \hbox{weighted parallelopiped average} + o(1) \quad\quad \hbox{(weighted gen. von Neumann)}\\
&\ll \hbox{unweighted parallelopiped average} + o(1). \quad\quad \hbox{(Cauchy-Schwarz)}
\end{align*}
The argument is based upon that used to prove \cite[Proposition 5.3]{gt-primes}, namely repeated application of the Cauchy-Schwarz inequality to replace various functions $g_i$ by $\nu$ (or $\nu+1$), followed by application of the polynomial forms condition (Definition \ref{polyform-def}) to replace the resulting polynomial averages of $\nu$ with $1+o(1)$.  The major new ingredient in the argument compared to \cite{gt-primes} will be the \emph{polynomial exhaustion theorem} (PET) induction scheme (used for instance in \cite{bl}) in order to estimate the polynomial average in \eqref{gvn-est} by a linear average similar to that treated in \cite[Proposition 5.3]{gt-primes}.  After using PET induction to achieve this linearization, the rest of the proof is broadly similar to that in \cite[Proposition 5.3]{gt-primes}, except for the fact that the shift parameters are restricted to be of size $M$ or $\sqrt{M}$ rather than $N$, and that there is also some additional averaging over short shift parameters of size $O(H)$.

The arguments are elementary and straightforward, but will require a rather large amount of new notation in order to keep track of all the weights and factors created by applications of the Cauchy-Schwarz inequality.  Fortunately, none of this notation will be needed in any other section; indeed, this section can be read independently of the rest of the paper (although it of course relies on the material in earlier sections, and also on Appendix \ref{gowers-sec}).

We begin with some simple reductions.  First observe (as in \cite{gt-primes}) that Lemma \ref{avglem} allows us to replace the hypotheses $|g_i| \leq 1 + \nu$ by the slightly stronger $|g_i| \leq \nu$, at the (acceptable) cost of worsening the implied constant in \eqref{gvn-est} by a factor of $2^k$.
Next, we claim that it suffices to find $d$, $t$, $c$ and $\vec Q \in \Z[\h_1,\ldots,\h_t,\W]^d$ for which we have the weaker estimate
\begin{equation}\label{polyav}
|\E_{m \in [M]} \int_X T^{P_1(Wm)/W} g_1 \ldots T^{P_k(Wm)/W} g_k| \ll \| g_1 \|_{U^{\vec Q([H]^t,W)}_{\sqrt{M}}}^c + o(1)
\end{equation}
(i.e. we only control the average using the norm of $g_1$, rather than the best norm of all the $g_i$).  Indeed, if we could show this, then by symmetry we could find $d_i$, $t_i$, and $\vec Q_i \in \Z[h_1,\ldots,\h_{t_i},\W]^{d_i}$ for $i = 1,\ldots,k$ such that
$$
|\E_{m \in [M]} \int_X T^{P_1(Wm)/W} g_1 \ldots T^{P_k(Wm)/W} g_k| \ll \| g_i \|_{U^{\vec Q_i([H]^{t_i},W)}_{\sqrt{M}}}^{c_i} + o(1)$$
whenever $1 \leq i \leq k$ and $\nu$ is pseudorandom.  The claim then follows by using Lemma
\ref{uk-concatenate} to obtain a local Gowers norm $U^{\vec Q([H]^t,W)}_{\sqrt{M}}$ which dominates each of
the individual norms $U^{\vec Q_i([H]^{t_i},W)}_{\sqrt{M}}$, and taking $c := \min_{1 \leq i \leq k} c_i$. (Note that the pointwise bound $|g_i| \leq \nu$ and the polynomial forms condition easily imply that the $U^{\vec Q_i([H]^{t_i},W)}_{\sqrt{M}}$ norm of $g_i$ is $O(1)$.)

It remains to prove \eqref{polyav}.  It should come as no surprise to the experts that this type of ``generalized von Neumann'' theorem will be proven via a large number of applications of van der Corput's lemma and the Cauchy-Schwarz inequality.  In order to keep track of the intermediate
multilinear expressions which arise during this process, it is convenient to prove a substantial generalization of this estimate.  We first need the notion of a polynomial system, and averages associated to such systems.

\begin{definition}[Polynomial system]\label{polysys}  A \emph{polynomial system} $\Sys$ consists of the following objects:
\begin{itemize}
\item An integer $D \geq 0$, which we call the \emph{number of fine degrees of freedom};
\item A non-empty finite index set $A$ (the elements of which we shall refer to as \emph{nodes} of the system);
\item A polynomial $R_\alpha \in \Z[\m, \h_1, \ldots, \h_D,\W]$ in $D+2$ variables and of degree at most $d_*$ attached to each node $\alpha \in A$;
\item An \emph{distinguished node} $\alpha_0 \in A$;
\item A (possibly empty) collection $A' \subset A \backslash \{\alpha_0\}$ of \emph{inactive nodes}.  The nodes in $A \backslash A'$ will be referred to as \emph{active}, thus for instance the distinguished node $\alpha_0$ is always active.
\end{itemize}
We say that a node $\alpha \in A$ is \emph{linear} if $R_\alpha - R_{\alpha_0}$ is at most linear in $\m$, thus the distinguished node is always linear.  We say that the entire system $\Sys$ is linear if every active node is linear.
We make the following non-degeneracy assumptions:
\begin{itemize}
\item If $\alpha, \beta$ are distinct nodes in $A$, then $R_\alpha - R_\beta$ is not constant in $\m,\h_1,\ldots,\h_D$. 
\item If $\alpha, \beta$ are distinct \emph{linear} nodes in $A$, then $R_\alpha - R_\beta$ is not constant in $\m$.
\end{itemize}
Given any two nodes $\alpha,\beta$, we define the \emph{distance} $d(\alpha,\beta)$ between the two nodes to be the $\m$-degree of the polynomial $R_\alpha - R_\beta$ (which is non-zero by hypothesis); thus this distance is symmetric, non-negative, and obeys the non-archimedean triangle inequality 
$$ d(\alpha,\gamma) \leq \max( d(\alpha,\beta), d(\beta,\gamma) ).$$
Note that $\alpha$ is linear if and only if $d(\alpha,\alpha_1) \leq 1$, and furthermore we have $d(\alpha,\beta)=1$ for any two distinct linear nodes $\alpha,\beta$.
\end{definition}

\begin{example}\label{simple-tree} Take $D := 0$, $A := \{1,2,3\}$, with $R_1 := 0$, $R_2 := \m$, and $R_3 := \m^2$ with distinguished node $3$.  Then the $3$ node is linear and the other two are non-linear.  (If the distinguished node was $1$ or $2$, the situation would be reversed.)
\end{example}

\begin{remark} The non-archimedean semi-metric is naturally identifiable with a tree whose terminal nodes are the nodes of $\Sys$, and whose intermediate nodes are balls with respect to this semi-metric; the distance between two nodes is then the height of their join.  It is this tree structure (and the distinction of nodes into active, inactive, and distinguished nodes) which shall implicitly govern the dynamics of the PET induction scheme which we shall shortly perform.  
We will however omit the details, as we shall not explicitly use this tree structure in this paper.
\end{remark}

\begin{definition}[Realizations and averages]  Let $\Sys$ be a polynomial system, and $\nu$ be a measure.  We define a \emph{$\nu$-realization} $\vec f = (f_\alpha)_{\alpha \in A}$ of $\Sys$ to be an assignment of a function $f_\alpha: X \to \R$ to each node $\alpha$ with the following properties:
\begin{itemize}
\item For any node $\alpha$, we have the pointwise bound $|f_\alpha| \leq \nu$.
\item For any inactive node $\alpha$, we have $f_\alpha = \nu$.
\end{itemize}
We refer to the function $f_{\alpha_0}$ attached to the distinguished node $\alpha_0$ as the \emph{distinguished function}.
We define the \emph{average} $\Lambda_\Sys(\vec f) \in \R$ of a system $\Sys$ and its $\nu$-realization $\vec f$ to be the quantity
$$ \Lambda_\Sys(\vec f) :=
\E_{h_1,\ldots,h_D \in [H]} \E_{m \in [M]} \int_X \prod_{\alpha \in A} T^{R_\alpha(m,h_1,\ldots,h_D,W)} f_\alpha.$$
\end{definition}

\begin{example} If $\Sys$ is the system in Example \ref{simple-tree}, then
\begin{equation}\label{lsys}
 \Lambda_\Sys(\vec f) = \E_{m \in [M]} \int_X f_1 T^m f_2 T^{m^2} f_3.
 \end{equation}
\end{example}

\begin{example}\label{desys} The average
$$ \E_{h, h' \in [H]} \E_{m \in [M]} \int_X \nu T^{m+h} f_2 T^{m+h'} f_2 T^{(m+h)^2} f_3 
T^{(m+h')^2} f_3$$
can be written in the form $\Lambda_\Sys(\vec f)$ with distinguished function $f_3$, where $\Sys$ is a system with $D := 2$, $A := \{1,2,2',3,3'\}$ with the $1$ node inactive and distinguished node $3$, with $R_1 := 0$, $R_2 := \m + \h_1$, $R_{2'} := \m + \h_2$, $R_3 := (\m + \h_1)^2$, $R_4 := (\m + \h_2)^2$, and $\vec f$ is given by
$f_1 := \nu$, $f_{2'} := f_2$, and $f_{3'} := f_3$. 
\end{example}

\begin{example}[Base example]\label{base-ex} Let $\Sys$ be the system with $D:=0$, $A := \{1,\ldots,k\}$, $\alpha_0 := 1$, $A' = \emptyset$ (thus all nodes are active), 
and $Q_i := P_i( \W \m ) / \W$.  We observe from \eqref{pij-nonconstant} that this is indeed a system.  Then $\vec f := (g_1,\ldots,g_k)$ is a $\nu$-realization of $\Sys$ with distinguished function $g_1$, and
$$ \Lambda_\Sys(\vec f) = \E_{m \in [M]} \int_X T^{P_1(Wm)/W} g_1 \ldots T^{P_k(Wm)/W} g_k.$$
This system $\Sys$ is linear if and only if the polynomials $P_i-P_j$ are all linear.
\end{example}

\begin{remark}[Translation invariance]\label{trans-remark}  Given a polynomial system $\Sys$ and a polynomial
$R \in \Z[\m,\h_1,\ldots,\h_t,\W]$, we can define the shifted polynomial system $\Sys-R$ by replacing each of the polynomials $R_\alpha$ by $R_\alpha-R$; it is easy to verify that this does not affect any of the characteristics of the system, and in particular we have $\Lambda_\Sys(\vec f) = \Lambda_{\Sys-R}(\vec f)$ for any $\nu$-realization $\vec f$ of $\Sys$ (and hence of $\Sys-R$).  This translation invariance gives us the freedom to set any single polynomial $R_\alpha$ of our choosing to equal $0$; indeed we shall exploit this freedom whenever we wish to use van der Corput's lemma or Cauchy-Schwarz to deactivate any given node.  
\end{remark}

The estimate \eqref{polyav} then follows immediately from

\begin{proposition}[Generalized von Neumann theorem for polynomial systems]\label{gvn-induct}  Let $\Sys$ be a polynomial system with distinguished node $\alpha_0$.  Then, if $\eta_1$ is sufficiently small depending on $\Sys, \alpha_0$, there exists $d, t \geq 0$, $\vec Q \in \Z[\h_1,\ldots,\h_t,\W]^d$ and $c > 0$ depending only on $\Sys, \alpha_0$, with $\vec Q$ having degree at most $d_*$, such that one has the bound
$$ |\Lambda_\Sys(\vec f)| \ll_\Sys \| f_{\alpha_0} \|_{U^{\vec Q([H]^t, W)}_{\sqrt{M}}}^c + o_\Sys(1)$$
whenever $\nu$ is a pseudorandom measure and $\vec f$ is a $\nu$-realization of $\Sys$ with distinguished function $f_{\alpha_0}$.
\end{proposition}

Indeed, one simply applies this proposition to Example \ref{base-ex} to conclude \eqref{polyav}.

It remains to prove Proposition \ref{gvn-induct}.  This will be done in three stages.  The first is the ``linearization'' stage, in which a weighted form of van der Corput's lemma and the polynomial forms condition are applied repeatedly (using the PET induction scheme) to reduce the proof of Proposition \ref{gvn-induct} to the case where the system $\Sys$ is linear.  The second stage is the ``parallelopipedization'' stage, in which one uses a weighted variant of the ``Cauchy-Schwarz-Gowers inequality'' to estimate the average $\Lambda_\Sys(f)$ associated to a linear system $\Sys$ by a weighted average of the distinguished function $f_{\alpha_0}$ over parallelopipeds.  Finally, there is a relatively simple ``Cauchy-Schwarz'' stage in which the polynomial forms condition is used one last time to replace the weights by the constant $1$, at which point the proof of Proposition \ref{gvn-induct} is complete.  We remark that the latter two stages (dealing with the linear system case) also appeared in \cite{gt-primes}; the new feature here is the initial linearization step, which can be viewed as a weighted variant of
the usual polynomial generalized von Neumann theorem (see e.g. \cite{bl}).  This linearization step is
also the step which shall use the fine shifts $h = O(H)$ (for reasons which will be clearer in Section \ref{polydual-sec}); this should be contrasted with the parallelopipedization step, which relies on coarse-scale shifts $m = O(\sqrt{M})$.  

\subsection{PET induction and linearization}\label{pet-sec}

We now reduce Proposition \ref{gvn-induct} to the linear case.  We shall rely heavily here on van der Corput's lemma, which in practical terms allows us to deactivate any given node at the expense of duplicating all the other nodes in the system.  Since this operation tends to increase the number of active nodes in the system, it is not immediately obvious that iterating this operation will eventually simplify the system.  To make this more clear we need to introduce the notion of a weight vector.

\begin{definition}[Weight vector]  A \emph{weight vector} is an infinite vector $\vec w = (w_1, w_2, \ldots, )$ of non-negative integers $w_i$, with only finitely many of the $w_i$ being non-zero.  Given two weight vectors $\vec w = (w_1, w_2, \ldots)$ and $\vec w' = (w'_1, w'_2, \ldots)$, we say that $\vec w < \vec w'$ if there exists $k \geq 1$ such that $w_k < w'_k$, and such that $w_i = w'_i$ for all $i > k$.  We say that a weight vector is \emph{linear} if $w_i = 0$ for all $i \geq 2$.
\end{definition}

It is well known that the space of all weight vectors forms a well-ordered set; indeed, it is isomorphic to the ordinal $\omega^\omega$.  In particular we may perform strong induction on this space.  The space of linear weight vectors forms an order ideal; indeed, a weight is linear if and only if it is less than $(0,1,0,\ldots)$.

\begin{remark} In this paper, we will only need to work with weight vectors $\vec w = (w_1, w_2, \ldots, )$ with $w_i$ vanishing for $i>d_*$, so the ordinal structure can in fact be reduced from $\omega^\omega$ to $\omega^{d_*}$ if desired.
\end{remark}

\begin{definition}[Weight]\label{weight}  Let $\Sys$ be a polynomial system, and let $\alpha$ be a node in $\Sys$ (in practice this will not be the distinguished node $\alpha_0$).  We say that two nodes $\beta,\gamma$ in $\Sys$ are \emph{equivalent} relative to $\alpha$ if $d(\beta,\gamma) < d(\alpha,\beta)$.  This is an equivalence relation on the nodes of $\Sys$, and the equivalence classes have a well-defined distance to $\alpha$.  We define the \emph{weight vector} $\vec w_{\alpha}(\Sys)$ of $\Sys$ relative to $\alpha$ by setting the $i^{\operatorname{th}}$ component for any $i \geq 1$ to equal the number of equivalence classes at distance $i$ from $\alpha$.
\end{definition}

\begin{example}\label{simplex} Consider the system in Example \ref{simple-tree}.  The weight of this system relative to the $1$ node is $(1,1,0,\ldots)$, whereas the weight of the system in Example \ref{desys} relative to the $2$ node is $(0,1,0,\ldots)$ (note that the inactive node $1$ is not relevant here, nor is the node $2'$ which has distance $0$ from $2$), which is a lower weight than that of the previous system.
\end{example}

The key inductive step in the reduction to the linear case is then

\begin{proposition}[Inductive step of linearization]\label{linearize}  Let $\Sys$ be a polynomial system with distinguished node $\alpha_0$ and a non-linear active node $\alpha$.  If $\eta_1$ is sufficiently small depending on $\Sys, \alpha_0, \alpha$, then there exists a polynomial system $\Sys'$ with distinguished node $\alpha'_0$ and an active node $\alpha'$ with $\vec w_{\alpha'}(\Sys') < \vec w_\alpha(\Sys)$ with the following property: given any pseudorandom measure $\nu$ and any $\nu$-realization $\vec f$ of $\Sys$, there exists a $\nu$-realization $\vec f'$ of $\Sys'$ with the same distinguished function (thus $f_{\alpha_0} = f'_{\alpha'_0}$) such that
\begin{equation}\label{lsysv}
 |\Lambda_\Sys(\vec f)|^2 \ll \Lambda_{\Sys'}(\vec f') + o_\Sys(1).
 \end{equation}
\end{proposition}

Indeed, given this proposition, a strong induction on the weight vector $\vec w_\alpha(\Sys)$ immediately implies that in order to prove Proposition \ref{gvn-induct}, it suffices to do so for linear systems (since these, by definition, are the only systems without non-linear active nodes).

Before we prove this proposition in general, it is instructive to give an example.

\begin{example} Consider the expression \eqref{lsys} with $f_1, f_2, f_3$ bounded pointwise by $\nu$.  We rewrite this expression as
$$ \Lambda_\Sys(\vec f) = \int_X f_1 \E_{m \in [M]} T^m f_2 T^{m^2} f_3$$
and thus by Cauchy-Schwarz \eqref{cauchy-schwarz}
$$ |\Lambda_\Sys(\vec f)|^2 \leq (\int_X \nu) (\int_X \nu |\E_{m \in [M]} T^m f_2 T^{m^2} f_3|^2 ).$$
By \eqref{measure} the first factor is $1+o(1)$.  Also, from van der Corput's lemma (Lemma \ref{VDC}) we have
$$ |\E_{m \in [M]} T^m f_2 T^{m^2} f_3|^2 \leq \E_{h,h' \in [H]} \E_{m \in [M]} T^{m+h} f_2 T^{m+h'} f_2 T^{(m+h)^2} f_3 T^{(m+h')^2} f_3 + o(1).$$
We may thus conclude a bound of the form \eqref{lsysv}, where $\Lambda_{\Sys'}(\vec f')$ is the quantity studied in Example \ref{desys}.  Note from Example \ref{simplex} that $\Sys'$ has a lower weight than $\Sys$ relative to suitably chosen nodes.
\end{example}

\begin{proof}[Proof of Proposition \ref{linearize}]  By using translation invariance (Remark \ref{trans-remark}) we may normalize $R_\alpha = 0$.  We split $A = A_0 \cup A_1$, where $A_0 := \{ \beta \in A: d(\alpha,\beta)=0\}$ and $A_1 := A \backslash A_0$. Since $\alpha$ is nonlinear, the distinguished node $\alpha_0$ lies in $A_1$.  Then $R_\beta$ is independent of $\m$ for all $\beta \in A_0$: $R_\beta(m,h_1,\ldots,h_d,W) = R_\beta(h_1,\ldots,h_d,W)$.  We can then write
$$ \Lambda_\Sys(\vec f) = \E_{h_1,\ldots,h_D \in [H]} \int_X F_{h_1,\ldots,h_D} \E_{m \in [M]} G_{m,h_1,\ldots,h_D}$$
where
$$ F_{h_1,\ldots,h_D} := \prod_{\beta \in A_0} T^{R_\beta(h_1,\ldots,h_d,W)} f_\beta$$
and
$$ G_{m,h_1,\ldots,h_D} := T^{R_\beta(m,h_1,\ldots,h_d,W)} f_\beta.$$
Since $|f_\beta|$ is bounded pointwise by $\nu$, we have $|F_{h_1,\ldots,h_D}| \leq H_{h_1,\ldots,h_D}$
where
\begin{equation}\label{hdelta}
 H_{h_1,\ldots,h_D} := \prod_{\beta \in A_0} T^{R_\beta(h_1,\ldots,h_D,W)} \nu
 \end{equation}
and thus by Cauchy-Schwarz
\begin{align*}
|\Lambda_\Sys(\vec f)|^2 &\leq (\E_{h_1,\ldots,h_D \in [H]} \int_X H_{h_1,\ldots,h_D})\\
& \quad\quad \times \E_{h_1,\ldots,h_D \in [H]} \int_X H_{h_1,\ldots,h_D} |\E_{m \in M} G_{m,h_1,\ldots,h_D}|^2.
\end{align*}
Since $\nu$ is pseudorandom and thus obeys the polynomial forms condition, we see from Definition \ref{polyform-def} and \eqref{hdelta} (taking $\eta_1$ sufficiently small) that
$$ \E_{h_1,\ldots,h_D \in [H]} \int_X H_{h_1,\ldots,h_D} = 1 + o_\Sys(1)$$
(note by hypothesis that the $R_\beta - R_{\beta'}$ are not constant in $\m, \h_1,\ldots,\h_D$).    Since we are always assuming $N$ to be large, the $o_\Sys(1)$ error is bounded.
Thus we reduce to showing that
$$ \E_{h_1,\ldots,h_D \in [H]} \int_X H_{h_1,\ldots,h_D} |\E_{m \in M} G_{m,h_1,\ldots,h_D}|^2
\ll \Lambda_{\Sys'}(\vec f') + o_\Sys(1)$$
for some suitable $\Sys'$, $\vec f'$.  But by the van der Corput lemma (Lemma \ref{VDC})
(using \eqref{nu-point} to get the upper bounds on $G_{m,h_1,\ldots,h_D}$), we have
$$ |\E_{m \in M} G_{m,h_1,\ldots,h_D}|^2 \ll \E_{h,h' \in [H]} \E_{m \in [M]} 
G_{m+h,h_1,\ldots,h_D} G_{m+h',h_1,\ldots,h_D} + o(1)$$
and so to finish the proof it suffices to verify that the expression
\begin{equation}\label{ehh}
 \E_{h_1,\ldots,h_D,h,h' \in [H]} \E_{m \in [M]} \int_X H_{h_1,\ldots,h_D}
G_{m+h,h_1,\ldots,h_D} G_{m+h',h_1,\ldots,h_D}
\end{equation}
is of the form $\Lambda_{\Sys'}(\vec f')$ for some suitable $\Sys'$ and $\vec f'$.  By inspection we see that we can construct $\Sys'$ and $\vec f'$ as follows:

\begin{itemize}
\item We have $D+2$ fine degrees of freedom, which we label $\h_1,\ldots,\h_D,\h,\h'$;
\item The nodes $A'$ of $\Sys'$ are $A' := A_0 \cup A_1 \cup A'_1$, where $A'_1$ is another copy of $A_1$ (disjoint from $A_0 \cup A_1$), with distinguished node $\alpha'_0 = \alpha_0 \in A_1$.
\item We choose the node $\alpha'$ to be an active node of $\Sys$ in $A_1$ which has minimal distance to $\alpha_0$.  (Note that $A_1$ always contains at least one active node, namely $\alpha_0$.)
\item If $\beta \in A_0$, then $\beta$ is inactive in $\Sys'$, with $R'_\beta(\m,\h_1,\ldots,\h_d,\h,\h',\W)$ $:=$ $R_\beta(\m,\h_1,\ldots,\h_d,\W)$ and $\vec f'_\beta := \nu$;
\item If $\beta \in A_1$, then $\beta$ is inactive in $\Sys'$ if and only if it is inactive in $\Sys$, with $R'_{\beta}(\m,\h_1,\ldots,\h_d,\h,\h',\W)$ $:=$ $R_{\beta}(\m+\h, \h_1,\ldots,\h_d,\W)$, and $\vec f'_\beta := \vec f_\beta$;
\item If $\beta' \in A'_1$ is the counterpart of some node $\beta \in A_1$, then $\beta'$ is inactive in $\Sys'$ if and only if $\beta$ is inactive in $\Sys$,
with $R'_{\beta'}(\m,\h_1,\ldots,\h_d,\h,\h',\W)$ $:=$ $R_{\beta}( \m+\h',\h_1,\ldots,\h_d,\W)$, and $\vec f'_{\beta'} := \vec f_\beta$.
\end{itemize}

It is then straightforward to verify that $\Sys'$ is a polynomial system, that $\vec f'$ is a realization of $\Sys'$, and that \eqref{ehh} is equal to $\Lambda_{\Sys'}(\vec f')$.  
It remains to show that $\vec w_{\alpha'}(\Sys') < \vec w_\alpha(\Sys)$.  Let $d$ be the degree in $\m$ of
$R_{\alpha} - R_{\alpha'}$, thus $d \geq 1$.  One easily verifies that the $i^{th}$ component of
$\vec w_{\alpha'}(\Sys')$ will be equal to that of $\vec w_\alpha(\Sys)$ for $i > d$, and equal to one less than that of $\vec w_\alpha(\Sys)$ when $i=d$ (basically due to the deactivation of all the nodes at $A_0$).  The claim follows.  (The behavior of these weight vectors for $i < d$ is much more complicated, but is fortunately not relevant due to our choice of ordering on weight vectors.)
\end{proof}

\subsection{Parallelopipedization}\label{par-sec}

By the preceding discussion, we see that to prove Proposition \ref{gvn-induct} it suffices to do so in the case where $\Sys$ is linear.  To motivate the argument let us first work through an unweighted example (with $\nu=1$).

\begin{example}[Unweighted linear case] Consider the linear average
$$ \Lambda_\Sys(\vec f) = \E_{h, h' \in [H]} \E_{m \in [M]} \int_X f_0 T^{hm} f_1 T^{h'm} f_2$$
with distinguished function $f_0$, and with $|f_0|$, $|f_1|$, $|f_2|$ bounded pointwise by $1$.  We introduce some new coarse-scale shift parameters $m_1, m_2 \in [\sqrt{M}]$.  By shifting $m$ to $m - m_1 - m_2$ one can express
the above average as
$$ \E_{h, h' \in [H]} \E_{m \in [M]} \int_X \E_{m_1,m_2 \in [\sqrt{M}]} 
f_0 T^{h(m-m_1-m_2)} f_1 T^{h'(m-m_1-m_2)} f_2 + o(1)$$
and then shifting the integral by $T^{hm_1 + h' m_2}$ we obtain
$$ \E_{h, h' \in [H]} \E_{m \in [M]} \int_X \E_{m_1,m_2 \in [\sqrt{M}]} 
T^{hm_1+h'm_2} f_0 T^{h(m-m_2)+h'm_2} f_1 T^{h'(m-m_1)+hm_1} f_2 + o(1).$$
The point is that the factor $T^{h(m-m_2)+h'm_2} f_1$ does not depend on $m_1$, while
$T^{h'(m-m_1)+hm_1} f_2$ does not depend on $m_2$.  One can then use the Cauchy-Schwarz-Gowers inequality
(see e.g. \cite[Corollary B.3]{linprimes}) to estimate this expression by
$$ (
\E_{h, h' \in [H]} \E_{m \in [M]} \int_X \E_{m_1,m'_1,m_2,m'_2 \in [\sqrt{M}]} 
T^{hm_1+h'm_2} f_0 T^{hm'_1+h'm_2} f_0 T^{hm_1+h'm'_2} f_0 T^{hm'_1+h'm'_2} f_0 )^{1/4} + o(1).$$
The main term here can then be recognized as a local Gowers norm \eqref{avgdef-2}.
\end{example}

Now we return to the general linear case.  Here we will need to address the presence of many additional weights which are all shifted versions of the measure $\nu$, which requires the repeated use of weighted Cauchy-Schwarz  inequalities.  See \cite[\S 5]{gt-primes} for a worked example of this type of computation.  Our
arguments here shall instead follow those of \cite[\S C]{linprimes}, in particular relying on the weighted generalized von Neumann inequality from that paper (reproduced here as Proposition \ref{wgn}).

We turn to the details.
To simplify the notation we write $\vec h := (h_1,\ldots,h_d)$ and $\vec \h := (\h_1,\ldots,\h_d)$.
We use the translation invariance (Remark \ref{trans-remark}) to normalize $R_{\alpha_0} = 0$.
We then split $A = \{\alpha_0\} \cup A_l \cup A_{nl}$, where $A_l$ consists of all the linear nodes, and $A_{nl}$ all the nonlinear (and hence inactive) nodes.  By the non-degeneracy assumptions in
Definition \ref{polysys}, we may write
$$ R_\alpha = b_\alpha \m + c_\alpha$$
for all $\alpha \in A_l$ and some $b_\alpha, c_\alpha \in \Z[\vec \h, \W]$ with the $b_\alpha$ all distinct and non-zero.  We can then write
$$ \Lambda_\Sys(\vec f) = \E_{\vec h\in [H]^D} \E_{m \in [M]} \int_X f_{\alpha_0}
(\prod_{\alpha \in A_{nl}} T^{R_\alpha(m,\vec h,W)} \nu)
\prod_{\alpha \in A_l} T^{b_\alpha m + c_\alpha} f_\alpha.$$
We introduce some new coarse-scale shift parameters $m_\alpha \in [\sqrt{M}]$ for $\alpha \in A_l$, thus the vector $\vec m := (m_\alpha)_{\alpha \in A_l}$ lies in $[\sqrt{M}]^{A_l}$.  We shift $m$ to $m - \sum_{\alpha \in A_l} m_\alpha$ and observe that
$$ \E_{m \in [M]} x_m = \E_{m \in [M]} x_{m - \sum_{\alpha \in A_l} m_\alpha} + o(1)$$
whenever $m_\alpha \in [\sqrt{M}]$ and $x_m \ll_\eps N^\eps$ for all $\eps > 0$.  Averaging this in $\vec m$ (cf. \eqref{xmh-1}) we obtain
$$ \E_{m \in [M]} x_m = \E_{\vec m \in [\sqrt{M}]^{A_l}} \E_{m \in [M]} x_{m - \sum_{\alpha \in A_l} m_\alpha} + o(1).$$
Applying this (and \eqref{nu-point}), and shifting the integral by the polynomial
\begin{equation}\label{q0def}
 Q_0 := \sum_{\alpha \in A_l} b_\alpha(\vec \h,\W) \m_\alpha,
 \end{equation}
we obtain
$$ \Lambda_\Sys(\vec f) = \E_{\vec h\in [H]^D} 
\E_{m \in [M]} \int_X 
\E_{\vec m \in [\sqrt{m}]^{A_l}} f_{\alpha_0, m, \vec h, \vec m, W}
\prod_{\alpha \in A_l} f_{\alpha, m, \vec h, \vec m,W} + o(1)$$
where 
\begin{align*}
f_{\alpha_0, \vec h, \vec m, W} &:= T^{Q_{\alpha_0}(\vec h, \vec m, W)} \left[f_{\alpha_0} 
\prod_{\alpha \in A_{nl}} T^{R_\alpha(m-\sum_{\alpha \in A_l} m_\alpha,\vec h,W)} \nu\right] \\
f_{\alpha,m,\vec h,\vec m,W} &:= T^{b_\alpha(\vec h,W) m + c_\alpha(\vec h,W) + \sum_{\beta \in A_l} 
(b_\beta(\vec h,W)-b_\alpha(\vec h,W)) m_\alpha} f_\alpha.
\end{align*}
The point of all these manipulations is that for each linear node $\alpha \in A_l$, $f_{\alpha,m,\vec h,\vec m,W}$ is independent of
the coarse-scale parameter $m_\alpha$.  Also observe the pointwise bound
$$ |f_{\alpha,m,\vec h, \vec m,W}| \leq \nu_{\alpha, m, \vec h, \vec m, W}$$
where
$$ \nu_{\alpha,m,\vec h,\vec m,W} := T^{b_\alpha(\vec h,W) m + c_\alpha(\vec h,W) + \sum_{\beta \in A_l} (b_\beta(\vec h,W)-b_\alpha(\vec h,W)) m_\beta} \nu.$$
By applying the weighted generalized von Neumann theorem (Proposition \ref{wgn}) in the $\vec m$ variables we thus have
\begin{equation}\label{ls}
 |\Lambda_\Sys(\vec f)| \leq \E_{\vec h\in [H]^D} 
\E_{m \in [M]} \int_X \| f_{\alpha_0, m, \vec h, \cdot, W} \|_{\Box^{A_l}(\nu)}
\prod_{\alpha \in A_l} \|\nu_{\alpha, m, \vec h, \cdot,W} \|_{\Box^{A_l \backslash \{\alpha\}}}^{1/2} + o(1)
\end{equation}
where the Gowers box norms $\Box^{A_l \backslash \{\alpha\}}$ and the weighted Gowers box norms $\Box^{A_l}(\nu)$ are defined\footnote{These norms will only make a brief appearance here; they are not used elsewhere in the main argument.} in Appendix \ref{gowers-sec}.
We now claim the estimate
\begin{equation}\label{boxest}
 \E_{\vec h\in [H]^D}
 \E_{m \in [M]} \int_X \|\nu_{\alpha, m, \vec h, \cdot,W} \|_{\Box^{A_l \backslash \{\alpha\}}}^{2^{|A_l|-1}} \ll 1
 \end{equation}
 for each $\alpha \in A_l$.
Indeed, the left-hand side can be expanded as
\begin{align*}
\E_{\vec h \in [H]^D} &\E_{\vec m^{(0)}, \vec m^{(1)} \in [H]^{A_l}}
 \E_{m \in [M]} \int_X
 \prod_{\omega \in \{0,1\}^{A_l \backslash \{\alpha\}}} \\
&\quad T^{b_\alpha(\vec h,W) m + c_\alpha(\vec h,W) + \sum_{\beta \in A_l} (b_\beta(\vec h,W)-b_\alpha(\vec h,W)) m_\beta^{(\omega_\beta)}} \nu.
\end{align*}
The distinctness of the $b_\beta$ ensures that the polynomial shifts of $\nu$ here are all distinct, and so by the polynomial forms condition (Definition \ref{polyform-def}) we obtain the claim (taking $\eta_1$ suitably small).

In view of \eqref{ls}, \eqref{boxest}, and H\"older's inequality, we see that
$$
 |\Lambda_\Sys(\vec f)| \ll_\Sys ( \E_{\vec h\in [H]^D} 
\E_{m \in [M]} \int_X \| f_{\alpha_0, m, \vec h, \cdot, W} \|_{\Box^{A_l}(\nu)}^{2^{|A_l|}} )^{1/2^{|A_l|}} + o(1).$$
Thus to prove Proposition \ref{gvn-induct} it suffices to show that
\begin{equation}\label{whf}
\E_{\vec h\in [H]^D} 
\E_{m \in [M]} \int_X \| f_{\alpha_0, m, \vec h, \cdot, W} \|_{\Box^{A_l}(\nu)}^{2^{|A_l|}}
=
\| f_{\alpha_0} \|_{U^{\vec Q([H]^D, W)}_{\sqrt{M}}}^{2^{|A_l|}} + o_\Sys(1)
\end{equation}
for some suitable $\vec Q \in \Z[\h_1,\ldots,\h_D,\W]$.  The left-hand side of \eqref{whf} can be expanded as a weighted average of $f_{\alpha_0}$ over parallelopipeds, or more precisely as
\begin{equation}\label{whf-1}
\E_{\vec h\in [H]^D} \E_{\vec m^{(0)}, \vec m^{(1)} \in [\sqrt{M}]^{A_l}} \int_X 
(\prod_{\omega \in \{0,1\}^{A_l}} T^{Q_0(\vec h, \vec m^{(\omega)}, W)} f_{\alpha_0})
w(\vec h, \vec m^{(0)}, \vec m^{(1)})
\end{equation}
where $\vec m^{(\omega)} := (m_\alpha^{(\omega_\alpha)})_{\alpha \in A_l}$ and
$w(\vec h, \vec m^{(0)}, \vec m^{(1)})$ is the weight
$$ w(\vec h, \vec m^{(0)}, \vec m^{(1)}) :=
\E_{m \in [M]}
\prod_{\omega \in \{0,1\}^{A_l}} 
\prod_{\alpha \in A_{nl}} T^{Q_0(\vec h, \vec m^{(\omega)}, W) + R_\alpha(m-\sum_{\alpha \in A_l} m_\alpha,\vec h,W)} \nu.$$

\subsection{Final Cauchy-Schwarz}\label{fcz-sec}

Let us temporarily drop the weight $w$ in \eqref{whf-1} and consider the unweighted average
$$ \E_{\vec h\in [H]^D} \E_{\vec m^{(0)}, \vec m^{(1)} \in [\sqrt{M}]^{A_l}} \int_X 
\prod_{\omega \in \{0,1\}^{A_l}} T^{Q_0(\vec h, \vec m^{(\omega)}, W)} f_{\alpha_0}.$$
Using \eqref{q0def} this is
$$ \E_{\vec h\in [H]^D} \E_{\vec m^{(0)}, \vec m^{(1)} \in [\sqrt{M}]^{A_l}} \int_X 
\prod_{\omega \in \{0,1\}^{A_l}} 
T^{\sum_{\alpha \in A_l} b_\alpha(\vec \h,\W) m^{(\omega_\alpha)}_\alpha} f_{\alpha_0}$$
which on comparison with \eqref{avgdef-2} is indeed of the form
$\| f_{\alpha_0} \|_{U^{\vec Q([H]^D, W)}_{\sqrt{M}}}^{2^{|A_l|}}$ for some\footnote{There is the (incredibly unlikely) possibility that $D=0$ or $D=1$, but by using the monotonicity of the Gowers norms (Lemma \ref{uk-concatenate}) one can easily increase $D$ to avoid this.}
$\vec Q \in \Z[\h_1,\ldots,\h_D,\W]$; note that because the $b_\alpha$ are non-zero, all the components of $\vec Q$ are non-zero.  Thus it suffices to show that
$$ \E_{\vec h\in [H]^D} \E_{\vec m^{(0)}, \vec m^{(1)} \in [\sqrt{M}]^{A_l}} \int_X 
(\prod_{\omega \in \{0,1\}^{A_l}} T^{Q_0(\vec h, \vec m^{(\omega)}, W)} f_{\alpha_0})
(w(\vec h, \vec m^{(0)}, \vec m^{(1)})-1) = o_\Sys(1).$$
Applying Cauchy-Schwarz \eqref{cauchy-schwarz} with the pointwise bound $|f_{\alpha_0}| \leq \nu$ we reduce to showing
$$ \E_{\vec h\in [H]^D} \E_{\vec m^{(0)}, \vec m^{(1)} \in [\sqrt{M}]^{A_l}} \int_X 
(\prod_{\omega \in \{0,1\}^{A_l}} T^{Q_0(\vec h, \vec m^{(\omega)}, W)} \nu)
(w(\vec h, \vec m^{(0)}, \vec m^{(1)})-1)^j = 0^j + o_\Sys(1)$$
for $j=0,2$ (with the usual convention $0^0 = 1$) which in turn will follow if we can show
$$ \E_{\vec h\in [H]^D} \E_{\vec m^{(0)}, \vec m^{(1)} \in [\sqrt{M}]^{A_l}} \int_X 
(\prod_{\omega \in \{0,1\}^{A_l}} T^{Q_0(\vec h, \vec m^{(\omega)}, W)} \nu)
w(\vec h, \vec m^{(0)}, \vec m^{(1)})^j = 1 + o_\Sys(1)$$
for $j=0,1,2$.  Let us just demonstrate this in the hardest case $j=2$, as it will be clear from the proof that the same argument also works for $j=0,1$ (as they involve fewer factors of $\nu$).  We expand the left-hand side as
\begin{align*}
&\E_{\vec h\in [H]^D} \E_{\vec m^{(0)}, \vec m^{(1)} \in [\sqrt{M}]^{A_l}} \E_{m,m' \in [M]} \int_X \\
&\quad (\prod_{\omega \in \{0,1\}^{A_l}} T^{Q_0(\vec h, \vec m^{(\omega)}, W)} \nu)\\
&\quad
\prod_{\omega \in \{0,1\}^{A_l}} 
\prod_{\alpha \in A_{nl}} 
T^{Q_0(\vec h, \vec m^{(\omega)}, W) + R_\alpha(m-\sum_{\alpha \in A_l} m_\alpha,\vec h,W)} \nu
T^{Q_0(\vec h, \vec m^{(\omega)}, W) + R_\alpha(m'-\sum_{\alpha \in A_l} m_\alpha,\vec h,W)} \nu
\end{align*}
One can then invoke the polynomial forms condition (Definition \ref{polyform-def}) one last time (again taking $\eta_1$ small enough) to verify that this is indeed $1 + o_\Sys(1)$.  Note that as every node in $A_{nl}$ is non-linear, the polynomials $R_\alpha$ have degree at least $2$, which ensures that the polynomials used to shift $\nu$ here are all distinct.  This concludes the proof of Proposition \ref{gvn-induct} in the linear case, and hence in general, and Theorem \ref{gvn} follows.
\endprf

\begin{remark}\label{multivar-rm}
One can define polynomial systems and weights (Definitions \ref{polysys}, \ref{weight}) 
for systems of  multivariable polynomials $R_{\alpha} \in \Z[\m_1,\ldots,\m_r,\h_1,\ldots,\h_D,W]$ (see for example \cite{leibman}). Following the steps of the PET induction (\ref{pet-sec}) and  
parallelopipedization (\ref{par-sec}) one
can  prove a multivariable version of the polynomial generalized von Neumann theorem (Theorem \ref{gvn}).
\end{remark}

\section{Polynomial dual functions}\label{polydual-sec}

This section and the next will be devoted to the proof of the structure theorem, Theorem \ref{struct-thm}.  In these sections we shall assume the notation of Section \ref{notation-sec}, and fix the bounded quantities $t \geq 0$, $d \geq 2$ and $\vec Q \in \Z[\h_1,\ldots,\h_t,\W]^d$.  As they are bounded we may permit all implicit constants in the $o()$ and $O()$ notation to depend on these quantities.  We also fix the pseudorandom measure $\nu$.
We shall abbreviate
$$ \|f\|_{U} := \|f\|_{U^{\vec Q([H]^t,W)}_{\sqrt{M}}}.$$
Roughly speaking, the objective here is to split any non-negative function bounded pointwise by $\nu$ to a non-negative function bounded pointwise by $1$, plus an error which is small in the $\|\|_U$ norm.  For technical reasons (as in \cite{gt-primes}) we will also need to exclude a small exceptional
set of measure $o(1)$, of which more will be said later.

Following \cite{gt-primes}, our primary tool for understanding the $U$ norm shall be via the concept of a \emph{dual function} of a function $f$ associated to this norm.

\begin{definition}[Dual function]  If $f: X \to \R$ is a function, we define the \emph{dual function}
$\D f: X \to \R$ by the formula
$$ \D f :=  \E_{\vec h \in [H]^t} 
\E_{\vec m^{(0)}, \vec m^{(1)} \in [\sqrt{M}]^d} 
\prod_{(\omega_1,\ldots,\omega_d) \in \{0,1\}^d \backslash \{0\}^d}
T^{\sum_{i \in [d]} (m^{(\omega_i)}_i - m^{(0)}_i) Q_i(\vec h, W)} f.
$$
where $m^{(i)} = (m^{(i)}_1,\ldots,m^{(i)}_d)$ for $i=0,1$.
\end{definition}

From \eqref{avgdef-2} and the translation invariance of the integral $\int_X$ we obtain the fundamental relationship
\begin{equation}\label{dualfn}
\| f \|_U^{2^d} = \int_X f \D f.
\end{equation}
Thus we have a basic dichotomy: either $f$ has small $U$ norm, or else it correlates with its own dual function\footnote{In the language of infinitary ergodic theory, it will be the dual functions which generate (in the measure-theoretic sense) the characteristic factor for the $U$ norm.  The key points will be that the dual functions are essentially bounded, and that $\nu-1$ is essentially orthogonal to the characteristic factor.}.  As in \cite{gt-primes}, it is the iteration of this dichotomy via a stopping time argument which shall power the proof of Theorem \ref{struct-thm}.

For future reference we observe the trivial but useful facts that $\D$ is monotone and homogeneous of degree $2^d-1$:
\begin{equation}\label{dualmono}
|f| \leq g \hbox{ pointwise} \implies |\D f| \leq \D g \hbox{ pointwise};
\end{equation}
\begin{equation}\label{dualhomog}
\D(\lambda f) = \lambda^{2^d-1} \D f \hbox{ for all } \lambda \in \R.
\end{equation}

We will need two key facts about dual functions, both of which follow primarily from the polynomial correlation condition.  The first, which is fairly easy, is that dual functions are essentially bounded.

\begin{proposition}[Essential boundedness of dual functions]\label{essbound}  
Let $f: X \to \R$ obey the pointwise bound $|f| \leq \nu + 1$.  Then for any integer $K \geq 1$ we have
the moment estimates
\begin{equation}\label{moment}
\int_X |\D f|^K (\nu+1) \leq 2 (2^{2^d-1})^K + o_K(1).
\end{equation}
In particular, if we define the \emph{global bad set}
\begin{equation}\label{omega-def}
 \Omega_0 := \{ x \in X: \D \nu(x) \geq 2^{2^d} \}
 \end{equation}
then we have the measure bound
\begin{equation}\label{measu}
 \int_X (\D \nu)^K 1_{\Omega_0} (\nu+1) = o_K(1)
 \end{equation}
for all $K \geq 0$, and the pointwise bound
\begin{equation}\label{dual-point}
|\D f| (1 - 1_{\Omega_0}) \leq 2^{2^d}.
\end{equation}
\end{proposition}

\begin{remark} In \cite{gt-primes}, the correlation conditions imposed on $\nu$ were strong enough that one could bound the dual function $\D f$ uniformly by $2^{2^d-1}+o(1)$, thus removing the need for a global bad set $\Omega_0$.  One could do something similar here by strengthening the correlation condition.  However, we were then unable to establish Theorem \ref{major}, i.e. we were unable to construct a measure concentrated on almost primes which obeyed this stronger correlation condition.  The basic difficulty is that the polynomials in $\vec Q$ could contain a number of common factors which could significantly distort functions such as $\D \nu$ at some rare points (such as the origin).  Fortunately, the presence of a small global bad set does not significantly impact our analysis (similarly to how sets of measure zero have no impact on ergodic theory), especially given how it does not depend on $f$.  In practice, $K$ will get as large as $1/\eta_6$, but no greater.
\end{remark}

\begin{proof}  We begin with \eqref{moment}.  By \eqref{dualmono}, \eqref{dualhomog} it suffices to show that
$$ \int_X \left[\D \frac{\nu+1}{2}\right]^K \frac{\nu+1}{2} = 1 + o_K(1);$$
in view of Lemma \ref{avglem}, it suffices to show that
$$ \int_X (\D \nu)^K \nu = 1 + o_K(1).$$
The left-hand side can be expanded as
\begin{align*} \E_{\vec h^{(1)}, \ldots, \vec h^{(K)} \in [H]^t} 
\int_X &\prod_{k \in [K]}
[
\E_{\vec m^{(0)}, \vec m^{(1)} \in [\sqrt{M}]^d}\\ 
&\quad \prod_{(\omega_1,\ldots,\omega_d) \in \{0,1\}^d \backslash \{0\}^d}
T^{\sum_{i \in [d]} (m^{(\omega_i)}_i - m^{(0)}_i) Q_i(\vec h^{(k)}, W)} \nu ]
\nu
\end{align*}
But this is $1+o_K(1)$ from \eqref{polycor-eq} (with $D=d$, $D' = 0$, $D'' = K t$, and $L=1$, with the $\vec Q_{j,k}$ and $\vec S_l$ vanishing).  This proves \eqref{moment}.  From Chebyshev's inequality this implies that
$$ \int_X (\D \nu)^{K'} 1_{\Omega_0} (\nu+1) \leq \frac{2 (2^{2^d-1})^{K'}}{2^K} + o_K(1)$$
for any $0 \leq K' < K$.  For fixed $K'$, the right-hand side can be made arbitrarily small by taking $K$ large, and then choosing $N$ large depending on $K$; thus the left-hand side is $o(1)$, which is \eqref{measu}.  Finally, \eqref{dual-point} follows from \eqref{dualmono} and \eqref{omega-def}.
\end{proof}

The global bad set $\Omega_0$ is somewhat annoying to deal with.  Let us remove it by defining the \emph{modified dual function} $\tilde \D f$ of $f$ as
$$ \tilde \D f := (1 - 1_{\Omega_0}) \D f.$$
Then Proposition \ref{essbound} and \eqref{dualfn} immediately imply

\begin{corollary}[Boundedness of modified dual function]\label{bound-mod}  
Let $f: X \to \R$ obey the pointwise bound $|f| \leq \nu + 1$.  Then $\tilde \D f$ takes
values in the interval 
\begin{equation}\label{I-def}
I := [-2^{2^d}, 2^{2^d}].
\end{equation}
Furthermore we have the correlation property
\begin{equation}\label{dualfn-mod}
\int_X f \tilde \D f = \| f \|_U^{2^d} + o(1).
\end{equation}
\end{corollary}

The second important estimate is easy to state, although non-trivial to prove:

\begin{proposition}[$\nu-1$ orthogonal to products of modified dual functions]\label{nuorthog}  
Let $1 \leq K \leq 1/\eta_6$ be an integer, and let $f_1,\ldots,f_K: X \to \R$ be functions with the pointwise bounds $|f_k| \leq \nu+1$ for all $k \in [K]$.  Then
\begin{equation}\label{feeble}
 \int_X \tilde \D f_1 \ldots \tilde \D f_K (\nu - 1) = o(1).
\end{equation}
\end{proposition}

\begin{remark} Note
that \eqref{moment} already gives an upper bound of $O(1)$ for \eqref{feeble}; the whole point is thus to extract enough cancellation from the $\nu-1$ factor to upgrade this bound to $o(1)$.
\end{remark}

The rest of the section is devoted to the proof of Proposition \ref{nuorthog}. The argument follows that of \cite[Section 6]{gt-primes}, and is based on a large number of applications of the Cauchy-Schwarz inequality, and the polynomial correlation condition (Definition \ref{polycor-def}).  The arguments here are not used again elsewhere in this paper, and so the rest of this section may be read independently of the remainder of the paper.

We begin with a very simple reduction: from \eqref{measu} we can replace the modified dual functions $\tilde \D f_i$ with their unmodified counterparts $\D f_i$.  
Our task is then to show
\begin{equation}\label{nuorthog-prod}
 \int_X \D f_1 \ldots \D f_K (\nu - 1) = o(1).
\end{equation}
 
\subsection{A model example}

Before we prove Proposition \ref{nuorthog} in general, it is instructive to work with a simple example to illustrate the idea.  Let us take an oversimplified toy model of the dual function $\D f$, namely
$$ \D f := \E_{h \in [H]} \E_{m \in [\sqrt{M}]} T^{mh} f.$$
This does not quite correspond to a local Gowers norm\footnote{However, the slight variant
$\E_{h \in [H]} \E_{m,m' \in [\sqrt{M}]} T^{(m-m')h} f$ does correspond to a (very simple) local Gowers norm, with $t=d=1$ and $\vec Q = (\h_1)$.}, but will serve as an illustrative model nonetheless.  Pick functions $f_1,\ldots,f_K$ with the pointwise bound $|f_k| \leq \nu$ for $k \in [K]$ and consider the task of showing \eqref{nuorthog-prod}. We expand the left-hand side as
\begin{equation}\label{duallhs}
 \E_{h_1,\ldots,h_K \in [H]} \E_{m_1,\ldots,m_K \in [\sqrt{M}]}
\int_X T^{m_1 h_1} f_1 \ldots T^{m_K h_K} f_K (\nu-1).
\end{equation}
Note that we cannot simply take absolute values and apply the pseudorandomness conditions, as these will give bounds of the form $O(1)$ rather than $o(1)$.  One could instead attempt to apply the Cauchy-Schwarz inequality many times (as in the previous section), however the fact that $K = O(1/\eta_6^2)$ could be very large compared to the pseudorandomness parameter $1/\eta_1$ defeats a naive implementation of this idea.  Instead, we must perform a change of variables to introduce two new parameters $n^{(0)}, n^{(1)}$ to average over (and which only requires a single Cauchy-Schwarz to estimate) rather than $K$ parameters (which would essentially require $K$ applications of Cauchy-Schwarz).

More precisely, we introduce two slightly less coarse-scale parameters $n^{(0)}, n^{(1)} \in [M^{1/4}]$ than $m_1,\ldots,m_K$.  Define the multipliers $\hat h_k := \prod_{k' \in [K] \backslash k} h_{k'}$, thus $\hat h_k = O( H^{K-1} )$, which is small compared to $M^{1/4}$ by the relative sizes of $\eta_7, \eta_6, \eta_2$.  Shifting each of the $m_k$ by $\hat h_k (n^{(1)}-n^{(0)})$ and using \eqref{nu-point}, we conclude that \eqref{duallhs} is equal to
$$\E_{h_1,\ldots,h_K \in [H]} \E_{m_1,\ldots,m_K \in [\sqrt{M}]}
\int_X \left[\prod_{k \in [K]} T^{(m_k+\hat h_k (n^{(1)}-n^{(0)})) h_k} f_k\right] (\nu-1) + o_K(1)
$$
for all $n^{(0)}, n^{(1)} \in [M^{1/4}]$.  Averaging over all $n^{(0)}, n^{(1)}$ and shifting the integral by 
$$n^{(0)} h_1 \ldots h_K = \hat h_1 n^{(0)} h_1 = \ldots = \hat h_K n^{(0)} h_K$$
we can thus write \eqref{duallhs} as
\begin{align*}
&\E_{h_1,\ldots,h_K \in [H]} \E_{m_1,\ldots,m_K \in [\sqrt{M}]} \E_{n^{(0)}, n^{(1)} \in [M^{1/4}]}\\
&\quad \int_X T^{n^{(1)}h_1 \ldots h_K} [T^{m_1 h_1} f_1 \ldots T^{m_K h_K} f_K] 
 T^{n^{(0)} h_1 \ldots h_K} (\nu-1) + o_K(1)
\end{align*}
 which we may factorize as
\begin{align*}
\E_{h_1,\ldots,h_K \in [H]} \int_X &\left[ \E_{n^{(1)} \in [M^{1/4}]} \prod_{k \in [K]} \E_{m \in [\sqrt{M}]}
T^{n^{(1)} h_1 \ldots h_k + m h_k} f_k \right]\\
& \left[\E_{n^{(0)} \in [M^{1/4}]} T^{n^{(0)} h_1 \ldots h_K} (\nu-1)\right]
+ o_K(1).
\end{align*}
By Cauchy-Schwarz it thus suffices to show that
$$ \E_{h_1,\ldots,h_K \in [H]} \int_X \left[ \E_{n^{(1)} \in [M^{1/4}]} \prod_{k \in [K]} \E_{m \in [\sqrt{M}]}
T^{n^{(1)} h_1 \ldots h_k + m h_k} f_k \right]^2 \ll_K 1$$
and
$$ \E_{h_1,\ldots,h_K \in [H]} \int_X 
\left[\E_{n^{(0)} \in [M^{1/4}]} T^{n^{(0)} h_1 \ldots h_K} (\nu-1)\right]^2 = o_K(1).$$
To prove the first estimate, we estimate $f_k$ by $\nu$ and expand out the square to reduce to showing
$$ \E_{h_1,\ldots,h_K \in [H]} \E_{n^{(1)}, n^{(2)} \in [M^{1/4}]}
\int_X \prod_{i=1}^2 \prod_{k \in [K]} \E_{m \in [\sqrt{M}]}
T^{n^{(i)} h_1 \ldots h_k + m h_k} \nu  \ll_K 1,$$
but this follows from the correlation condition \eqref{polycor-eq} (for $\eta_1$ small enough\footnote{It is important to note however that $\eta_1$ does not have to be small relative to $K$ or to parameters such as $\eta_7$.}).  To prove the second estimate, we again expand out the square and reduce to showing
$$ \E_{h_1,\ldots,h_K \in [H]} \int_X 
\left[\E_{n^{(0)} \in [M^{1/4}]} T^{n^{(0)} h_1 \ldots h_K} \nu\right]^j = 1 + o_K(1).$$
for $j=0,1,2$, which will again follow from \eqref{polycor-eq} for $\eta_1$ small enough.

\subsection{Conclusion of argument}

Now we prove \eqref{nuorthog-prod} in the general case.  We may take $d \geq 1$, since the $d=0$ case follows from \eqref{measure}.
We expand the left-hand side as
$$ \E_{\vec h} \E_{\vec m}
\int_X \left[\prod_{k \in [K]} 
\prod_{\omega \in \{0,1\}^d \backslash \{0\}^d}
T^{\sum_{i \in [d]} (m^{(\omega_i)}_{i,k} - m^{(0)}_{i,k}) Q_i(\vec h^{(k)}, W)} f_k\right] (\nu-1).
$$
where $\omega = (\omega_1,\ldots,\omega_d)$ we use the abbreviations
\begin{align*}
\E_{\vec h} &:= \E_{\vec h^{(1)}, \ldots, \vec h^{(K)} \in [H]^{t}} \\
\E_{\vec m} &:= \E_{m^{(\omega)}_{i,k} \in [\sqrt{M}] \forall \omega \in \{0,1\}, i \in [d], k \in [K]}.
\end{align*}
We introduce moderately coarse-scale parameters $n^{(0)}_i, n^{(1)}_i \in [M^{1/4}]$ for $i \in [d]$ and the multipliers
$$ \hat h_{k,i} = \hat h_{k,i}(\vec h, W) := \prod_{k' \in [K] \backslash \{k\}} Q_i(\vec h^{(k')}, W).$$
Observe that $\hat h_{k,i} = O( H^{O(K)} )$, which will be much smaller than $M^{1/4}$ by the relative sizes of $\eta_7, \eta_6, \eta_2$.
Shifting each $m^{(1)}_{i,k}$ by $\hat h_{k,i} (n^{(1)}_i - n^{(0)}_i)$ and using \eqref{nu-point}, we can then rewrite \eqref{nuorthog-prod} as
$$ \E_{\vec h} \E_{\vec m}
\int_X [\prod_{k \in [K]} 
\prod_{\omega \in \{0,1\}^d \backslash \{0\}^d}
T^{\sum_{i \in [d]} (m^{(\omega_i)}_{i,k} - m^{(0)}_{i,k} + \hat h_{k,i} (n^{(\omega_i)}_i-n^{(0)}_i)) Q_i(\vec h^{(k)}, W)} f_k] (\nu-1) + o_K(1)
$$
for any $n^{(0)}_1,\ldots,n^{(0)}_d,n^{(1)}_1,\ldots,n^{(1)}_d \in [M^{1/4}]$.
Now from construction we have $\hat h_{k,i} Q_i(\vec h^{(k)}, W) =  b_i$, where
$$ b_i = b_i(\vec h, W) := \prod_{k \in [K]} Q_i(\vec h^{(k)}, W)$$
(note that $b_i \neq 0$ by hypothesis on $\vec Q$)
and on averaging over the $n$ variables, we can write the left-hand side of \eqref{nuorthog-prod} as
$$ \E_{\vec h, \vec m, \vec n^{(0)}, \vec n^{(1)}_1} \int_X
\prod_{\omega \in \{0,1\}^d}
T^{\sum_{i \in [d]} (n^{(\omega_i)}_i-n^{(0)}_i) b_i} g_{\omega,\vec h,\vec m}
+ o_K(1)$$
where $\vec n^{(i)} = (n^{(i)}_0,\ldots,n^{(i)}_d)$ will be understood to range over $[M^{1/4}]^d$ for $i=0,1$,
$$ g_{\omega,\vec h, \vec m} := 
\prod_{k \in [K]}
T^{\sum_{i \in [d]} (m^{(\omega_i)}_{i,k} - m^{(0)}_{i,k}) Q_i(\vec h^{(k)}, W)} f_k.$$
for $\omega \in \{0,1\}^d \backslash \{0\}^d$, and
$$ g_{\{0\}^d,\vec h, \vec m} := \nu - 1.$$
Shifting the integral by $T^{\sum_{i \in [d]} n^{(0)}_i b_i}$ we can rewrite this as
$$ \E_{\vec h, \vec m}  \int_X
\E_{\vec n^{(0)}, \vec n^{(1)}}
\prod_{\omega \in \{0,1\}^d}
T^{\sum_{i \in [d]} n^{(\omega_i)}_i b_i} g_{\omega,\vec h,\vec m}
+ o_K(1).$$
Now use the Cauchy-Schwarz-Gowers inequality \eqref{gcz} to obtain the pointwise estimate
\begin{align*}
&|\E_{\vec n^{(0)}, \vec n^{(1)}} 
\prod_{\omega \in \{0,1\}^d} T^{\sum_{i \in [d]} n^{(\omega_i)}_i b_i} g_{\omega,\vec h,\vec m}|\\
&\quad \leq 
\prod_{\omega' \in \{0,1\}^d} 
(\E_{\vec n^{(0)}, \vec n^{(1)}}
\prod_{\omega \in \{0,1\}^d} T^{\sum_{i \in [d]} n^{(\omega_i)}_i b_i} g_{\omega',\vec h,\vec m})^{1/2^d}.
\end{align*}
By H\"older's inequality, we thus see that to prove \eqref{nuorthog-prod} it suffices to show that the quantity
\begin{equation}\label{evac}
 \E_{\vec h, \vec m}  \int_X
\E_{\vec n^{(0)}, \vec n^{(1)}}
\prod_{\omega \in \{0,1\}^d}
T^{\sum_{i \in [d]} n^{(\omega_i)}_i b_i} g_{\omega',\vec h,\vec m} 
\end{equation}
is $O_K(1)$ when $\omega' \in \{0,1\}^d \backslash \{0\}^d$ and is $o_K(1)$ when $\omega' = 0^d$.

Let us first deal with the case when $\omega' \neq 0^d$.  Our task is to show that
$$
 \E_{\vec h, \vec m}  \int_X
\E_{\vec n^{(0)}, \vec n^{(1)}}
\prod_{\omega \in \{0,1\}^d}
T^{\sum_{i \in [d]} n^{(\omega_i)}_i b_i}  
\prod_{k \in [K]}
T^{\sum_{i \in [d]} (m^{(\omega'_i)}_{i,k} - m^{(0)}_{i,k}) Q_i(\vec h^{(k)}, W)} f_k = O_K(1).$$
We can bound $f_k$ pointwise by $\nu$, and factorise the left-hand side as
$$
 \E_{\vec h, \vec n^{(0)}, \vec n^{(1)}}
\int_X
\prod_{k \in [K]}
\E_{\vec m^{(0)}, \vec m^{(1)} \in [\sqrt{M}]^d}
\prod_{\omega \in \{0,1\}^d}
T^{\sum_{i \in [d]} (m^{(\omega'_i)}_i - m^{(0)}_i) Q_i(\vec h^{(k)}, W) + n^{(\omega_i)}_i b_i} \nu. $$
But this is $1 + o_K(1) = O_K(1)$ by \eqref{polycor-def} with $L=0$ (here we use the fact that the $b_i$ are non-zero polynomials of $\vec h$ and $W$).  Here we need $\eta_1$ to be sufficiently small depending on $t,d,\vec Q$ but not on $K$.

Finally, we have to deal with the case $\omega = 0^d$.  Since $g_{\omega',\vec h,\vec m} = \nu-1$ and $b_i = b_i(\vec h, W)$ are independent of $W$, we can rewrite \eqref{evac} as
$$
 \E_{\vec h, \vec n^{(0)}, \vec n^{(1)}} \int_X
\prod_{\omega \in \{0,1\}^d}
T^{\sum_{i \in [d]} n^{(\omega_i)}_i b_i(\vec h,W)} (\nu-1)$$
and so by the binomial formula it suffices to show that
$$
\E_{\vec h, \vec n^{(0)}, \vec n^{(1)}} \int_X
\prod_{\omega \in A}
T^{\sum_{i \in [d]} n^{(\omega_i)}_i b_i(\vec h,W)} \nu = 1 + o(1)$$
for all $A \subset \{0,1\}^d$.  But this follows from the polynomial correlation condition \eqref{polycor-eq} (with $K=0$), again taking $\eta_1$ sufficiently small depending on $t,d,\vec Q$.  This concludes the proof of \eqref{nuorthog-prod}, and hence Proposition \ref{nuorthog}.

\section{Proof of structure theorem}\label{struct-sec}

We can now complete the proof of the structure theorem by using the arguments of
\cite[\S 7-8]{gt-primes} more or less verbatim.  In fact these arguments can be abstracted as follows.

\begin{theorem}[Abstract structure theorem]\label{abstract}  Let $I$ be an interval bounded by $O(1)$.  Let $\nu: X \to \R^+$ be any measure, and let $f \mapsto \tilde \D f$ be a (nonlinear) operator obeying the following properties:
\begin{itemize}
\item If we have the pointwise bound $|f| \leq \nu+1$, then $\tilde \D f: X \to I$ takes values in $I$, in particular
\begin{equation}\label{dfb}
\tilde \D f = O(1).
\end{equation}
\item If $1 \leq K \leq 1/\eta_6$, and $f_1,\ldots,f_K: X \to \R$ are functions with the pointwise bound $|f_k| \leq \nu+1$ for all $k \in [K]$, then \eqref{feeble} holds.
\end{itemize}
Then for any $g: X \to \R^+$ with the pointwise bound $0 \leq g \leq \nu$, there exist functions $g_{U^\perp}, g_U: X \to \R$ obeying
the estimates \eqref{gdecomp}, \eqref{struct-bound}, \eqref{most-mass}, and
\begin{equation}\label{gud}
|\int_X g_U \tilde \D g_U| \leq \eta_4.
\end{equation}
\end{theorem}

Indeed, Theorem \ref{struct-thm} immediately follows by applying Theorem \ref{abstract}, with \eqref{unstruct-bound} following from \eqref{gud} and \eqref{dualfn-mod}.

In the remainder of this section we prove Theorem \ref{abstract}.  Henceforth we fix $I, \nu, \tilde \D$ obeying the hypotheses of the theorem. 

\subsection{Factors}  

As in \cite{gt-primes} we shall recall the very useful notion of a \emph{factor} from ergodic theory,
though for our applications we actually only need the finitary version of this concept.

Let us set $\X$ to be the probability space $\X=(X,\B_X,\mu_X)$, where $X = \Z_N$, $\B_X = 2^X$ is the power set of $X$, and $\mu_X$ is the uniform probability measure on $\X$.  We define a \emph{factor}\footnote{In infinitary ergodic theory one also requires the probability spaces $\X, \Y$ to be invariant under the shift $T$, and for the factor map $\pi$ to respect the shift.  In the finitary setting it is unrealistic to demand these shift-invariances, for if $N$ were prime then this would mean that there were no non-trivial factors whatsoever.  While there are concepts of ``approximate shift-invariance'' which can be used as a substitute, see \cite{tao:ergodic}, we will fortunately not need to use them here, as the remainder of the argument does not even involve the shift $T$ at all.} to a quadruple $\Y = (Y,\B_Y,\mu_Y,\pi_Y)$, where $(Y,\B_Y,\mu_Y)$ is a probability space (thus $\B_Y$ is a $\sigma$-algebra on $Y$ and $\mu_Y$ is a probability measure on $\B_Y$) together with a measurable map $\pi:X \rightarrow Y$ such that $(\pi_Y)_*\mu_X=\mu_Y$, or in other words $\mu_X( \pi_Y^{-1}(E) ) = \mu_Y(E)$ for all $E \in \B_Y$.  The factor map $\pi_Y$ induces the pullback map $\pi_Y^*:L^2(\Y) \rightarrow L^2(\X)$ and its adjoint $(\pi_Y)_*:L^2(\X) \rightarrow L^2(\Y)$, where $L^2(\X)$ is the usual Lebesgue space of square-integrable functions on $\X$. 
We refer to the projection $\pi_Y^* (\pi_Y)_*: L^2(\X) \to L^2(\X)$ as the \emph{conditional expectation operator}, and denote $\pi_Y^* (\pi_Y)_*(f)$ by $\E(f|\Y)$; this is a linear self-adjoint orthogonal projection from $L^2(\X)$ to $\pi_Y^* L^2(\Y)$.
 
The conditional expectation operator is in fact completely determined by the $\sigma$-algebra
$\pi_Y^{-1}(\B_Y) \subset \B_X$.  As $X$ is finite (with every point having positive measure), $\pi_Y^{-1}(\B_Y)$ is generated by a partition of $X$ into atoms (which by abuse of notation we refer to as atoms of the factor $\Y$), and the conditional expectation is given explicitly by the formula
$$ \E( f|\Y )(x) = \E_{y \in \B(x)} f(y)$$
where $\B(x)$ is the unique atom of $\pi^{-1}(\B_Y)$ which contains $x$.  We refer to the number of atoms\footnote{It would be more natural to work instead with the \emph{entropy} of $\Y$ rather than the complexity, but the entropy is a slightly more technical concept and so we have avoided its use here for simplicity.} of $\Y$ as the \emph{complexity} of the factor $\Y$.  By abuse of notation we say that a function $f: X \to \R$ is \emph{measurable with respect to $\Y$} if it is measurable with respect to $\pi_Y^{-1}(\B_Y)$, or equivalently if it is constant on all atoms of $\Y$.  Thus for instance $(\pi_Y)^* L^q(\Y)$ consists of the functions in $L^q(\X)$ which are measurable with respect to $\Y$.

If $\Y = (Y,\B_Y,\mu_Y,\pi_Y)$ and $\Y' = (Y',\B_{Y'},\mu_{Y'},\pi_{Y'})$ are two factors, we may form their join $\Y \vee \Y' = (Y \times Y', \B_Y \times \B_{Y'}, \mu_Y \times \mu_{Y'}, \pi_Y \oplus \pi_{Y'})$ in the obvious manner; note that the atoms of $\Y \vee \Y'$ are simply the non-empty intersections of atoms of $\Y$ with atoms of $\Y'$, and so any function which is measurable with respect to $\Y$ or $\Y'$ is automatically measurable with respect to $\Y \vee \Y'$.

Note that any function $f: \X \to \R$ automatically generates a factor $(\R, \B_\R, f_* \mu_X, f)$, where $\B_\R$ is the Borel $\sigma$-algebra,which is the minimal factor with respect to which $f$ is (Borel-) measurable.  In our finitary setting it turns out we need a discretized version of this construction, which we give as follows.

\begin{proposition}[Each function generates a factor]\label{sigma-gen}
For any function $G: X \to I$ there exists a factor $\Y(G)$ with the following properties:

\begin{itemize}

\item \textup{($G$ lies in its own factor)} For any factor $\Y'$, we have
\begin{equation}\label{trivial}
G = \E( G | \Y(G) \vee \Y') + O(\eta_4^2).
\end{equation}

\item \textup{(Bounded complexity)} $\Y(G)$ has at most $O_{\eta_4}(1)$ atoms.

\item \textup{(Approximation by continuous functions of $G$)}  If $A$ is any atom in $\Y(G)$, then there exists a polynomial $\Psi_A: \R \to \R$ of degree $O_{\eta_5}(1)$ and coefficients $O_{\eta_5}(1)$ such that 
\begin{equation}\label{psi-contract}
\Psi_{A}(x) \in [0,1] \hbox{ for all } x \in I
\end{equation}
and
\begin{equation}\label{nup-bound}
\int_X |1_A - \Psi_{A}(G)| (\nu + 1) \ll \eta_5.
\end{equation}

\end{itemize}

\end{proposition}

\begin{proof}  This is essentially \cite[Proposition 7.2]{gt-primes}, but we shall give a complete proof here for the convenience of the reader.  

We use the probabilistic method. Let $\alpha$ be a real number in the interval $[0,1]$, chosen at random.  We then define the factor
$$ \Y(G) := ( \R, \B_{\eta_4^2,\alpha}, G_* \mu_X, G )$$
where $\B_{\eta_4^2,\alpha}$ is the $\sigma$-algebra on the real line $\R$ generated by the intervals $[(n+\alpha)\eta_4^2, (n+\alpha+1)\eta_4^2)$ for $n \in \Z$.  This is clearly a factor of $\X$, with atoms
$A_{n,\alpha} := G^{-1}([(n+\alpha)\eta_4^2, (n+\alpha+1)\eta_4^2))$.  Since $G$ ranges in $I$, and we allow constants to depend on $I$, it is clear that there are at most $O_{\eta_4}(1)$ non-empty atoms and that $G$ fluctuates by at most $O(\eta_4^2)$ on each atom, which yields the first two desired properties.  It remains to verify that with positive probability, the approximation by continuous functions property holds for all atoms $A_{n,\alpha}$.  By the union bound, it suffices to show that each individual atom $A_{n,\alpha}$ has the approximation property with probability $1 - O(\eta_5)$.

By the Weierstrass approximation theorem, we can find for each $\alpha$ a polynomial $\Psi_{A_{n,\alpha}}$ obeying \eqref{psi-contract} which is equal to
$1_{[(n+\alpha)\eta_4^2, (n+\alpha+1)\eta_4^2)} + O(\delta)$ outside of the set
$$ E_{n,\alpha} := [(n+\alpha-\eta_5^2)\eta_4^2, (n+\alpha+\eta_5^2)\eta_4^2] \cup
[(n+\alpha+1-\eta_5^2)\eta_4^2, (n+\alpha+1+\eta_5^2)\eta_4^2].$$
Simple compactness arguments allow us to take $\Psi_{A_{n,\alpha}}$ to have degree $O_{\eta_5}(1)$ and coefficients $O_{\eta_5}(1)$.  Since
$$ 1_{A_{n,\alpha}} = 1_{[(n+\alpha)\eta_4^2, (n+\alpha+1)\eta_4^2)}(G),$$
we thus conclude (from \eqref{measure}) that
$$ \int_X |1_A - \Psi_{A_{n,\alpha}}(G)| (\nu + 1) \ll \eta_5 + \int_X 1_{E_{n,\alpha}}(G) (\nu+1).$$
By Markov's inequality, it thus suffices to show that
$$ \int_0^1 [\int_X 1_{E_{n,\alpha}}(G) (\nu+1)]\ d\alpha \ll \eta_5^2.$$
But this follows from Fubini's theorem, \eqref{measure}, and the elementary pointwise estimate
$$ \int_0^1 1_{E_{n,\alpha}}(G)\ d\alpha \ll \eta_5^2.$$
\end{proof}

Henceforth we set $\Y(G)$ to be the factor given by the above proposition.  A key consequence of the hypotheses of Theorem \ref{abstract} is that $\nu-1$ is well distributed with respect to any finite combination of these factors:

\begin{proposition}[$\nu$ uniformly distributed wrt dual function factors]\label{nu-uniform}  Let $K \geq 1$ be an integer with $K = O_{\eta_4}(1)$, and let $f_1,\ldots,f_K: X \to \R$ be functions with the pointwise bounds $|f_k| \leq \nu+1$ for all $1 \leq k \leq K$.  Let $\Y := \Y(\tilde \D f_1) \vee \ldots \vee
  \Y(\tilde \D f_K)$.  Then we have
\begin{equation}\label{tide}
 \tilde \D f_k = \E( \tilde \D f_k | \Y ) + O(\eta_4^2)
\end{equation} 
for all $k \in [K]$, and we have a $\Y$-measurable set $\Omega \subset X$ obeying the smallness bound
\begin{equation}\label{om-small}
\int_X 1_\Omega (\nu+1) \ll_{\eta_4} \eta_5^{1/2}  
\end{equation}
and we have the pointwise bound
\begin{equation}\label{nu-dist}
 |(1 - 1_\Omega) \E( \nu-1|\Y )| \leq O_{\eta_4}(\eta_5^{1/2}).
\end{equation} 
\end{proposition}

\begin{proof} We repeat the arguments from \cite[Proposition 7.3]{gt-primes}.  The claim \eqref{tide}
follows immediately from \eqref{trivial}, so we turn to the other two properties.  Since each $\Y(\tilde \D f_k)$ is generated by $O_{\eta_4}(1)$ atoms, $\Y$ is generated by $O_{\eta_4,K}(1) = O_{\eta_4}(1)$ atoms.
Call an atom $A$ of $\Y$ \emph{small} if $\int_X 1_A (\nu+1) \leq \eta_5^{1/2}$, and let $\Omega$ be
the union of all the small atoms, then $\Omega$ is clearly $\Y$-measurable and obeys \eqref{om-small}.
It remains to prove \eqref{nu-dist}, or equivalently that
$$ \frac{ \int_X 1_A (\nu - 1)}{ \int_X 1_A } = \E_{y \in A} \nu(y)-1 \ll_{\eta_4} \eta_5^{1/2} + o(1)$$
for all non-small atoms $A$.  

Fix $A$. Since $A$ is not small, we have
$$\int_X 1_A (\nu - 1) + 2\int_X 1_A  = \int_X 1_A (\nu+1) > \eta_5^{1/2}.$$
Hence it will suffice to show that
$$\int_X 1_A (\nu - 1) \ll_{\eta_4} \eta_5 + o(1).$$
On the other hand, since $A$ is the intersection of atoms $A_1,\ldots,A_K$ from $\Y(\tilde \D f_1)$, $\ldots$, $\Y(\tilde \D f_K)$, we see from Proposition \ref{sigma-gen} and an easy induction argument that there exists a polynomial $\Psi: \R^K \to \R$ of degree $O_{\eta_5}(1)$ and coefficients
$O_{\eta_5}(1)$ which maps $I^K$ into $[0,1]$ such that
$$ \int_X |1_A - \Psi( \tilde \D f_1, \ldots, \tilde \D f_K )| (\nu + 1 ) \ll_{\eta_4} \eta_5.$$
In particular
$$ \int_X (1_A - \Psi( \tilde \D f_1, \ldots, \tilde \D f_K )) (\nu - 1) \ll_{\eta_4} \eta_5.$$
On the other hand, by decomposing $\Psi$ into monomials and using \eqref{feeble} (assuming $\eta_6$ sufficiently small depending on $\eta_5$) we have
$$ \int_X \Psi( \tilde \D f_1, \ldots, \tilde \D f_K ) (\nu - 1) = o(1)$$
and the claim follows (we can absorb the $o(1)$ error by taking $N$ large enough).
\end{proof}

\subsection{The inductive step}

The proof of the abstract structure theorem proceeds by a stopping time argument.  To clarify this
argument we introduce a somewhat artificial definition.

\begin{definition}[Structured factor]  A \emph{structured factor} is a tuple $\Y_K = (\Y_K$, $K$, $F_1$, $\ldots$, $F_K$, $\Omega_K)$, where $K \geq 0$ is an integer, $F_1,\ldots,F_K: X \to \R$ are functions with the pointwise bound $|F_k| \leq \nu+1$ for all $k \in [K]$, $\Y_K$ is the factor $\Y_K := \Y_K(F_1) \vee \ldots \vee \Y_K(F_K)$, and $\Omega_K \subset X$ is a $\Y_K$-measurable set.  We refer to $K$ as the \emph{order} of the structured factor, and $\Omega_K$ as the \emph{exceptional set}.  We say that the structured factor has \emph{noise level} $\sigma$ for some $\sigma > 0$ if we have the smallness bound
$$ \int_X 1_{\Omega_K} (\nu+1) \leq \sigma$$
and the pointwise bound
\begin{equation}\label{epoint}
 |(1 - 1_{\Omega_K}) \E( \nu-1|\Y_K )| \leq \sigma.
 \end{equation}
If $g: X \to \R$ is the function in Theorem \ref{abstract}, we define the \emph{energy} ${\mathcal E}_g(\Y_K)$ of the structured factor $Y$ relative to $g$ to be the quantity
$$ {\mathcal E}_g(\Y_K) := \int_X (1-1_{\Omega_K}) \E(g|\Y_K)^2.$$
\end{definition}

If $\Y_K$ has noise level $\sigma \leq 1$, then since $g$ is bounded in magnitude by $\nu$, then
observe that
\begin{equation}\label{gyk}
 |(1-1_{\Omega_K}) \E(g|\Y_K)| \leq (1-1_{\Omega_K}) (\E(\nu-1|\Y_K)+1) \leq 1+\sigma \leq 2
\end{equation}
and so we conclude the energy bound
\begin{equation}\label{egyk}
 0 \leq {\mathcal E}_g(\Y_K) \leq 4.
\end{equation}
This will allow us to apply an \emph{energy increment
argument} to obtain Theorem \ref{abstract}.  More precisely, Theorem \ref{abstract} is obtained from the following inductive step.

\begin{proposition}[Inductive step]\label{induct}  Let $\Y_K = (\Y_K, K, F_1, \ldots, F_K, \Omega_K)$ be a structured factor of order $K$ and noise level $0 < \sigma < \eta_4^4$.
Then, if we
\begin{equation}\label{fka}
 F_{K+1} := \frac{1}{1+\sigma} (1 - 1_{\Omega_K}) (g - \E(g|\Y))
 \end{equation}
and we suppose that
\begin{equation}\label{fka2}
 |\int_X F_{K+1} \tilde \D F_{K+1}| > \eta_4
\end{equation} 
then there exists a structured factor $\Y_{K+1} = (\Y_{K+1}, K+1, F_1, \ldots, F_K, F_{K+1},\Omega_{K+1})$ of order $K+1$ and noise level $\sigma + O_{\eta_4}(\eta_5^{1/2})$ with the energy increment property
\begin{equation}\label{energy-inc}
 {\mathcal E}_g(\Y_{K+1}) \geq {\mathcal E}_g(\Y_{K}) + c \eta_4^2
\end{equation}
for some constant $c > 0$ (depending only on $I$).
\end{proposition}

Let us assume Proposition \ref{induct} for the moment and deduce Theorem \ref{abstract}.  
Starting with a trivial structured factor $\Y_0$ of order $0$ and quality $0$ and iterating 
Proposition \ref{induct} repeatedly (and using \eqref{egyk} to prevent the iteration for proceeding
for more than $4/c\eta_4^2 = O_{\eta_4}(1)$ steps), we may find a structured factor $\Y_K$ of order $K = O_{\eta_4}(1)$ and
noise level 
\begin{equation}\label{sigbound}
\sigma = O_{\eta_4}(\eta_5^{1/2}) < \eta_4^4
\end{equation}
such that the function $F_{K+1}$ defined in \eqref{fka}
obeys the bound
$$ |\int_X F_{K+1} \tilde \D F_{K+1}| \leq \eta_4.$$
If we thus set $g_U := F_{K+1}$ and
$$ g_{U^\perp} := \frac{1}{1+\sigma} (1 - 1_{\Omega_K}) \E(g|\Y)$$
then we easily verify \eqref{gdecomp} and \eqref{gud}, while \eqref{struct-bound} follows from \eqref{epoint}, since $\E(g|\Y) \leq 1 + \E(\nu-1|\Y)$.  To prove \eqref{most-mass},
we see from \eqref{sigbound} that it suffices to show that
$$ \int_X (1 - 1_{\Omega_K}) \E(g|\Y) = \int_X g - O_{\eta_4}(\eta_5^{1/2}).$$
Since $\Omega_K$ is $\Y$-measurable, the left-hand side is $\int_X g - \int_X 1_{\Omega_K} g$.
But the claim then follows from \eqref{epoint} and \eqref{sigbound}.  This proves
Theorem \ref{abstract}.

It remains to prove Proposition \ref{induct}.  
Set $\Y_{K+1} := \Y \vee \Y(\tilde \D F_{K+1}) = \Y(\tilde \D F_1) \vee \ldots \vee \Y(\tilde \D F_{K+1})$. Now from Proposition \ref{nu-uniform} we can find a $\Y_{K+1}$-measurable set $\Omega$
obeying the smallness bound \eqref{om-small} and the pointwise bound
\begin{equation}\label{nup}
 |(1 - 1_\Omega) \E( \nu-1|\Y_{K+1} )| \leq O_{\eta_4}(\eta_5^{1/2}).
\end{equation}
Now set $\Omega_{K+1} := \Omega_K \cup \Omega$.  This is still $\Y_{K+1}$-measurable and
$\int_X \Omega_{K+1} \leq \sigma + O_{\eta_4}(\eta_5^{1/2})$; from \eqref{nup} we thus conclude
that $\Y_{K+1}$ has noise level $\sigma + O_{\eta_4}(\eta_5^{1/2})$.  Thus the only thing left to
verify is the energy increment property \eqref{energy-inc}.

From \eqref{fka}, \eqref{fka2} we have
\begin{equation}\label{omg}
 |\int_X (1 - 1_{\Omega_K}) (g - \E(g|\Y_K)) \tilde \D F_{K+1}| \geq \eta_4 - O(\eta_4^2).
 \end{equation}
Now from \eqref{dfb}, the pointwise bound $0 \leq g \leq \nu$, \eqref{gyk}, and \eqref{om-small} 
we have
\begin{align*}
|\int_X (1_{\Omega_{K+1}} - 1_{\Omega_K}) (g - \E(g|\Y_K)) \tilde \D F_{K+1}|
&\leq O(\int_X (1_{\Omega_{K+1}} - 1_{\Omega_K}) (\nu + 1))  \\
&\leq O_{\eta_4}(\eta_5^{1/2}) = O(\eta_4^2)
\end{align*}
and hence by \eqref{omg}
$$
 |\int_X (1 - 1_{\Omega_{K+1}}) (g - \E(g|\Y_K)) \tilde \D F_{K+1}| \geq \eta_4 - O(\eta_4^2).$$
Next, from \eqref{tide}, the pointwise bound $0 \leq g \leq \nu$ and \eqref{measure} we have
\begin{align*}
 |\int_X (1 - 1_{\Omega_{K+1}}) (g - \E(g|\Y_K)) (\tilde \D F_{K+1} - \E( \tilde \D F_{K+1}|\Y_{K+1} ))|
 &\leq \int_X (\nu + 1) O(\eta_4^2) \\
 &= O(\eta_4^2) 
\end{align*}
and thus
$$
 |\int_X (1 - 1_{\Omega_{K+1}}) (g - \E(g|\Y_K)) \E(\tilde \D F_{K+1}|\Y_{K+1}) | \geq \eta_4 - O(\eta_4^2).$$
Since $\Omega_{K+1}$, $\E(g|\Y_K)$, and $\E(\tilde \D F_{K+1}|\Y_{K+1})$ are already $\Y_{K+1}$-measurable, we conclude
$$
 |\int_X (1 - 1_{\Omega_{K+1}}) (\E(g|\Y_{K+1}) - \E(g|\Y_K)) \E(\tilde \D F_{K+1}|\Y_{K+1}) | \geq \eta_4 - O(\eta_4^2).$$
By \eqref{dfb} and Cauchy-Schwarz we conclude that
\begin{equation}\label{pythag}
 \int_X (1 - 1_{\Omega_{K+1}}) |\E(g|\Y_{K+1}) - \E(g|\Y_K)|^2 \geq 2c\eta_4^2 - O(\eta_4^3)
\end{equation}
for some $c > 0$.

To pass from this to \eqref{energy-inc}, first observe from \eqref{gyk} and \eqref{om-small} that
$$ \int_X (1_{\Omega_{K+1}} - 1_{\Omega_K}) \E(g|\Y_K)^2 \ll_{\eta_4} \eta_5^{1/2}$$
and so by the triangle inequality and \eqref{egyk}, \eqref{energy-inc} will follow from the estimate
$$ \int_X (1 - 1_{\Omega_{K+1}}) \E(g|\Y_{K+1})^2
\geq \int_X (1 - 1_{\Omega_{K+1}}) \E(g|\Y_{K})^2 + 2c \eta_4^2 - O(\eta_4^3).$$
Using the identity
$$ \E(g|\Y_{K+1})^2 = \E(g|\Y_{K})^2 + |\E(g|\Y_{K+1}) - \E(g|\Y_K)|^2
+ 2 \E(g|\Y_K) (\E(g|\Y_{K+1}) - \E(g|\Y_K))$$
and \eqref{pythag}, it will suffice to show that
$$ \int_X (1 - 1_{\Omega_{K+1}}) \E(g|\Y_K) (\E(g|\Y_{K+1}) - \E(g|\Y_K)) \ll \eta_4^3.$$
Now observe that $\E(g|\Y_{K+1}) - \E(g|\Y_K)$ is orthogonal to all $\Y_K$-measurable functions, and in particular
$$ \int_X (1 - 1_{\Omega_{K}}) \E(g|\Y_K) (\E(g|\Y_{K+1}) - \E(g|\Y_K)) = 0.$$
Thus it suffices to show that
$$ \int_X (1_{\Omega_{K+1}} - 1_{\Omega_{K}}) \E(g|\Y_K) (\E(g|\Y_{K+1}) - \E(g|\Y_K)) \ll \eta_4^3.$$
Since everything here is $\Y_{K+1}$-measurable, we may replace $\E(g|\Y_{K+1})$ by $g$.  Using \eqref{gyk} it then suffices to show
$$ \int_X (1_{\Omega_{K+1}} - 1_{\Omega_{K}}) |g - \E(g|\Y_K)| \ll \eta_4^3.$$
But this follows from the pointwise bound $0 \leq g \leq \nu$, from \eqref{gyk}, and \eqref{om-small}.
This concludes the proof of Proposition \ref{induct}, which in turn implies Theorem \ref{abstract} and thus Theorem \ref{struct-thm}.

\endprf

\section{A pseudorandom measure which majorizes the primes}\label{sec8}

In the remainder of the paper we prove Theorem \ref{major}, which constructs the pseudorandom measure $\nu$ which will pointwise dominate the function $f$ defined in \eqref{fdef}.  As in all 
previous sections we are using the notation from Section \ref{notation-sec} to define quantities such
as $W, R, M, b$.

The measure $\nu$ can in fact be described explicitly, following \cite{tao-gy-notes}, \cite{host-survey}, \cite{linprimes}.  Let $\chi: \R \to \R$ is a fixed smooth even function which vanishes outside of the interval $[-1,1]$ and obeys
the normalization 
\begin{equation}\label{chinorm}
\int_0^1 |\chi'(t)|^2\ dt = 1,
\end{equation}
but is otherwise arbitrary\footnote{This differs slightly from the majorant introduced by Goldston and Y{\i}ld{\i}r{\i}m \cite{goldston-yildirim-old1} and used in \cite{gt-primes}; in our current notation, the majorants from those papers corresponds to the case $\chi(t) := \max(1-|t|,0)$.  It turns out that choosing $\chi$ to be smooth allows for some technical simplifications, at the (acceptable) cost of lowering $\eta_2 = \frac{\log R}{\log N}$ slightly.}.  We then define $\nu$ by the formula
\begin{equation}\label{nudef}
\nu(x) = \nu_{\chi}(x) := \frac{\phi(W)}{W} \log R \left(\sum_{m|Wx+b} \mu(m) \chi(\frac{\log m}{\log R})\right)^2
\end{equation}
for $x \in [N]$, where the sum is over all positive integers $m$ which divide $Wx+b$, and $\mu(m)$ is the \emph{M\"obius function} of $m$, defined as $(-1)^k$ when $m$ is the product of $k$ distinct primes for some $k \geq 0$, and zero otherwise (i.e. zero when $m$ is divisible by a non-trivial square).

\begin{remark}
The definition of $\nu$ may seem rather complicated, but its behavior is in fact rather easily controlled, at least at ``coarse scales'' (averaging $x$ over intervals of length greater than a large power of $R$), by sieve theory techniques, and in particular by a method of Goldston and Y{\i}ld{\i}r{\i}m \cite{goldston-yildirim-old1}, though in the paper here we exploit the smoothness of the cutoff $\chi$ (as in \cite{host-survey}, \cite{coates}, \cite{linprimes}) to avoid the need for multiple contour integration, relying on the somewhat simpler Fourier integral expansion instead.  
For instance, at such scales it is known from these methods that the average value of $\nu$ is $1 + o(1)$ (see e.g. \cite{host-survey}, \cite{coates}), and more generally a large family of linear correlations of $\nu$ with itself
are also $1 + o(1)$ (see \cite{gt-primes}, \cite{linprimes}).  Thus one can view $\nu$ as being close to $1$ in a weak (averaged) sense, though of course in a pointwise sense $\nu$ will fluctuate tremendously.
\end{remark}

It is easy to verify the pointwise bound $f(x) \leq \nu(x)$.  Indeed, from \eqref{fdef}, \eqref{nudef} it suffices to verify that
$$ \sum_{m|Wx+b} \mu(m) \chi(\frac{\log m}{\log R}) = 1$$
whenever $x \in [\frac{1}{2} N]$ and $Wx+b \in A$.  But this is clear since $Wx+b$ is prime and greater than $R$.  It is also easy to verify the bound \eqref{nu-point}, using the elementary result that the number of divisors of an integer $n$ is $O_\eps(n^\eps)$ for any $\eps > 0$.

The remaining task is to verify that $\nu$ obeys both the polynomial forms condition \eqref{polyform} and the polynomial correlation condition \eqref{polycor-eq} (note that \eqref{measure} follows from \eqref{polyform}).  We can of course take $N$ to be large compared with the parameters $\eta_0,\ldots,\eta_7$ and with the parameters $D',D'',K,\eps$ (in the case of \eqref{polycor-eq}) as the claim is trivial otherwise.

We begin with a minor reduction designed to eliminate the ``wraparound'' effects caused by working in the cyclic group $X = \Z/N\Z$ rather than the interval $[N]$.  Let us define the truncated domain $X'$ to be the interval $X' := \{ x \in \Z: \sqrt{N} \leq x \leq N - \sqrt{N} \}$ (say).  From \eqref{nu-point} we can replace the average in $X$ with the average in $X'$ in both \eqref{polyform} and \eqref{polycor-eq} while only incurring an error of $o(1)$ or $o_{D',D'',K}(1)$ at worst.  The point of restricting to $X'$ is that all the shifts which occur in \eqref{polyform} and \eqref{polycor-eq} have size at most $O(M^{d_*} W^{O(1/\eta_1)})$ or $O_{D',D'',K}(M^d W^{O(1/\eta_1)})$, because of the hypotheses on the degree and coefficients of the polynomials and because all convex bodies are contained in a ball $B(0,M^2)$.  By choice of $M$, these shifts are thus less than $\sqrt{N}$ and so we do not encounter any wraparound issues.  Thus \eqref{polyform} is now equivalent to
\begin{equation}\label{polyform-2}
 \E_{\vec h \in \Omega \cap \Z^d} \E_{x \in X'} \prod_{j \in [J]} \nu(x + Q_1(\vec h)) = 1 + o_{\eps}(1)
\end{equation}
and \eqref{polycor-eq} is similarly equivalent to
\begin{equation}\label{polycor-eq-2}
\begin{split}
 &\E_{\vec n \in \Omega' \cap \Z^{D'}} \E_{\vec h \in \Omega'' \cap \Z^{D''}} \E_{x \in X'}  \\
&\left[ \prod_{k \in [K]} \E_{\vec m \in \Omega \cap \Z^D} \prod_{j \in [J]} \nu(x + \vec P_j(\vec h) \cdot \vec m + \vec Q_{j,k}(\vec h) \cdot \vec n) \right] \prod_{l \in [L]} \nu( x + \vec S_l(\vec h) \cdot \vec n) \\
&\quad = 1 + o_{D',D'',K,\eps}(1)
\end{split}
\end{equation}
where $\nu$ is now viewed as a function on the integers rather than on $X = \Z/N\Z$, defined by \eqref{nudef}.

We shall prove \eqref{polyform-2} and \eqref{polycor-eq-2} in Section \ref{pfc-sec} and Section \ref{pcc-sec} respectively.  Before we do so, let us first discuss what would happen if
we tried to generalize these averages by considering the more general expression
\begin{equation}\label{test}
 \E_{\vec x \in \Omega \cap \Z^D} \prod_{j \in [J]} \nu( P_j( \vec x ) )
\end{equation}
where $D,J \geq 0$ are integers, $\Omega$ is a convex body in $\R^D$, and $P_1,\ldots,P_J \in \Z[\x_1,\ldots,\x_D]$ are polynomials of bounded degree and whose coefficients are not too large (say of size $O(W^{O(1)})$).  In light of the linear correlation theory, one would generally expect
these polynomial correlations to also be $1+o(1)$ as long as the polynomials $P_1,\ldots,P_J$ were suitably ``distinct'' and that the range $\Omega$ is suitably large.

There will however be some technical issues in establishing such a statement.  For sake of 
exposition let us just discuss the case $J=1$, so that we are averaging a single factor $\nu(P(\vec x))$ for some polynomial $P$ of $D$ variables $\vec x = (x_1,\ldots,x_D)$.  Even in this simple case, two basic problems arise.  

The first problem is that $\nu$ is
not perfectly uniformly distributed modulo $p$ for all primes $p$.  The ``$W$-trick'' of using $Wx+b$ instead of $x$ in \eqref{nudef} (and renormalizing by $\frac{\phi(W)}{W}$ to compensate) does guarantee a satisfactory uniform distribution of $\nu$ modulo $p$ for small primes $p < w$.  However for larger primes $p \geq w$, it turns out that $\nu$ will generally avoid
the residue class $\{x: Wx+b = 0\ \mod\ p\}$ and instead distribute itself uniformly among the other $p-1$ residue classes.  This corresponds to the basic fact that primes (and almost primes) are mostly coprime to any given modulus $p$.  Because of this, the
expected value of an expression such as $\nu(P(x))$ will increase from $1$ to roughly $(1-\frac{1}{p})^{-1}$ if we know that $WP(x)+b$ is coprime to $p$, and conversely it will drop to essentially zero if we know that $WP(x)+b$ is divisible by $p$.  These two effects will essentially balance each other out, provided that the algebraic variety $\{ x \in F_p^D: WP(x)+b = 0 \}$ has the expected density of $\frac{1}{p} + O(\frac{1}{p^{3/2}})$ (say) over the finite field affine space $F_p^D$.  The famous result of Deligne \cite{deligne}, \cite{deligne2}, in which the Weil conjectures were proved, establishes this when $WP+b$ was non-constant and is absolutely irreducible modulo $p$ (i.e. irreducible over the algebraic closure of $F_p$).  However, there can be some ``bad'' primes $p \geq w$ for which this irreducibility fails; a particularly ``terrible'' case arises when $p$ divides the polynomial $WP+b$, in which case the variety has density $1$ in $F_p^D$ and
the expected value of \eqref{test} drops to zero.  This reflects the intuitive fact that $WP(x)+b$ is much less likely to be
prime or almost prime if $WP+b$ itself is divisible by some prime $p$.  The other bad primes $p$ do not cause such a severe change in the expectation \eqref{test}, but can modify the expected answer of $1 + o(1)$ by a factor of $1+O(\frac{1}{p}) = \exp( O(\frac{1}{p}) )$, leading to a final value which is something like $\exp( O( \sum_{p \hbox{ bad}} \frac{1}{p} ) + o(1) )$.
In most cases this expression will be in fact very close to one, because of the restriction $p \geq w$.  However, the
(very slow) divergence of the sum $\sum_p \frac{1}{p}$ means that there are some exceptional cases in which averages such as \eqref{test} are unpleasantly large.  For instance, for any fixed $h \neq 0$, the average value of $\nu(x)\nu(x+h)$ over sufficiently coarse scales turns out to be $\exp( O(\sum_{p \geq w: p|h} \frac{1}{p}) + o(1) )$, which can be arbitrarily large in the (very rare) case that $h$ contains many prime factors larger than $w$, the basic problem being that the algebraic variety $\{ x \in F_p: (Wx+b)(W(x+h)+b) = 0 \}$, which is normally empty, becomes unexpectedly large when $p \geq w$ and $p$ divides $h$.  This phenomenon was already present in \cite{gt-primes}, leading in particular to the rather technical ``correlation condition'' for $\nu$.

The second problem, which is a new feature in the polynomial case compared with the previous linear theory, is that we will not
necessarily be able to average \emph{all} of the parameters $x_1,\ldots,x_D$ over coarse scales (e.g. at scales $O(M)$, $O(\sqrt{M})$ or $O(M^{1/4})$.  Instead, some of the parameters will be only averaged over fine scales such as $O(H)$.  At these scales, the elementary sieve theory methods we are employing cannot estimate the expression \eqref{test} directly; indeed, the problem then becomes analogous to that of understanding the distribution of primes in short intervals, which is notoriously difficult.  Fortunately, we can proceed by first fixing the fine-scale parameters and using the sieve theory methods to compute the averages over the coarse-scale parameters rather precisely, leading to certain tractable divisor sums over ``locally bad primes'' which can then be averaged over fine scales.  Here we will rely on a basic heuristic from algebraic 
geometry, which asserts that a ``generic'' slice of an algebraic variety by a linear subspace will have the same codimension as the original variety.  In our context, this means that a prime which is ``globally good'' with respect to many parameters, will also be ``locally good'' when freezing one or more parameters, for ``most'' choices of such parameters.  We will phrase the precise versions of these statements as a kind of ``combinatorial Nullstellensatz'' (cf. \cite{alon}) in Appendix \ref{variety-sec}.
This effect lets us deal with the previous difficulty that the sum of $\frac{1}{p}$ over bad primes can occasionally be very large.

We have already mentioned the need to control the density of varieties such as $\{ x \in F_p^D: WP(x)+b = 0 \}$, which in general requires the Weil conjectures as proven by Deligne.  Fortunately, for the application to polynomial progressions, the polynomials $P$ involved will always be linear in at least one of the coarse-scale variables.  This makes the density of the algebraic variety much easier to compute, provided that the coefficients in this linear representation do not degenerate (either by the linear coefficient vanishing, or by the linear and constant coefficients sharing a common factor).  
Thus we are able to avoid using the Weil conjectures. In fact we will be able to proceed by rather elementary algebraic methods, without using modern tools from arithmetic geometry; see Appendix \ref{variety-sec}.

\subsection{Notation}\label{notation}

We now set out some notation which will be used throughout the proof of \eqref{polyform-2} and \eqref{polycor-eq-2}. If $p$ is a prime, we use $F_p$ to denote the finite field of $p$ elements.

If $P,Q$ lie in some ring $R$, we use $P|Q$ to denote the statement that $Q$ is a multiple of $P$.
An element of a ring is a \emph{unit} if it is invertible, and \emph{irreducible}\footnote{We shall reserve the term \emph{prime} for the rational primes $2,3,5,7,\ldots$ to avoid confusion.} if it is not a unit, and cannot be written as the product of two non-units.  A ring is a \emph{unique factorization domain} if every element is uniquely expressible as a finite product of irreducibles, up to permutations and units.  
If $P_1,\ldots,P_J$ lie in a unique factorization domain, we say that $P_1,\ldots,P_J$ are \emph{jointly coprime} (or just \emph{coprime} if $J=2$) if there exists no irreducible which divides all the $P_1,\ldots,P_J$, and \emph{pairwise coprime} if each pair $P_i, P_j$ is coprime for $1 \leq i < j \leq J$; thus pairwise coprime implies jointly coprime, but not conversely.

As observed by Hilbert, if $R$ is a unique factorization domain, then so is $R[\x]$ (thanks to the Euclidean algorithm).  In particular, $F_p[\x_1,\ldots,\x_D]$ and $\Z[\x_1,\ldots,\x_D]$ are unique factorization domains (with units
$F_p \backslash \{0\}$ and $\{-1,+1\}$ respectively).  

Every polynomial in $R[\x_1,\ldots,\x_D]$ can of course be viewed as a function from $R^D$ to $R$.
If $P \in \Z[\x_1,\ldots,\x_D]$ is a polynomial and $N \geq 1$, we write $P\ \mod\ N$ for the associated polynomial
in $\Z_N[\x_1,\ldots,\x_D]$ formed by projecting all the coefficients onto the ring $\Z_N$, thus $P\ \mod\ N$ can be viewed as a function from $\Z_N^D$ to $\Z_N$.   Note that this projection may alter the property of two or more polynomials being jointly or pairwise coprime; the precise analysis of when this occurs will in fact be a major focus of our arguments here.

It will be convenient to introduce the modified exponential function
$$ \EXP(x) := \max( e^x - 1, 0 )$$
thus $\EXP(x) \sim x$ when $x$ is non-negative and small, while $\EXP(X) \sim e^x$ for $x$ large.  Observe the elementary
inequalities
\begin{equation}\label{exp-bound}
 \EXP(x+y) \leq \EXP(2x) + \EXP(2y); \quad \EXP(x)^K \ll_K \EXP(Kx) 
\end{equation}
for any $x,y \geq 0$ and $K \geq 1$.

\section{Local estimates}

Before we estimate correlation estimates for $\nu$ on the integers, we first need to consider the analogous problem modulo $p$.  To formalize this problem we introduce the following definition.

\begin{definition}[Local factor]\label{local-def}
Let $P_1,\ldots,P_J: \Z[\x_1,\ldots,\x_D]$ be polynomials with integer coefficients.  For any prime $p$, we define
the (principal) \emph{local factor}
$$ c_p(P_1,\ldots,P_J) := \E_{x \in F_p^D} \prod_{j \in [J]} 1_{P_j(x)=0\ \mod \ p}.$$
We also define the \emph{complementary local factor}
$$ \overline{c_p}(P_1,\ldots,P_J) := \E_{x \in F_p^D} \prod_{j \in [J]} 1_{P_j(x) \neq 0\ \mod \ p}.$$
\end{definition}

\begin{examples} If $P_1,\ldots,P_J$ are homogeneous linear forms on $F_p^D$, with total rank $r$, then $c_p(P_1,\ldots,P_J) = p^{-r}$.  If the forms are independent (thus $J=r$), then $\overline{c_p}(P_1,\ldots,P_J) = (1 - \frac{1}{p})^J$.  If $D=1$, then the local factor $c_p( \x^2 + 1 )$ equals $2/p$ when $p = 1\ \mod\ 4$ and equals $0$ when $p = 3\ \mod 4$, by quadratic reciprocity.  (When $p=2$, it is equal to $1/p$.)   More generally, the Artin reciprocity law \cite{artin} relates Artin characters to certain local factors.  Deligne's celebrated proof \cite{deligne}, \cite{deligne2} of the Weil conjectures implies (as a very special case) that $c_p(P) = 1/p + O_{k,D}(1/p^{3/2})$ whenever $P \in \Z[\x_1,\ldots,\x_D]$ determines a non-singular projective algebraic variety over $F_p$.  For instance, if $P = \x_2^2 - \x_1^3 - a \x_1 - b$, so that $P$ determines an elliptic curve,
with discriminant $\Delta= -16(4a^3+27b^2)$ coprime to $p$, then $c_p(P) = 1/p + O(1/p^{3/2})$ (a classical result of Hasse).  The Birch and Swinnerton-Dyer conjectures, if true, would provide more precise information (though not of upper bound type) on the error term in this case.
\end{examples}

\begin{remark}
The factor $c_p$ denotes the proportion of points on $F_p^D$ which lie on the algebraic variety determined by the polynomials $P_1,\ldots,P_J$, while the complementary factor $\overline{c_p}$ is the proportion of points in $F_p^D$ for which all the $P_1,\ldots,P_J$ are coprime to $p$.  Clearly these factors lie between $0$ and $1$; for instance when $J=0$ we have $c_p = 1$ and $\overline{c_p} = 0$. Our interest is to estimate $c_p$ for higher values of $J$.  This will be of importance when we come to the ``global'' estimates for $\prod_{j \in [J]} \nu(P_j(x))$ over various subsets of $\Z^d$; heuristically, the average value of this expression should be approximately the product of the complementary factors $\overline{c_p}$ as $p$ ranges over primes.
\end{remark}

From the inclusion-exclusion principle we have the identity
\begin{equation}\label{exclude}
\overline{c_p}(P_1,\ldots,P_J) = \sum_{S \subseteq \{1,\ldots,J\}} (-1)^{|S|}
c_p( (P_j)_{j \in S} )
\end{equation}
and so we can estimate the complementary local factors using the principal local factors.

As mentioned earlier, the precise estimation of $c_p(P_1,\ldots,P_J)$ for general $P_1, \ldots, P_J$ is intimately connected to a number of deep results in arithmetic geometry such as the Weil conjectures and the Artin reciprocity law.  Fortunately, for our applications we
will only need to know the $\frac{1}{p}$ coefficient of $c_p(P_1,\ldots,P_J)$ and can neglect lower order terms.  Also, we will be working in the case where each of the polynomials $P_j$ are linear in at least one of the co-ordinates $x_1,\ldots,x_D$ of $x$ and are ``non-degenerate'' in the other co-ordinates.  In such a simplified context, we will be able to control $c_p$ quite satisfactorily using only arguments from elementary algebra.  To state the results, we first need the notion of a prime $p$ being good, bad, or terrible with respect to a collection of polynomials:

\begin{definition}[Good prime]\label{bad-def} Let $P_1,\ldots,P_J \in \Z[\x_1,\ldots,\x_D]$ be a collection of polynomials.  We say that a prime $p$ is \emph{good} with respect to $P_1,\ldots,P_J$ if the following hold:
\begin{itemize}
\item The polynomials $P_1\ \mod\ p,\ldots,P_J\ \mod\ p$ are pairwise coprime.
\item For each $1 \leq j \leq J$, there exists a co-ordinate $1 \leq i_j \leq D$ for which we have the linear behavior
\begin{align*}
 P_j(\x_1,\ldots,\x_D) &= P_{j,1}(\x_1,\ldots,\x_{i_j-1}, \x_{i_j+1},\ldots, \x_{D}) \x_{i_j} \\
 &\quad + P_{j,0}(\x_1,\ldots,\x_{i_j-1},\x_{i_j+1},\ldots,\x_{D})\ \mod\ p 
\end{align*}
where $P_{j,1}, P_{j,0} \in F_p[\x_1,\ldots,\x_{i_j-1},\x_{i_j+1},\ldots,\x_D]$ are such that $P_{j,1}$ is non-zero and coprime to $P_{j,0}$.
\end{itemize}
We say that a prime is \emph{bad} if it is not good.  We say that a prime is \emph{terrible} if at least one of the $P_j$ vanish identically modulo $p$ (i.e. all the coefficients are divisible by $p$).  Note that terrible primes are automatically bad.
\end{definition}

Our main estimate on the local factors is then as follows.

\begin{lemma}[Local estimates]\label{local-cor}  Let $P_1,\ldots,P_J \in \Z[\x_1,\ldots,\x_D]$ have degree at most $d$, let $p$ be a prime, and let $S \subset \{1,\ldots,J\}$.  Then:
\begin{itemize}
\item[(a)] If $|S| = 0$, then $c_p( (P_j)_{j \in S} ) = 1$.
\item[(b)] If $|S| \geq 1$ and $p$ is not terrible, then $c_p( (P_j)_{j \in S} ) = O_{d,D,J}(\frac{1}{p})$.
\item[(c)] If $|S| = 1$ and $p$ is good, then $c_p( (P_j)_{j \in S} ) = \frac{1}{p} + O_{d,D,J}(\frac{1}{p^2})$.
\item[(d)] If $|S| > 1$ and $p$ is good, then $c_p( (P_j)_{j \in S} ) = O_{d,D,J}(\frac{1}{p^2})$.
\item[(e)] If $p$ is terrible, then $\overline{c_p}(P_1,\ldots,P_J) = 0$.
\item[(f)] If $p$ is not terrible, then $\overline{c_p}(P_1,\ldots,P_J) = 1 + O_{d,D,J}(\frac{1}{p})$.
\end{itemize}
\end{lemma}

The proof of this lemma involves only elementary algebra, but we defer it to Appendix \ref{variety-sec} so as not to disrupt the flow of the argument.

\begin{remark} From \eqref{exclude} and Lemma \ref{local-cor}(acd) we also have
$$ \overline{c_p}(P_1,\ldots,P_J) = 1 - \frac{J}{p} + O_{d,D,J}(\frac{1}{p^2})$$
when $p$ is good.  In practice we shall need a more sophisticated version of this fact, when certain complex weights $p^{-\sum_{j \in S} z_j}$ are inserted into the right-hand side of \eqref{exclude}; see Lemma \ref{euler}.
\end{remark}

\section{Initial correlation estimate}\label{init-sec}

To prove \eqref{polyform-2} and \eqref{polycor-eq-2} we shall need the following initial estimate which handles general polynomial averages of $\nu$ over large scales, but with an error term that can get large if there are many ``bad'' primes present.  More precisely, this section is devoted to proving

\begin{proposition}[Correlation estimate]\label{corrprop} Let $P_1,\ldots,P_J \in \Z[\x_1,\ldots,\x_D]$ have degree at most $d$ for some $J,D,d \geq 0$.  Let $\Omega$ be a convex body in $\R^D$ of inradius at least $R^{4J+1}$. 
Let $\primes_b$ be the set of primes $w \leq p \leq R^{\log R}$ which are bad with respect to $WP_1+b,\ldots,WP_J+b$, and let $\primes_t \subset \primes_b$ be the set of primes $w \leq p \leq R^{\log R}$ which are terrible (as defined in Definition \ref{bad-def}).
Then 
\begin{equation}\label{numess}
 \E_{x \in \Omega \cap \Z^D} \prod_{j \in [J]} \nu(P_j(x)) = 1_{\primes_t = \emptyset} + o_{D,J,d}(1) + 
 O_{D,J,d}( \EXP( O_{D,J,d}( \sum_{p \in \primes_b} \frac{1}{p} ) ) ).
 \end{equation}
\end{proposition}

\begin{remark} We only expect this estimate to be useful when the number of bad primes is finite.  This is equivalent to requiring that the polynomials $P_1,\ldots,P_J$ are coprime, and each one is linear in at least one variable.  Because the sum $\sum_p \frac{1}{p}$ is (very slowly) divergent (see \eqref{log-chebyshev}), the last error term can be unpleasantly 
large on occasion, but in practice we will be able to introducing averaging over additional parameters which will make the effect of the error small on average, the point being that the sets $\primes_t, \primes_b$ are generically rather small.  The radius $R^{4J+1}$ is not best possible, but to lower it too much would require some deep analytical number theory estimates such as the Bombieri-Vinogradov inequality, which we shall avoid using here.  The upper bound $R^{\log R}$ (which was not present in earlier work) can also be lowered, but for our purposes any bound which is subexponential in $R$ will suffice.
\end{remark}

\begin{remark} All the primes $p < w$ will be bad (but not terrible); however their contribution will be almost exactly canceled by the $\frac{\phi(W)}{W}$ term present in $\nu$ and we do not need to include them into $\primes_b$.   Even a single terrible prime will cause the main term $1_{\primes_t = \emptyset}$ to vanish (basically because one of the $P_j(x)$ will now be inherently composite and so will unlikely to have a large value of $\nu$), which will make asymptotics difficult; however, terrible primes are no worse than merely bad primes for the purposes of \emph{upper} bounds.
\end{remark}

\begin{proof}[Proof of Proposition \ref{corrprop}]  Throughout this proof we fix $D$, $J$, $d$, and allow the implied constants in the $O()$ and $o()$ notation to depend on these parameters.  We will also always assume $R$ to be sufficiently large depending on $D,J,d$.

We expand out the left-hand side using \eqref{nudef} as
\begin{equation}\label{disc0}
\begin{aligned}
\left[\frac{\phi(W)}{W} \log R\right]^J \sum_{m_1,m'_1,\ldots,m_J,m'_J \geq 1}
&\left(\prod_{j \in [J]} \mu(m_j) \mu(m'_j) \chi(\frac{\log m_j}{\log R}) \chi(\frac{\log m'_j}{\log R}) \right) \\
&\E_{x \in \Omega \cap \Z^D} \prod_{j = 1}^J 1_{\lcm(m_j,m'_j) | WP_j(x)+b}.
\end{aligned}
\end{equation}
Here of course $\lcm()$ denotes least common multiple.  Note that the presence of the $\mu$ and $\chi$ factors allows us to restrict $m_1,\ldots,m'_J$ to be square-free and at most $R$.

The first task is to eliminate the role of the convex body $\Omega$, taking advantage of the large inradius assumption.
Let $M := \lcm(m_1,m'_1,\ldots,m_j,m'_j)$, thus $M$ is square-free and at most $R^{2J}$.  The function
$x \mapsto 1_{\lcm(m_j,m'_j) | WP_j(x)+b}$ is periodic with respect to the lattice $M \cdot \Z^D$, and thus can be meaningfully defined on the group $\Z_M^D$.  Applying Corollary \ref{avg} (recalling that $\Omega$ is assumed to have inradius at least $R^{4J+1}$), we thus have
$$ \begin{aligned}
&\E_{x \in \Omega \cap \Z^D} \prod_{j = 1}^J 1_{\lcm(m_j,m'_j) | WP_j(x)+b} \\
&= \left(1 + O(\frac{1}{R^{2J+1}}) \right)
\E_{y \in \Z_M^D} \prod_{j = 1}^J 1_{\lcm(m_j,m'_j) | WP_j(y)+b}. \end{aligned}$$
Let us first dispose of the error term $O(\frac{1}{R^{2J+1}})$.  The contribution of this term to \eqref{disc0} can be crudely
bounded by $O(R^{-2J-1})$, and so the contribution of this term to \eqref{disc0} can be crudely bounded by
$$ O\left( \left[\frac{\phi(W)}{W} \log R\right]^J \sum_{1 \leq m_1,m'_1,\ldots,m_J,m'_J \leq R} R^{-2J-1} \right) \ll \frac{\log^J R}{R}  = o(1).$$
Thus we may discard this error,
and reduce to showing that
\begin{equation}\label{disc2}
\begin{aligned}
&\left[\frac{\phi(W)}{W} \log R\right]^J \sum_{m_1,m'_1,\ldots,m_J,m'_J \geq 1}\\
& \quad \left(\prod_{j \in [J]} \mu(m_j) \mu(m'_j) \chi(\frac{\log m_j}{\log R}) \chi(\frac{\log m'_j}{\log R}) \right)
    \alpha_{\lcm(m_1,m'_1),\ldots,\lcm(m_J,m'_J)} \\
& =1_{\primes_t = \emptyset} + o( 1 ) + O\left( \EXP\left( O( \sum_{p \in \primes_b} \frac{1}{p} ) \right) \right)
\end{aligned}
\end{equation}
where $\alpha_{\lcm(m_1,m'_1),\ldots,\lcm(m_J,m'_J)}$ is the local factor
$$\alpha_{\lcm(m_1,m'_1),\ldots,\lcm(m_J,m'_J)} := \E_{y \in \Z_M^D} \prod_{j = 1}^J 1_{\lcm(m_j,m'_j) | WP_j(y)+b}.$$
Observe from the Chinese remainder theorem that $\alpha$ is multiplicative, so that if $\lcm(m_j,m'_j) = \prod_p p^{r_j}$ then
$$ \alpha_{\lcm(m_1,m'_1),\ldots,\lcm(m_J,m'_J)} = \prod_p \alpha_{p^{r_1}, \ldots, p^{r_J}}$$
(note that all but finitely many of the terms in the product are $1$.  If the $m_1,\ldots,m'_J$ are squarefree then
the $r_j$ are either zero or one, and we simplify further to
$$ \alpha_{\lcm(m_1,m'_1),\ldots,\lcm(m_J,m'_J)} = \prod_p c_p( ( WP_j + b )_{r_j = 1} )$$
where the local factors $c_p$ are defined in Definition \ref{local-def} in the appendices, and the dummy variable $j$ is ranging over all indices for which $r_j = 1$.  Also note that the $m_j, m'_j$ are bounded by $R$, we may certainly restrict the primes $p$ to be less than $R^{\log R}$ without difficulty.

The next step is to replace the $\chi$ factors by terms which are multiplicative in the $m_j,m'_j$.
Since $\chi$ is smooth and compactly supported, we have the Fourier expansion 
\begin{equation}\label{fourier}
 e^x \chi(x) = \int_{-\infty}^{\infty} \varphi(\xi) e^{-i x \xi}\ d\xi
\end{equation}
for some smooth, rapidly decreasing function $\varphi(\xi)$ (so in particular $\varphi(\xi) = O_A( (1+|\xi|)^{-A})$ for any $A > 0$).  For future reference, we observe that \eqref{fourier} and the hypotheses on $\chi$ will imply the identity
\begin{equation}\label{varphi-ident}
\int_{-\infty}^\infty \int_{-\infty}^\infty \frac{(1+it) (1+it')}{2+it+it'} \varphi(t) \varphi(t')\ dt dt' = 1
\end{equation}
(see \cite[Lemma D.2]{linprimes}, or the proof of \cite[Proposition 2.2]{coates}).  

We follow the arguments in \cite{tao-gy-notes}, \cite{host-survey}, \cite{coates}, \cite{linprimes}, except that for
technical reasons (having to do with the terrible primes) we will be unable to truncate the $\xi$ variables.
From \eqref{fourier} we have
\begin{align*}
\chi\left( \frac{\log m_j}{\log R} \right) &= \int_{-\infty}^\infty m_j^{-z_j}\ d\xi_j \\
\chi\left( \frac{\log m'_j}{\log R} \right) &= \int_{-\infty}^\infty (m'_j)^{-z'_j}\ d\xi'_j
\end{align*}
where we adopt the notational conventions
$$ z_j := (1+\xi_j)/\log R; \quad z'_j := (1+\xi'_j)/\log R.$$
Our task is thus to show that
\begin{equation}\label{disc3}
\begin{aligned}
&\left[\frac{\phi(W)}{W} \log R\right]^J \sum_{m_1,m'_1,\ldots,m_J,m'_J \geq 1} 
\int_{-\infty}^\infty \dots \int_{-\infty}^\infty \\
&\quad \left[\prod_{j \in [J]} \mu(m_j) \mu(m'_j) m_j^{-z_j} (m'_j)^{-z'_j} \varphi(\xi_j) \varphi(\xi'_j)\ d\xi_j d\xi'_j\right] 
\prod_{p \leq R^{\log R}} c_p( ( WP_j + b )_{r_j = 1} ) \\
&\quad\quad = 1 + o( 1 ) + O\left( \EXP\left( O( \sum_{p \in \primes_b} \frac{1}{p} ) \right) \right).
\end{aligned}
\end{equation}
The left-hand side can be factorized as
$$
\int_{-\infty}^\infty \dots \int_{-\infty}^\infty \left[\frac{\phi(W)}{W} \log R\right]^J \prod_{p \leq R^{\log R}} E_p \prod_{j \in [J]} \varphi(\xi_j) \varphi(\xi'_j)\ d\xi_j d\xi'_j$$
where $E_p = E_p(z_1,\ldots,z'_J)$ is the Euler factor
$$ E_p := \sum_{m_1,\ldots,m'_J \in \{1,p\}}
\prod_{j \in [J]} \mu(m_j) \mu(m'_j) m_j^{-z_j} (m'_j)^{-z'_j} c_p( ( WP_j + b )_{r_j = 1} ).$$
Note that if the $z_j,z'_j$ were zero, then this would just be the complementary factor
$\overline{c_p}( WP_1+b,\ldots,WP_J+b )$ defined in Definition \ref{local-def}; see \eqref{exclude}.  Of course,
$z_j, z'_j$ are non-zero.  To approximate $E_p$ in this case we introduce the Euler factor
$$ E'_p := \prod_j \frac{ (1 - 1/p^{1+z_j}) (1 - 1/p^{1+z'_j})}{1 - 1/p^{1+z_j+z'_j}}.$$
Note that $E'_p$ never vanishes.

\begin{lemma}[Euler product estimate]\label{euler} We have
$$ \prod_{p \leq R^{\log R}} \frac{E_p}{E'_p} = \left(1_{\primes_t = \emptyset} + o(1) + O\left( \EXP\left(O(\sum_{p \in \primes_b} \frac{1}{p})\right) \right) \right) \left(\frac{W}{\phi(W)}\right)^J.$$
\end{lemma}

\begin{proof}  For $p < w$, we directly compute (since $w$ is slowly growing compared to $R$, and $WP_j+b$ is equal modulo $p$ to $b$, which is coprime to $p$) that
$$ E_p = 1 + o(1) \hbox{ and } E'_p = (1 - \frac{1}{p})^J + o(1)$$
and hence (again because $w$ is slowly growing)
$$ \prod_{p<w} \frac{E_p}{E'_p} = (1 + o(1)) \left(\frac{W}{\phi(W)}\right)^J.$$
Thus it will suffice to show that
$$ \prod_{w \leq p \leq R^{\log R}} \frac{E_p}{E'_p} = 1_{\primes_t = \emptyset} + o(1) + O\left( \EXP\left(O(\sum_{p \in \primes_b} \frac{1}{p})\right) \right).$$
For $p$ terrible, Lemma \ref{local-cor} gives the estimate
$$ E_p \ll \frac{1}{p} \ll \frac{1}{p} E'_p$$
and so it will suffice to show
$$ \prod_{w \leq p \leq R^{\log R}, \hbox{ not terrible}} \frac{E_p}{E'_p} = 1 + o(1) + O\left( \EXP\left(O(\sum_{p \in \primes_b} \frac{1}{p})\right) \right) .$$
For $p$ bad but not terrible, Lemma \ref{local-cor} gives the crude estimate
$$ E_p = 1 + O(1/p) = \exp( O(1/p) ) E'_p$$
and thus
$$ \prod_{w \leq p \leq R^{\log R}, \hbox{ bad but not terrible}} \frac{E_p}{E'_p} = 1 + O\left( \EXP\left(O(\sum_{p \in \primes_b} \frac{1}{p})) \right) \right).$$
Thus it suffices to show that
$$ \prod_{w \leq p \leq R^{\log R}, \hbox{ good}} \frac{E_p}{E'_p} = 1 + o(1).$$
Since the product $\prod_p (1 + O(\frac{1}{p^2}))$ is convergent, and $w$ goes to infinity, it in fact suffices to show that
$$ E_p = (1 + O(\frac{1}{p^2})) E'_p$$
for all good primes larger than $w$.  But this easily follows from Lemma \ref{local-cor} and Taylor expansion (recall that  the real parts of $z_j, z'_j$ are $1/\log R > 0$).
\end{proof}

Now we use the theory of the Riemann zeta function. From \eqref{euler-formula-2} we have
$$  \prod_{p \leq R^{\log R}} E'_p = (1 + o(1)) \frac{ \zeta(1 + z_j + z'_j) }{\zeta(1+z_j) \zeta(1 + z'_j) }.$$
On the other hand, from \eqref{zeta-bounds} we have
$$ \frac{1}{\zeta(1+(1+i\xi)/\log R)} = (1 + o( (1 + |\xi|)^2 ) ) \frac{1+i\xi}{\log R} $$
and
$$ \zeta(1+(1+i\xi)/\log R) = (1 + o( (1 + |\xi|)^2 ) ) \frac{\log R}{1+i\xi}$$
for any real $\xi$, and hence
$$ \prod_{p \leq R^{\log R}} E'_p =
(1 + o( (1 + |\xi_j| + |\xi'_j|)^{6} ) )
\frac{1}{\log R} \frac{(1 + i \xi_j) (1 + i \xi'_j)}{2 + i \xi_j + i \xi'_j}.$$
Applying Lemma \ref{euler}, we conclude
$$ \begin{aligned}
&\left[\frac{\phi(W)}{W} \log R\right]^J \prod_{p \leq R^{\log R}} E_p = \\
&\left(1_{\primes_t = \emptyset} + o\left(\prod_{j \in [J]} (1 + |\xi_j| + |\xi'_j|)^{6} \right) + O( \EXP\left(O\left(\sum_{p \in \primes_b} \frac{1}{p})\right) \right)  \right) \prod_{j \in [J]} \frac{(1 + i \xi_j) (1 + i \xi'_j)}{2 + i \xi_j + i \xi'_j}.
\end{aligned}$$
Thus by the triangle inequality, to show \eqref{disc3} it will suffice to show that
$$
\int_\R \dots \int_\R \prod_{j \in [J]} \varphi(\xi_j) \varphi(\xi'_j) \frac{(1 + i \xi_j) (1 + i \xi'_j)}{2 + i \xi_j + i \xi'_j}\ d\xi_j d\xi'_j = 1$$
and
$$
\int_I \dots \int_I \prod_{j \in [J]} |\varphi(\xi_j)| |\varphi(\xi'_j)| (1 + |\xi_j| + |\xi'_j|)^6 \frac{|1 + i \xi_j| |1 + i \xi'_j|}{|2 + i \xi_j + i \xi'_j|}\ d\xi_j d\xi'_j = O(1).$$
But the first estimate follows from \eqref{varphi-ident}, while the second estimate follows from the rapid decrease of $\varphi$.  This proves Proposition \ref{corrprop}.
\end{proof}

To illustrate the above proposition, let us specialize to the case of monic linear polynomials of one variable (this case was essentially treated in \cite{gt-primes} or \cite{goldston-yildirim-old1}).

\begin{corollary}[Correlation condition]\label{corr} Let $h_1,\ldots,h_J$ be integers, and let $I \subset \R$ be an interval of length at least $R^{4J+1}$.  Then
$$
 \E_{x \in I \cap \Z} \prod_{j \in [J]} \nu(x+h_j) = 1 + o_{D,J,d}(1) + O_{D,J,d}\left( \EXP\left( O_{D,J,d}( \sum_{p \in \primes_b} \frac{1}{p} ) \right) \right)
$$
where 
$$\primes_b := \{w \leq p \leq R^{\log R}: p | h_j - h_{j'} \hbox{ for some } 1 \leq j < j' \leq J\}.$$
\end{corollary}

\begin{proof} Apply Proposition \ref{corrprop} with $P_j(x) := x+h_j$ and $\Omega := I$.  Then there are no terrible primes, and the only bad primes larger than $w$ are those which divide $h_j - h_{j'}$ for some $1 \leq j < j' \leq J$.
\end{proof}

This can already be used to derive the ``correlation condition'' in \cite{gt-primes}; a similar aplication of Proposition \ref{corrprop} also gives the ``linear forms condition'' from that paper. We will also need the following variant of the above estimate:

\begin{corollary}[Correlation condition on progressions]\label{corr2} Let $h_1,\ldots,h_J$ be integers, let $q \geq 1$, let $a \in \Z_q$, and let $I \subset \R$ be an interval of length at least $q R^{4J+1}$. Then
$$
 \E_{x \in I \cap \Z; x = a\ \mod\ q} \prod_{j \in [J]} \nu(x+h_j) = 
 O_{D,J,d}\left( \exp\left( O_{D,J,d}( \sum_{p \in \primes_b} \frac{1}{p} ) \right) \right) $$
where 
\begin{equation}\label{pbdef}
 \primes_b := \{w \leq p \leq R^{\log R}: p | h_j - h_{j'} \hbox{ for some } 1 \leq j < j' \leq J\}
\cup \{ p \geq w: p | q \}.
\end{equation}
\end{corollary}

\begin{proof} Apply Proposition \ref{corrprop} with $P_j(x) := (qx+a)+h_j$ and with $\Omega := \{ x \in \R: qx+a \in I\}$.  Then the bad primes are those which divide $h_j - h_{j'}$ or which divide $q$.  (There are terrible primes if $a$ and $q$ are not coprime, but this will not affect the upper bound.  One can get more precise estimates as in Corollary \ref{corr}, but we will not need them here.)
\end{proof}

This in turn implies

\begin{corollary}[Correlation condition with periodic weight]\label{corr3} Let $h_1,\ldots,h_J$ be integers, let $q \geq 1$, let $I \subset \R$ be an interval of length at least $q R^{4J+1}$, and let $f: \Z \to \R^+$ be periodic modulo $q$ (and thus definable on $\Z_q$. Then
$$
 \E_{x \in I \cap \Z} f(x) \prod_{j \in [J]} \nu(x+h_j) =  O_{D,J,d}\left( (\E_{y \in \Z_q} f(y)) 
 \exp\left( O_{D,J,d}( \sum_{p \in \primes_b} \frac{1}{p} ) \right) \right) $$
where $\primes_b$ was defined in \eqref{pbdef}.
\end{corollary}

\begin{proof}
The left-hand side can be bounded by
$$ \E_{y \in \Z_q} f(y) \E_{x \in I \cap \Z; x = a \ \mod\ q} \prod_{j \in [J]} \nu(x+h_j) $$
simply because the set $\{I \cap \Z; x = a\ \mod\ q\}$ has cardinality roughly $\frac{1}{q}$ that of $I \cap \Z$ (by the hypotheses on the length of $I$).  The claim then follows from Corollary \ref{corr2}.
\end{proof}
 
\section{The polynomial forms condition}\label{pfc-sec}

In this section we use the above correlation estimates to prove the polynomial forms condition \eqref{polyform-2}.  We begin with a preliminary bound in this direction.

\begin{theorem}[Polynomial forms condition]\label{pfc} Let $M,D,d,J \geq 0$ and $\eps > 0$.  Let $P_1,\ldots,P_J \in \Z[\m_1,\ldots,\m_D]$ be polynomials of degree $d$ with all coefficients of size at most $W$.  Let $I \subset \R$ be any interval of length at least $R^{4J+1}$, and let $\Omega \subset \R^D$ be any convex body of inradius at least $R^\eps$.  Then
\[\begin{aligned} 
&\E_{x \in I \cap \Z; \vec m \in \Omega \cap \Z^D} \prod_{j \in [J]} \nu(x + P_j(\vec m))\\
&= 1 + o_{D,d,\eps,J}\left(\exp\left( O_{D,d,\eps,J}(\sum_{p \in \primes_b} \frac{1}{p} ) \right)\right)
\end{aligned}\]
where $\primes_b$ denote the set of all $w \leq p \leq \R^{\log R}$ which are ``globally bad'' in the sense that
$p | P_j-P_{j'}$ for some $1 \leq j < j' \leq J$.
\end{theorem}

\begin{proof} Let us fix $D,d,J,\eps$, and allow all implied constants to depend on these quantities.
From Corollary \ref{corr} we have
$$\E_{x \in I \cap \Z} \prod_{j \in [J]} \nu(x + P_j(\vec m)) = 1 + o(1)
 + O\left( \EXP\left( O( \sum_{p \in \primes_{\vec m}} \frac{1}{p} ) \right) \right)$$
 for all $\vec m \in \Omega' \cap \Z^D$, where $\primes_{\vec m}$ are the collection of primes 
 $w \leq p \leq R^{\log R}$ such that $p | P_j(\vec m) - P_{j'}(\vec m)$ for some $1 \leq j < j' \leq d$.  Thus it suffices to show that
$$ \E_{\vec m \in \Omega \cap \Z^D} \EXP\left( O( \sum_{p \in \primes_{\vec m}} \frac{1}{p} ) \right)= o\left(\exp\left(O\left(\sum_{p \in \primes_b} \frac{1}{p} \right) \right)\right).$$
Applying \eqref{exp-bound}, we reduce to showing that
$$ \E_{\vec m \in \Omega \cap \Z^D} \EXP\left( O( \sum_{p \in \primes_{\vec m} \backslash \primes_b} \frac{1}{p} ) \right) = o(1).$$
Applying Lemma \ref{explog} it suffices to show that
$$ \sum_{w \leq p \leq R^{\log R}; p \not \in \primes_b} \frac{\log^{O(1)} p}{p}
\E_{\vec m \in \Omega \cap \Z^D} 1_{p \in \primes_{\vec m}} = o(1).$$
From \eqref{plog-chebyshev}, \eqref{psum-chebyshev} it will suffice to establish the bounds
$$ \E_{\vec m \in \Omega \cap \Z^D} 1_{p \in \primes_{\vec m}} =
O( \frac{1}{p} ) + O( \frac{1}{R^\eps} ) $$
for any $w \leq p \leq R^{\log R}$ with $p \not \in\primes_b$ (note that $\log(R^{\log R})^{O(1)} = o(R^{\eps})$).  By the triangle inequality it suffices to show that
$$ \E_{\vec m \in \Omega \cap \Z^D} 1_{p| P_j(\vec m) - P_{j'}(\vec m)} =
O( \frac{1}{p} ) + O( \frac{1}{R^\eps} ) $$
for all $1 \leq j < j' \leq J$.

Fix $j,j'$. Observe that the property $p| P_j(\vec m) - P_{j'}(\vec m)$ is periodic in each component of $\vec m$ of period $p$, and can thus meaningfully be defined for $\vec m \in F_p^{D}$.  Applying Corollary \ref{avg} (for $p \ll R^\eps$) or Lemma \ref{cover} (for $p \gg R^\eps$) it will thus suffice to show the bound
$$\E_{m_1 \in A_1, \ldots, m_D \in A_D} 1_{p| P_j(m_1,\ldots,m_D) - P_{j'}(m_1,\ldots,m_D)} = O(\frac{1}{M})$$
for all subsets $A_1,\ldots,A_D$ in $F_p$ of size at least $M \geq 1$ for some $M$.  But since $p \not \in \primes_b$, the polynomial $P_j - P_{j'}$ does not vanish modulo $p$, and this claim follows from Lemma \ref{babynull}.
\end{proof}

We can improve the error term if the coefficients of the polynomials are not too large:

\begin{corollary}[Polynomial forms condition, again]\label{pfc2} Let $M,D,d,J \geq 0$ and $\eps > 0$.  Let $P_1,\ldots,P_J \in \Z[\m_1,\ldots,\m_D]$ be distinct polynomials of degree $d$ with all coefficients of size at most $W^M$.  Let $I \subset \R$ be any interval of length at least $R^{4J+1}$, and let $\Omega \subset \R^D$ be any convex body of inradius at least $R^\eps$.  Then
$$ \E_{x \in I \cap \Z; \vec m \in \Omega \cap \Z^D} \prod_{j \in [J]} \nu(x + P_j(\vec m)) = 1 + o_{D,d,\eps,J,M}(1).$$
\end{corollary}

\begin{proof}
Let $\primes_b$ denote the set of all $w \leq p \leq R^{\log R}$ such that
$p | P_j-P_{j'}$ for some $1 \leq j < j' \leq J$.  Since $P_j-P_{j'}$ is nonzero, this $p$ must then divide a non-zero difference of two of the coefficients of the $P_j$, which is $O(W^M)$.  Thus the total product of all such $p$ is at most $O(W^{O(1)})$, and hence by Lemma \ref{divbound} we have $\sum_{p \in \primes_b} \frac{1}{p} = o(1)$.   The claim now follows from Theorem \ref{pfc}.
\end{proof}

From Corollary \ref{pfc2}, the desired estimate \eqref{polyform-2} quickly follows.

\section{The polynomial correlation condition}\label{pcc-sec}

Now we use the estimates from Section \ref{init-sec} to prove the polynomial correlation condition \eqref{polycor-eq-2}.  It will suffice to prove the following estimate.

\begin{theorem}[Polynomial correlation condition]\label{pcc} Let $B,D,D',D'',d,J,K,L \geq 0$ and $\eps > 0$.  For any $j \in [J]$, $k \in [K]$, $l \in [L]$, let
\begin{align*}
\vec P_j & \in \Z[ \h_1, \ldots, \h_{D''} ]^D \\
Q_{j,k,} & \in \Z[ \h_1,\ldots,\h_{D''} ]^{D'} \\
S_{l} &\in \Z[\h_1,\ldots,\h_{D''}]^{D'}
\end{align*}
be polynomials obeying the following conditions.  
\begin{itemize}
\item For any $1 \leq j < j' \leq d$ and $1 \leq k \leq K$, the vector-valued polynomials 
$(\vec P_j,\vec Q_{j,k})$ and $(\vec P_{j'},\vec Q_{j',k})$ are not parallel.
\item The coefficients of the $P_{j,d}$ and $S_{l,d'}$ are bounded in magnitude by $W^B$.
\item The vector-valued polynomials $\vec S_l$ are distinct as $l$ varies in $[L]$.
\end{itemize}  
Let $I \subset \R$ be any interval of length at least $R^{4L+1}$, let $\Omega \subset \R^D$ be a bounded convex body of inradius at least $R^{8J+2}$, and let $\Omega' \subset \R^{D'}$ and $\Omega'' \subset \R^{D''}$ have inradii at least $R^\eps$.  Suppose also that $\Omega''$ is contained in the ball $B(0,R^B)$.  Then
\begin{equation}\label{hugemess2}
\begin{aligned}
\E_{x \in I, \vec n \in \Omega' \cap \Z^{D'}, \vec h \in \Omega'' \cap \Z^{D''}} 
&\left[ \prod_{k \in [K]} 
\E_{\vec m \in \Omega \cap \Z^D} \prod_{j \in [J]} \nu(x + \vec P_{j}(\vec h) \cdot \vec m + \vec Q_{j,k}(\vec h) \cdot \vec n) \right] \\
& \prod_{l \in [L]} \nu( x + \vec S_{l}(\vec h) \cdot \vec n )\\
&\quad = 1 + o_{B,D,D',D'',d,J,K,L,\eps}(1).
\end{aligned}
\end{equation}
\end{theorem}

\begin{proof}
We repeat the same strategy of proof as in the preceding section.  We fix $B,D,D',D'',d,J,K,L,\eps$ and
allow implicit constants to depend on these parameters.  Thus for instance the right-hand side of \eqref{hugemess2} is now simply $1+o(1)$.

We begin by fixing $k, x, \vec h, \vec n$ and 
considering a single average
$$ \E_{\vec m \in \Omega \cap \Z^D} \prod_{j \in [J]} \nu(x + \vec P_j(\vec h) \cdot \vec m + \vec Q_{j,k}(\vec h) \cdot \vec n) .$$
By Proposition \ref{corrprop}, this average is
\begin{equation}\label{precrude}
 1_{\primes_t[k,x,\vec h,\vec n] = \emptyset} + o(1) +  O\left( \EXP\left( O\left( \sum_{p \in \primes_b[k,x,\vec h,\vec n]} \frac{1}{p} \right) \right) \right)
\end{equation}
where $\primes_t[k,x,\vec h,\vec n]$ are the collection of primes $w \leq p \leq R^{\log R}$ which are terrible with respect to the linear polynomials
\begin{equation}\label{wlin}
 W \times [x + \vec P_j(\vec h) \cdot \vec \m + \vec Q_{j,k}(\vec h) \cdot \vec n] + b \in \Z[ \m_1,\ldots,\m_D] \hbox{ for } j \in [J],
\end{equation}
and $\primes_b[k,x,\vec h,\vec n]$ are the collection of primes which are bad.  We can thus express
$$ \prod_{k \in [K]} \E_{\vec m \in \Omega \cap \Z^D} \prod_{j \in [J]} \nu(x + \vec P_j(\vec h) \cdot \vec m + \vec Q_{j,k}(\vec h) \cdot \vec n)  $$
using \eqref{exp-bound} as
$$ 1_{\primes_t[k,x,\vec h,\vec n] = \emptyset \hbox{ for all } k} + o(1) + 
O\left( \EXP\left( O\left( \sum_{p \in \bigcup_{k=1}^K \primes_b[k,x,\vec h,\vec n]} \frac{1}{p} \right) \right) \right)$$
which we estimate crudely by
$$ 1 + o(1) + O\left( \sum_{k \in [K]} \sum_{p \in\primes_t[k,x,\vec h,\vec n]} 1 \right)
+ O\left( \EXP\left( O\left( \sum_{p \in \bigcup_{k \in [K]} \primes_b[k,x,\vec h,\vec n] \backslash \primes_t[k,x,\vec h,\vec n]} \frac{1}{p} \right) \right) \right).$$

Observe that if $p \geq w$ is terrible for \eqref{wlin}, then
$p | \vec P_j(\vec h)$ and $p | x + \vec Q_{j,k}(\vec h) \cdot \vec n$
for some $j \in [J]$, while if $p \geq w$ is bad but not terrible,
then 
$$p | (\vec P_j(\vec h), \vec Q_{j,k}(\vec h) \cdot \vec n) \wedge 
(\vec P_{j'}(\vec h), \vec Q_{j',k}(\vec h) \cdot \vec n)$$
for some $1 \leq j < j' \leq J$, where $\wedge$ denotes the wedge product on $D+1$-dimensional space.  
Thus we may estimate the preceding sum (using \eqref{exp-bound}) by
\begin{align*}
1 &+ o(1) + \sum_{j \in [J]} \sum_{k \in [K]} O( \sum_{w \leq p \leq R^{\log R}: p | \vec P_j(\vec h), x + \vec Q_{j,k}(\vec h) \cdot \vec n} 1 )^{1/2}\\
&+ \sum_{1 \leq j < j' \leq J} 
O\left( \EXP\left( O\left( \sum_{w \leq p \leq R^{\log R}: p | (\vec P_j(\vec h), \vec Q_{j,k}(\vec h) \cdot \vec n) \wedge 
(\vec P_{j'}(\vec h), \vec Q_{j',k}(\vec h) \cdot \vec n)} \frac{1}{p} \right) \right) \right)^{1/2}.
\end{align*}
At this point we pause to remove some ``globally bad'' primes.  Let $\primes_b$ denote the primes $w \leq p \leq R^{\log R}$ which divide $(\vec P_j, \vec Q_{j,k}) \wedge (\vec P_{j'}, \vec Q_{j'})$
for some $1 \leq j < j' \leq J$ and $k \in [K]$ (note this is now the wedge product in $D+D'$ dimensions).  Because the wedge products $(\vec P_j, \vec Q_{j,k}) \wedge (\vec P_{j'}, \vec Q_{j'})$ are
non-zero and have coefficients $O(W^{O(1)})$, the product of all these primes is $O( W^{O(1)} )$, and hence by Lemma \ref{divbound} we have $\sum_{p \in \primes_b} \frac{1}{p} = o(1)$.  Thus we may safely delete these primes from the expression inside the $\EXP()$.  If we then apply Lemma \ref{explog}, we can bound the above sum as
$$ \begin{aligned}
& 1 + o(1) + \sum_{j \in [J]} O\left( \sum_{w \leq p \leq R^{\log R}: p | \vec P_j(\vec h), x + \vec Q_{j,k}(\vec h) \cdot \vec n} 1 \right)\\
&+ \sum_{1 \leq j < j' \leq J} O\left(\sum_{w \leq p \leq R^{\log R}; p \not \in \primes_b} \frac{\log^{O(1)} p}{p}
1_{p | (\vec P_j(\vec h), \vec Q_{j,k}(\vec h) \cdot \vec n) \wedge 
(\vec P_{j'}(\vec h), \vec Q_{j',k}(\vec h) \cdot \vec n)}\right)^{1/2}.
\end{aligned}$$
Inserting this bound into \eqref{hugemess2} and using Cauchy-Schwarz, we reduce to showing the bounds
\begin{align}
&\E_{x \in I; \vec n \in \Omega'; \vec h \in \Omega''} \prod_{l \in [L]} \nu( x + \vec S_l(\vec h) \cdot \vec n ) 
= 1 + o(1) \label{bound-1}\\
&\E_{x \in I; \vec n \in \Omega'; \vec h \in \Omega''} 
\sum_{w \leq p \leq R^{\log R}: p | \vec P_j(\vec h),x + \vec Q_{j,k}(\vec h) \cdot \vec n}  \prod_{l \in [L]} \nu( x + \vec S_l(\vec h) \cdot \vec n ) = o(1) \label{bound-3}\\
&\E_{x \in I; \vec n \in \Omega'; \vec h \in \Omega''} 
\sum_{w \leq p \leq R^{\log R}; p \not \in \primes_b} 
\frac{\log^{O(1)} p}{p} 1_{p | (\vec P_j(\vec h), \vec Q_{j,k}(\vec h) \cdot \vec n) \wedge 
(\vec P_{j'}(\vec h), \vec Q_{j',k}(\vec h) \cdot \vec n)} \nonumber\\
& \qquad  \qquad \prod_{l \in [L]} \nu( x + \vec S_l(\vec h) \cdot \vec n ) 
= o(1) \label{bound-4}
\end{align}
for all $1 \leq j < j' \leq J$ and $k \in [K]$.

The bound \eqref{bound-1} already follows from Corollary \ref{pfc2} and the hypotheses on $S_{l,d'}$.  Now we turn to \eqref{bound-3}.  We rewrite the left-hand side as
\begin{equation}\label{b3-lhs}
\sum_{w \leq p \leq R^{\log R}} \E_{\vec n \in \Omega'; \vec h \in \Omega''}
\left( 1_{p | \vec P_j(\vec h)}
\E_{x \in I} 1_{x = -\vec Q_{j,k}(\vec h) \cdot \vec n\ \mod\ p} 
\prod_{l \in [L]} \nu( x + \vec S_l(\vec h) \cdot \vec n ) \right).
\end{equation}
Let us first consider the contributions of the primes $p$ which are larger than $R^{4L+1}$.  In this case we bound
$\nu$ extremely crudely by $O( R^2 \log R )$ (taking absolute values in \eqref{nudef}), to bound the inner expectation of \eqref{b3-lhs} by
$O( \frac{1}{R^{4L+1}} (R^2 \log R)^L ) = o( R^{-1/2} )$
(say), and  to show that
$$ \E_{\vec h \in \Omega''} \sum_{R^{4L} \leq p \leq R^{\log R}} 1_{p | \vec P_j(\vec h)} \ll R^{1/2}.$$
But from the bounds on $\Omega''$ and $\vec h$ we see that $\vec P_j(\vec h) = O(R^{O(1)})$, and so at most $O(1)$ primes $p$ can contribute to the sum for each $\vec h$.  The claim follows.

Now we consider the contributions of the primes $p$ beteen $w$ and $R^{4L+1}$.  We can then apply Corollary \ref{corr3}
and estimate the inner expectation of \eqref{b3-lhs} by
$$ \frac{1}{p} O\left( \exp\left( O\left( \sum_{1 \leq l < l' \leq L} \sum_{w \leq p' \leq R^{\log R}: p' | \vec S_l(\vec h) - \vec S_{l'}(\vec h) \frac{1}{p'} } \right) \right) \right)$$
which by Lemma \ref{explog} can be bounded by
$$ O(\frac{1}{p}) + \sum_{1 \leq l < l' \leq L} \sum_{w \leq p' \leq R^{\log R}} \frac{\log^{O(1)} p'}{pp'} O(1_{p' | \vec S_l(\vec h) - \vec S_{l'}(\vec h)}).$$
The contribution to \eqref{b3-lhs} can thus be bounded by the sum of 
$$ \sum_{w \leq p \leq R^{\log R}} \frac{1}{p} O( \E_{\vec h \in \Omega''} 1_{p | \vec P_j(\vec h)} ) $$
and
$$ \sum_{1 \leq l < l' \leq L} \sum_{w \leq p,p' \leq R^{\log R}} \frac{\log^{O(1)} p'}{pp'}
\E_{\vec h \in \Omega''} 1_{p | \vec P_j(\vec h)} 1_{p' | \vec S_l(\vec h) - \vec S_{l'}(\vec h)}.$$
Now by hypothesis, the vector-valued polynomials $\vec P_j$ and $\vec S_l - \vec S_{l'}$ are non-zero, thus by Lemma \ref{babynull} we have
$$ \E_{\vec h \in A_1 \times \ldots \times A_{D''}} 1_{p' | \vec S_l(\vec h) - \vec S_{l'}(\vec h)} \ll \frac{1}{M}$$
whenever $A_1,\ldots,A_{D''} \subset F_{p'}$ have cardinality at least $M$.  Applying Corollary \ref{avg} and Lemma \ref{cover}
we conclude that
$$ \E_{\vec h \in \Omega''} 1_{p' | \vec S_l(\vec h) - \vec S_{l'}(\vec h)} \ll \frac{1}{p'}  +  \frac{1}{R^\eps}.$$
A similar argument gives
$$ \E_{\vec h \in \Omega''} 1_{p | \vec P_j(\vec h)} \ll \frac{1}{p} + \frac{1}{R^\eps}$$
and hence by Cauchy-Schwarz
$$
\E_{\vec h \in \Omega''}
 1_{p | \vec P_j(\vec h)} 1_{p' | \vec S_l(\vec h) - \vec S_{l'}(\vec h)}
\ll \frac{1}{(pp')^{1/2}}  +  \frac{1}{R^{\eps/2}} .$$
Applying all of these bounds, we can thus bound the total contribution of this case to \eqref{b3-lhs} by
$$ \sum_{w \leq p \leq R^{\log R}} \frac{1}{p} [ O(\frac{1}{p} ) + O( \frac{1}{R^\eps} ) ] +
\sum_{w \leq p,p' \leq R^{\log R}} \frac{\log^{O(1)} p'}{pp'}
\left[ O( \frac{1}{(pp')^{1/2}} ) + O( \frac{1}{R^{\eps/2}} ) \right] $$
which is $o(1)$ by \eqref{plog-chebyshev}, \eqref{psum-chebyshev}.

Finally, we consider \eqref{bound-4}.  We first apply Theorem \ref{pfc} to bound
$$\E_{x \in I} \prod_{l \in [L]} \nu( x + \vec S_l(\vec h) \cdot \vec n )  = 
O\left( \exp\left( O\left( \sum_{1 \leq l < l' \leq L} \sum_{w \leq p' \leq R^{\log R}:
p' | [\vec S_l(\vec h) - \vec S_{l'}(\vec h)] \cdot \vec n} \frac{1}{p'} \right) \right) \right).
$$
Once again we must extract out the ``globally bad'' primes.  Let $\prod'_b$ denote all the primes
$w \leq p' \leq R^{\log R}$ which divide $\vec S_l - \vec S_{l'}$ for some $1 \leq l < l' \leq L$.
Since these polynomials are non-zero and have coefficients $O(W^{O(1)})$, the product of all the primes in $\prod'_b$ is
$O(W^{O(1)})$, and hence by Lemma \ref{divbound} as before these primes contribute only $o(1)$ and can be discarded.
If we then apply Lemma \ref{explog}, we can bound the preceding expression by
$$
O(1) + \sum_{1 \leq l < l' \leq L} \sum_{w \leq p' \leq R^{\log R}: p' \not \in \prod'_b} \frac{\log^{O(1)} p'}{p'}
O( 1_{p' | [(\vec S_l(\vec h) - \vec S_{l'}(\vec h))] \cdot \vec n} ).
$$
Thus to prove \eqref{bound-4} it suffices to show the estimates
$$
 \sum_{w \leq p \leq R^{\log R}; p \not \in \primes_b} \frac{\log^{O(1)} p}{p} 
\E_{\vec n, \vec h} 1_{p | (\vec P_j(\vec h), \vec Q_{j,k}(\vec h) \cdot \vec n) \wedge 
(\vec P_{j'}(\vec h), \vec Q_{j',k}(\vec h) \cdot \vec n)} = o(1)
$$
and
$$\begin{aligned}
&\sum_{w \leq p,p' \leq R^{\log R}; p \not \in \primes_b,   p' \not \in \prod'_b}  \frac{\log^{O(1)} p \log^{O(1)} p'}{pp'}  \\
&  \qquad \qquad \qquad  \E_{\vec n, \vec h} 1_{p | (\vec P_j(\vec h), \vec Q_{j,k}(\vec h) \cdot \vec n) \wedge 
(\vec P_{j'}(\vec h), \vec Q_{j',k}(\vec h) \cdot \vec n)} 
1_{p' | (\vec S_l(\vec h) - \vec S_{l'}(\vec h)) \cdot \vec n}
= o(1)
\end{aligned}$$
for all $1 \leq j < j' \leq J$ and $1 \leq l < l' \leq L$, where $ \vec n, \vec h$ are averaged over  
$\Omega' \cap \Z^{D'},  \Omega'' \cap \Z^{D''}$ respectively.   Applying Cauchy-Schwarz to the expectation in latter inequality and then factorising the double sum, we see that it will suffice to show that
\begin{equation}\label{b4a}
 \sum_{w \leq p \leq R^{\log R}; p \not \in\primes_b} \frac{\log^{O(1)} p}{p} 
\left( \E_{\vec n, \vec h} 1_{p | (\vec P_j(\vec h), \vec Q_{j,k}(\vec h) \cdot \vec n) \wedge 
(\vec P_{j'}(\vec h), \vec Q_{j',k}(\vec h) \cdot \vec n)} \right)^{1/2} = o(1)
\end{equation}
and
\begin{equation}\label{b4b}
 \sum_{w \leq p' \leq R^{\log R}; p \not \in \prod_{b'}} \frac{\log^{O(1)} p'}{p'} 
\left( \E_{\vec n, \vec h} 1_{p | (\vec S_l(\vec h) - \vec S_{l'}(\vec h)) \cdot \vec n} \right)^{1/2} = o(1).
\end{equation}
Since $p \not \in \primes_b$, we observe from
Lemma \ref{babynull}, Corollary \ref{avg}, Lemma \ref{cover} that
$$ \E_{\vec h} 1_{p | (\vec P_j(\vec h), \vec Q_{j,k}(\vec h) \cdot \vec n) \wedge 
(\vec P_{j'}(\vec h), \vec Q_{j',k}(\vec h) \cdot \vec n)}  \ll \frac{1}{p}  + \frac{1}{R^\eps} $$
and the claim \eqref{b4a} now follows from \eqref{plog-chebyshev}, \eqref{psum-chebyshev}.  The estimate \eqref{b4b} is proven
similarly.  This (finally!) completes the proof of Theorem \ref{pcc}.
\end{proof}

\appendix

\section{Local Gowers uniformity norms}\label{gowers-sec}

In this appendix we shall collect a number of elementary inequalities based on the Cauchy-Schwarz
inequality, including several related to Gowers-type uniformity norms.

The formulation of the Cauchy-Schwarz inequality which we shall rely on is
\begin{equation}\label{cauchy-schwarz}
|\E_{a \in A, b \in B} f(a) g(a,b)|^2 \leq
(\E_{a \in A} F(a)) (\E_{a \in A} F(a) |\E_{b \in B} g(a,b)|^2)
\end{equation}
whenever $f: A \to \R$, $F: A \to \R^+$, $g: A \times B \to \R$ are functions on non-empty finite sets $A,B$ with the pointwise bound $|f| \leq F$.

One well known consequence of Cauchy-Schwarz is the \emph{van der Corput lemma}, which allows one to estimate a coarse-scale average of a function $f$ by coarse-scale averages of ``derivatives'' of $f$ over short scales.  The precise formulation we need here is as follows.

\begin{lemma}[van der Corput]\label{VDC} Let $N$, $M$ and $H$ be as in Section \ref{notation-sec}.
Let $(x_m)_{m \in \Z}$ be a sequence of real numbers obeying the bound
\begin{equation}\label{xmh}
 x_m \ll_\eps N^\eps 
 \end{equation}
for any $\eps > 0$ and $m \in \Z$.  Then we have
\begin{equation}\label{xmh-1}
\E_{m \in [M]} x_m = \E_{h \in [H]} \E_{m \in [M]} x_{m+h} + o(1)
\end{equation}
and
\begin{equation}\label{xmh-2}
|\E_{m \in [M]} x_m|^2 \ll \E_{h, h' \in [H]} \E_{m \in [M]} x_{m+h} x_{m+h'} + o(1).
\end{equation}
\end{lemma}

\begin{proof}
From \eqref{xmh} we see that
$$ \E_{m \in [M]} x_m = \E_{m \in [M]} x_{m+h} + o(1)$$
for all $h \in [H]$; averaging over all $h$ and rearranging we obtain \eqref{xmh-1}.
Applying \eqref{cauchy-schwarz} we conclude
$$ |\E_{m \in [M]} x_m| \ll \E_{m \in [M]} |\E_{h \in [H]} x_{m+h}|^2 + o(1)$$
and \eqref{xmh-2} follows.
\end{proof}

We will use Lemma \ref{VDC} in only one place, namely Proposition \ref{linearize}, which is the key inductive step needed to estimate a polynomial average by a collection of linear averages.

Next, we recall some Cauchy-Schwarz-Gowers inequalities, which can be found for instance in \cite[Appendix B]{linprimes}.  Let $X$ be a finite non-empty set.
If $A$ is a finite set and $f: X^A \to \R$, define the \emph{Gowers box norm} $\|f\|_{\Box^A}$ as
\begin{equation}\label{boxdef}
 \| f \|_{\Box^A} := \left(\E_{m^{(0)}, m^{(1)} \in X^A}
\prod_{\omega \in \{0,1\}^A} 
f( (m^{(\omega_\alpha)}_\alpha)_{\alpha \in A} )\right)^{1/2^{|A|}}
\end{equation}
where $\omega = (\omega_\alpha)_{\alpha \in A}$ and $m^{(i)} = (m^{(i)}_\alpha)_{\alpha \in A}$ for $i=0,1$.
This is indeed a norm\footnote{If $|A|=0$ then $\|f\|_{\Box^A} = f(\emptyset)$, while if $|A|=1$ then $\|f\|_{\Box^A} = |\E f|$.} for $|A| \geq 2$.  It obeys the \emph{Cauchy-Schwarz-Gowers inequality}
\begin{equation}\label{gcz}
 |\E_{\vec m^{(0)}, \vec m^{(1)} \in X^A}
\prod_{\omega \in \{0,1\}^A} 
f_\omega( (m^{(\omega_\alpha)}_\alpha)_{\alpha \in A} )| \leq \prod_{\omega \in \{0,1\}^A}
|\|f_\omega\|_{\Box^A}|.
\end{equation}
We shall also need a weighted variant of this inequality.

\begin{proposition}[Weighted generalized von Neumann inequality]\label{wgn} Let $A$ be a non-empty finite set, and let $f: X^A \to \R$ be a function.  For every $\alpha \in A$, let $f_\alpha: X^{A \backslash \{\alpha\}} \to \R$ and $\nu_\alpha: X^{A \backslash \{\alpha\}} \to \R^+$ be functions with the pointwise bound $|f_\alpha| \leq \nu_\alpha$.  Then we have
$$ 
|\E_{\vec m \in X^A} f(\vec m) \prod_{\alpha \in A} f_\alpha(\vec m|_{A \backslash \{\alpha\}})|
\leq |\| f \|_{\Box^A(\nu)}| \prod_{\alpha \in A} \|\nu_\alpha\|_{\Box^{A \backslash \{\alpha\}}}^{1/2}$$
where $\vec m|_{A \backslash \{\alpha\}}$ is the restriction of $\vec m \in X^A$ to $X^{A \backslash \{\alpha\}}$, and $\|f\|_{\Box^A(\nu)}$ is the weighted Gowers box norm of $f$, defined by the formula
\begin{align*}
\| f \|_{\Box^A(\nu)}^{2^{|A|}} &:= \E_{\vec m^{(0)}, \vec m^{(1)} \in X^A}
\left[\prod_{\omega \in \{0,1\}^A} f( (m^{(\omega_\alpha)}_\alpha)_{\alpha \in A} )\right]\\
&\quad \times
\prod_{\alpha \in A} \prod_{\omega^{(\alpha)} \in \{0,1\}^{A \backslash \{\alpha\}}}
\nu_\alpha\left( ( m^{(\omega^{(\alpha)}_\beta)}_\beta )_{\beta \in A \backslash \{\alpha\}} \right).
\end{align*}
\end{proposition}

\begin{proof}  This is a special case of \cite[Corollary B.4]{linprimes}, in which all functions $f_B$, $\nu_B$ associated to subsets $B$ of $A$ of cardinality $|A|-2$ or less are set equal to one.
\end{proof}

The local Gowers norm $U^{a_1,\ldots,a_d}_{\sqrt{M}}$ defined in \eqref{local-gowers} is related to the above Gowers box norms by the obvious identity
\begin{equation}\label{fum}
 \|f\|_{U^{a_1,\ldots,a_d}_{\sqrt{M}}} = (\E_{x \in X} \|F_x\|_{\Box^d}^{2^d})^{1/2^d}
 \end{equation}
for any $f: X \to \R$, where for each $x \in X$ the function $F_x: [\sqrt{M}]^d \to \R$ is defined by
$$ F_x( m_1,\ldots,m_d ) := f( x + a_1 m_1 + \ldots + a_d m_d ).$$
In particular, since $\Box^d$ is norm for $d \geq 2$, we easily verify from Minkowski's inequality that $U^{a_1,\ldots,a_d}_{\sqrt{M}}$ is a norm also when $d \geq 2$.  This in turn implies that 
the averaged local Gowers norms $U^{\vec Q([H]^t,W)}_{\sqrt{M}}$ are also indeed norms.

Now we introduce the concept of concatenation of two or more averaged local Gowers norms.  If
$\vec Q \in \Z[\h_1,\ldots,\h_t,\W]^d$ and $\vec Q' \in \Z[\h'_1,\ldots,\h'_{t'},\W]^{d'}$ are a $d$-tuple and $d'$-tuple of polynomials respectively, we define the \emph{concatenation} $\vec Q \oplus \vec Q' \in \Z[\h_1,\ldots,\h_t,\h'_1,\ldots,\h'_t,\W]^{d+d'}$ to be the $d+d'$-tuple of polynomials whose first $d$ components are those of $\vec Q$ (using the obvious embedding of $\Z[\h_1,\ldots,\h_t,\W]$ into $\Z[\h_1,\ldots,\h_t,\h'_1,\ldots,\h'_{t'},\W]$) and the last $d'$ components are those of $\vec Q'$ (using the obvious embedding of
$\Z[\h'_1,\ldots,\h'_{t'},\W]$ into $\Z[\h_1,\ldots,\h_t,\h'_1,\ldots,\h'_t,\W]$).  One can similarly define concatenation of more than two tuples of polynomials in the obvious manner.

The key lemma concerning concatenation is as follows.

\begin{lemma}[Domination lemma]\label{uk-concatenate}  Let $k \geq 1$.  For each $1 \leq i \leq k$, let
$t_i \geq 0$ be an integer and $\vec Q_i \in \Z[\h_1,\ldots,\h_{t_i},\W]^{d_i}$ be a polynomial.  Let $t := t_1 + \ldots + t_k$, $d := d_1 + \ldots + d_k$, and let $\vec Q \in \Z[\h_1,\ldots,\h_{t},\W]^d$ be the concatenation of all the $\vec Q_i$.  Then we have
$$ \| g \|_{U^{\vec Q_i([H]^{t_i},W)}_{\sqrt{M}}} \leq \|g\|_{U^{\vec Q([H]^t,W)}_{\sqrt{M}}}$$
for all $1 \leq i \leq k$ and $g: X \to \R$.
\end{lemma}

\begin{proof}  By induction we may take $k=2$.  By symmetry it thus suffices to show that
$$ \| g \|_{U^{\vec Q([H]^{t},W)}_{\sqrt{M}}}
 \leq \|g\|_{U^{\vec Q \oplus \vec Q'([H]^{t+t'},W)}_{\sqrt{M}}}$$
for any $g: X \to \R$ and any $\vec Q \in \Z[\h_1,\ldots,\h_t,W]^d$ and $\vec Q' \in \Z[\h'_1,\ldots,\h'_{t'},W]^{d'}$.  We may take $d' \geq 1$ as the case $d'=0$ is trivial.
From \eqref{avg-def} and H\"older's inequality it suffices to prove the estimate
$$ \| g \|_{U^{a_1,\ldots,a_d}_{\sqrt{M}}} \leq \|g\|_{U^{a_1,\ldots,a_d,a'_1,\ldots,a'_{d'}}_{\sqrt{M}}}$$
for all $a_1,\ldots,a_d,a'_1,\ldots,a'_d \in \Z$. Applying \eqref{fum}, we see that it suffices to prove the monotonicity formula\footnote{This is of course closely connected with the monotonicity of the Gowers $U^d$ norms, noted for instance in \cite{gt-primes}.}
$$ \| f \|_{\Box^d([\sqrt{M}])} \leq \| f \|_{\Box^{d+d'}([\sqrt{M}])}$$
for any $f: [\sqrt{M}]^d \to \R$, where we extend $f$ to $[\sqrt{M}]^{d+d'}$ by adding $d'$ dummy variables, thus
$$ f( m_1,\ldots,m_d, m_{d+1},\ldots,m_{d+d'} ) := f(m_1,\ldots,m_d).$$
But this easily follows by raising both sides to the power $2^{d}$ and using the Cauchy-Schwarz-Gowers
inequality \eqref{gcz} for the $\Box^{d+d'}$ norm (setting $2^d$ factors equal to $f$, and the other $2^{d+d'}-2^d$ factors equal to $1$).
\end{proof}

\section{Uniform polynomial Szemer\'edi theorem}\label{uniform-sec}

In this appendix we use the Furstenberg correspondence principle and the Bergelson-Leibman theorem \cite{bl} to prove the quantitative polynomial Szemer\'edi theorem, Theorem \ref{pSZ-quant}. 
The arguments here are reminiscent of those in \cite{bmht}; see also \cite{tao-transference} for another argument in a similar spirit.

Firstly, observe to prove Theorem \ref{pSZ-quant} it certainly suffices to do so in the case when $g$ is an indicator function $1_E$, since in the general case one can obtain a lower bound $g \geq \frac{\delta}{2} 1_E$ where $E := \{ x \in X: g(x) \geq \delta/2\}$, which must have measure at least $\delta/2 - o(1)$.

Fix $P_1,\ldots,P_k$ and $\delta$, and suppose for contradiction that Theorem \ref{pSZ-quant} failed.  Then (by the axiom of choice\footnote{It is not difficult to rephrase this argument so that the axiom of choice is not used; we leave the details to the interested reader.  The weak sequential compactness of probability measures which we need later in this section can also be established by an Arzela-Ascoli type diagonalization argument which avoids the axiom of choice.  On the other hand, the only known proof of the multi-dimensional Bergelson-Leibman theorem does use the axiom of choice, and so the main result of this paper also currently requires this axiom.  However, it is expected that the Bergelson-Leibman theorem (and hence our result also) will eventually be provable by means which avoid using this axiom.  For instance, for the one-dimensional Bergelson-Leibman theorem one can use the Gowers-Host-Kra seminorm characteristic factors as in \cite{nikos-kra}, which do not require choice.}) we can find a sequence of $N$ going to infinity, and a sequence of indicator functions $1_{E_N}: \Z/N\Z \to \R$ of density 
\begin{equation}\label{mudens}
|E_N|/N \geq \delta - o(1)
\end{equation}
such that
$$ \lim_{N \to \infty} \E_{m \in [M]} \int_{\Z/N\Z} \prod_{i \in [k]} T_N^{P_i(Wm)/W} 1_{E_N} = 0$$
where $T_N$ is the shift on $\Z/N\Z$ and $N$ is always understood to lie along the sequence (recall that $M$ and $W$ both depend on $N$).

The next step is to use an averaging argument (dating back to Varnavides \cite{var}) to deal with the fact that $M$ is growing rather rapidly in $N$.  Let $B \geq 1$ be an integer, then for $N$ sufficiently large we have
$$ \E_{m \in [M]} \int_{\Z/N\Z} \prod_{i \in [k]} T_N^{P_i(Wm)/W} 1_{E_N} \leq \frac{1}{B^3}$$
and hence
$$ \E_{m \in [M/B]} \int_{\Z/N\Z} \prod_{i \in [k]} T_N^{P_i(Wbm)/W} 1_{E_N} \ll \frac{1}{B^2} $$
for all $b \in [B]$.  In particular
$$ \E_{m \in [M/B]} \sum_{b=1}^B \int_{\Z/N\Z} \prod_{i \in [k]} T_N^{P_i(Wbm)/W} 1_{E_N} \ll \frac{1}{B} $$
and hence by the pigeonhole principle (and the axiom of choice) we can find $m_N \in [M/B]$ for all sufficiently large $N$ such that
$$ \sum_{b=1}^B \int_{\Z/N\Z} \prod_{i \in [k]} T_N^{P_i(Wbm_N)/W} 1_{E_N} \ll \frac{1}{B} $$
and hence
$$ \lim_{N \to \infty} \int_{\Z/N\Z} \prod_{i \in [k]} T_N^{P_i(Wbm_N)/W} 1_{E_N} = 0$$
for each $b \geq 1$.

We now eliminate the $W$ and $m_N$ dependence by ``lifting'' the one-dimensional shift to many dimensions.  Let $d$ be the maximum degree of the $P_1,\ldots,P_k$, then we may write
$$ P_i(Wbm_N)/W = \sum_{j \in [d]} W^{j-1} m_N^j b^j c_{i,j}$$
for some integer constants $c_{i,j}$.  Thus, if we set $T_{N,j} := T_N^{W^{j-1} m_N^j}$, we have
\begin{equation}\label{emn}
\lim_{N \to \infty} \E_{m \in [M]} \int_{\Z/N\Z} \prod_{i \in [k]} (\prod_{j \in [d]} T^{c_{i,j} b^j}_{N,j}) 1_{E_N} = 0 \hbox{ for all } b \geq 1.
\end{equation}
Now we use the Furstenberg correspondence principle to take a limit.  Let $\Omega$ be the product space
$\Omega := \{0,1\}^{\Z^d}$, endowed with the usual product $\sigma$-algebra and with the standard commuting shifts $T_1,\ldots,T_d$ defined by
$$ T_j ( ( \omega_n )_{n \in \Z^d} ) := (\omega_{n-e_j} )_{n \in \Z^d} \hbox{ for } j \in [d]$$
where $e_1,\ldots,e_d$ is the standard basis for $\Z_d$.  We can define a probability measure $\mu_N$ on this space by $\mu_N := \E_{x \in \Z/N\Z} \mu_{N,x}$, where $\mu_{N,x}$ is the Dirac measure on the point
$$ ( 1_{T_1^{n_1} \ldots T_d^{n_d} x \in E_N} )_{n \in \Z^d}.$$
One easily verifies that $\mu_N$ is invariant under the commuting shifts $T_1,\ldots,T_d$.  Also if we let $A \subset \Omega$ be the cylinder set 
$$ A := \{ (\omega_n)_{n \in \Z^d}: \omega_0 = 1 \}$$
then we see from \eqref{mudens}, \eqref{emn} that $\mu_N(A) \geq \delta-o(1)$ and
$$ \lim_{N \to \infty} \mu_N( \bigcap_{i \in [k]} (\prod_{j \in [d]} T_j^{c_{i,j} b^j}) A) = 0$$
for all $b \geq 1$.  By weak sequential compactness, we may after passing to a subsequence assume that the measures $\mu_N$ converge weakly to another probability measure $\mu$ on $\Omega$, which is thus
translation-invariant and obeys the bounds
$$ \mu(A) \geq \delta > 0$$
and
$$ \mu( \bigcap_{i \in [k]} (\prod_{j \in [d]} T_j^{c_{i,j} b^j}) A) = 0 \hbox{ for all } b \geq 1.$$
But this contradicts the multidimensional Bergelson-Leibman recurrence theorem \cite[Theorem $A_0$]{bl}.  This contradiction concludes the proof of Theorem \ref{pSZ-quant}.

\section{Elementary convex geometry}\label{convex-sec}

In this paper we shall frequently be averaging over sets of the form $\Omega \cap \Z^D$, where $\Omega \subset \R^D$ is a convex body.  It is thus of interest to estimate the size of such sets.  Fortunately we will be able to do this using only very crude estimates (we only need the main term in the asymptotics, and do not need the deeper theory of error estimates).
We shall bound the geometry of $\Omega$ using the inradius $r(\Omega)$; this is more or less dual to the approach in \cite{linprimes}, which uses instead the circumradius.

Observe that the cardinality of $\Omega \cap \Z^D$ equals the Lebesgue measure of the Minkowski sum $(\Omega \cap \Z^D) + [-1/2,1/2]^D$ of $\Omega \cap \Z^D$ with the unit cube $[-1/2,1/2]^D$.  The latter set differs from $\Omega$ only on the $\sqrt{D/2}$-neighborhood ${\mathcal N}_{\sqrt{D}/2}(\partial \Omega)$ of the boundary $\partial \Omega$.  We thus have
the \emph{Gauss bound}
$$
 |\Omega \cap \Z^D| = \mes(\Omega) + O( \mes( {\mathcal N}_{\sqrt{D}/2}(\partial \Omega) ) )
$$
where $\mes()$ denotes Lebesgue measure.  By dilation and translation, we thus have
\begin{equation}\label{gauss}
 |\Omega \cap (m \cdot \Z^D + a)| = m^{-D} [ \mes(\Omega) + O( \mes( {\mathcal N}_{m\sqrt{D}/2}(\partial \Omega) ) )]
\end{equation}
for any $m > 0$ and $a \in \R^D$.

Now we estimate the boundary term in terms of the inradius $r(\Omega)$ of $\Omega$.

\begin{lemma}[Gauss bound]\label{gausslemma}  Suppose that $\Omega \subset \R^D$ is a convex body.   Then for any $0 < r < r(\Omega)$ we have
$$ \mes( {\mathcal N}_r( \partial \Omega ) ) \ll_D \frac{r}{r(\Omega)} \mes( \Omega ).$$
\end{lemma}

\begin{proof} We may rescale $r(\Omega)=1$ (so $0 < r < 1$), and translate so that $\Omega$ contains the open unit ball $B(0,1)$.  Elementary convex geometry then shows that for a sufficiently large constant $C_D > 0$, we have
$$ B(x,r) \subset \Omega \hbox{ whenever } x \in (1 - C_D r) \cdot \Omega$$
and
$$ B(x,r) \cap \Omega = \emptyset \hbox{ whenever } x \not \in (1 + C_D r) \cdot \Omega.$$
This shows that
$$ {\mathcal N}_r( \partial \Omega ) \subset [(1 + C_D r) \cdot \Omega] \backslash [(1 - C_D r) \cdot \Omega]$$
and the claim follows.
\end{proof}

From Lemma \ref{gausslemma} and \eqref{gauss} we conclude that
\begin{equation}\label{gauss-lat}
|\Omega \cap (m \cdot \Z^D + a)| = \left(1 + O_D\left( \frac{m}{r(\Omega)} \right)\right) m^{-D} \mes(\Omega)
\end{equation}
whenever $0 < m \leq r(\Omega)$ and $a \in \Z^D$.  As a corollary we obtain

\begin{corollary}[Equidistribution of residue classes]\label{equi} Let $m \geq 1$ be an integer and $a \in \Z_m^D$, and let $\Omega \subset \R^D$ be a convex body.  If $r(\Omega) \geq C_D m$ for some sufficiently large constant $C_D > 0$, then we have
$$ \E_{x \in \Omega \cap \Z^D} 1_{x \in m \cdot \Z^D + a} = \left(1 + O_D\left( \frac{m}{r(\Omega)} \right)\right) m^{-D}.$$
\end{corollary}

This lets us average $m$-periodic functions on convex bodies as follows.

\begin{corollary}[Averaging lemma]\label{avg} Let $m \geq 1$ be an integer, and let $f: \Z^D \to \R^+$ be a non-negative $m$-periodic function (thus $f$ can also be identified with a function on $\Z_m^D$).  Let $\Omega \subset \R^D$ be a convex body.  If $r(\Omega) \geq C_D m$ for some sufficiently large constant $C_D > 0$, then we have
$$ \E_{x \in \Omega \cap \Z^D} f(x) = \left(1 + O_D\left( \frac{m}{r(\Omega)} \right)\right) \E_{y \in \Z_m^D} f.$$
\end{corollary}

\begin{proof} We expand the left-hand side as
$$ \E_{x \in \Omega \cap \Z^D} f(x) = \sum_{y \in \Z_m^D} f(y) \E_{x \in \Omega \cap \Z^D} 1_{x \in m \cdot \Z^D + y}$$
and apply Corollary \ref{equi}.
\end{proof}

Corollary \ref{avg} is no longer useful when the period $m$ is large compared to the inradius $r(\Omega)$.  In such cases we shall need to rely instead on the following cruder estimate.

\begin{lemma}[Covering inequality]\label{cover}  Let $\Omega \subset \R^D$ be a convex body with $r(\Omega) > C_D$ for some large constant $C_D > 1$, and let $f: \Z^D \to \R^+$ be an arbitrary function.  Then
$$ \E_{x \in \Omega \cap \Z^D} f(x) \ll_D  \sup_{y \in \R^D} \E_{x \in y + [-r(\Omega),r(\Omega)]^D \cap \Z^D} f(x).$$
\end{lemma}

\begin{proof} From \eqref{gauss-lat} we have $|\Omega \cap \Z^D| \sim \mes(\Omega)$, so it suffices to show that
$$ \sum_{x \in \Omega \in \Z^D} f(x) \ll_D 
 \mes(\Omega) \sup_{y \in \R^D} \E_{x \in y + [-r(\Omega),r(\Omega)]^D \in \Z^D} f(x) .$$
This will follow if we can cover $\Omega$ by $O_D( \mes(\Omega) / r(\Omega)^D )$ translates of the cube
$[-r(\Omega),r(\Omega)]^D$.

By rescaling and translating, we reduce to verifying the following fact: if $\Omega$ is a convex body containing $B(0,1)$,
then $\Omega$ can be covered by $O_D( \mes(\Omega) )$ translates of $[-1,1]^D$.  To see this, we use a covering argument of Ruzsa \cite{ruzsa-group}.  First observe that because the cube $[-1/2,1/2]^D$ is contained in a dilate of $B(0,1)$ (and hence $\Omega$) by $O_D(1)$, the Minkowski sum $\Omega + [-1/2,1/2]^D$ is also contained in an $O_D(1)$-dilate of $\Omega$ and thus has volume $O_D(\mes(\Omega))$.  Now let
$x_1 + [-1/2,1/2]^D, \ldots, x_N + [-1/2,1/2]^D$ be a maximal collection of disjoint shifted cubes with $x_1,\ldots,x_N \in \Omega$, then by the previous volume bound we have $N \ll_D \mes(\Omega)$.  But by maximality we see that the cubes
$x_1 + [-1,1]^D, \ldots, x_N + [-1,1]^D$ must cover $\Omega$, and the claim follows.
\end{proof}

\section{Counting points of varieties over $F_p$}\label{variety-sec}

Let $R$ be an arbitrary ring, and let $P_1,\ldots,P_J \in R[\x_1,\ldots,\x_D]$ be polynomials.  Our interest here is to control the ``density'' of the (affine) algebraic variety
$$ \{ (x_1,\ldots,x_D) \in R^D: P_j(x_1,\ldots,x_D) = 0 \hbox{ for all } 1 \leq j \leq J \}$$
and more precisely to estimate quantities such as
\begin{equation}\label{avg-def}
 \E_{ x_1 \in A_1, \ldots, x_D \in A_D } \prod_{j \in [J]} 1_{P_j(x_1,\ldots,x_D) = 0} 
\end{equation}
for certain finite non-empty subsets $A_1,\ldots,A_D \subset R$ (typically the $A_i$ will either be all of $R$, or some
arithmetic progression).  We are particularly interested in the case when $R$ is the finite field $F_p$, but in order to also encompass the case of the integers $\Z$ (and of polynomial rings over $F_p$ or $\Z$), we shall start by working in the more general context of a unique factorization domain.

Of course, the proper way to do this would be to use the tools of modern algebraic geometry, for instance using the concepts of generic point and algebraic dimension of varieties.  Indeed, the results in this appendix are ``morally trivial'' if one uses the fact that the codimension of an algebraic variety is preserved under restriction to generic subspaces.
However, to keep the exposition simple we have chosen a very classical, pedestrian and elementary approach, to emphasise that the facts from algebraic geometry which we will need are not very advanced.

From the factor theorem (which is valid over any unique factorization domain) we have

\begin{lemma}[Generic points of a one-dimensional polynomial]\label{factor-one} Let $P \in R[\x]$ be a polynomial of one variable of degree  at most $d$ over a ring $R$.  If $P \neq 0$, then $P(x) \neq 0$ for all but at most $d$ values of $x \in R$.
\end{lemma}

As a corollary, we obtain

\begin{corollary}[Generic points of a multi-dimensional polynomial]\label{factor-two} 
Let $P \in R[\x_1,\ldots,\x_D]$ be a polynomial of $D$ variables of degree at most $d$ over a ring $R$.  If $P \neq 0$, then $P(\cdot,x_D) \neq 0$ for all but at most $d$ values of $x_D \in R$, where $P(\cdot,x_D) \in F[\x_1,\ldots,\x_{D-1}]$ is the polynomial of $D-1$ variables formed from $P$ by replacing $\x_D$ by $x_D$.
\end{corollary}

\begin{proof} View $P$ as a one-dimensional polynomial of $\x_D$ with coefficients in the ring $R[\x_1,\ldots,\x_{D-1}]$ (which contains $R$), and apply Lemma \ref{factor-one}.
\end{proof}

As a consequence, we obtain a ``baby combinatorial Nullstellensatz'' (cf. \cite{alon}):

\begin{lemma}[Baby Nullstellensatz]\label{babynull}
Let $P \in R[\x_1,\ldots,\x_D]$ be a polynomial of $D$ variables of degree at most $d$ over a ring $R$.  Let $A_1,\ldots,A_D$ be finite subsets of $R$ with $|A_1|,\ldots,|A_D| \geq M$ for some $M > 0$. If $P \neq 0$, then
$$ \E_{x_1 \in A_1, \ldots, x_D \in A_D } 1_{P(x_1,\ldots,x_D)=0}  \leq \frac{Dd}{M} \ll_{D,d} \frac{1}{M}.$$
\end{lemma}

\begin{proof} We induct on $D$.  The case $D=0$ is vacuous.  Now suppose $D \geq 1$ and the claim has already been proven for $D-1$.  By Corollary \ref{factor-two}, we have $P(\cdot,x_D) \neq 0$ for all but at most $d$ values of $x_D \in A_D$.  The exceptional values of $x_D$ can contribute at most $\frac{d}{M}$, while the remaining values of $x_D$ will contribute at most $\frac{(D-1)d}{M}$ by the induction hypothesis.  This completes the induction.
\end{proof}

This gives us a reasonable upper bound on the quantity \eqref{avg-def} in the case $J=1$, which then trivially implies the same
bound for $J > 1$.  However, we expect to do better than $\frac{1}{M}$ type bounds for higher $J$ when the polynomials $P_1,\ldots,P_J$ are jointly coprime.  To exploit the property of being coprime we recall the classical \emph{resultant} of two polynomials.

\begin{definition}[Resultant]  Let $R$ be a ring, let $d, d' \geq 1$, and let
$$ P = a_0 + a_1 \x + \ldots + a_d \x^d; \quad Q = b_0 + b_1 + \ldots + b_{d'} \x^{d'}$$
be two polynomials in $R[\x]$ of degree at most $d,d'$ respectively.  
Then the \emph{resultant} $\Res_{d,d'}(P,Q) \in R$ is defined to be the determinant of the
$(d + d') \times (d + d')$ matrix whose rows are the coefficients in $R$ of the polynomials $P, \x P, \ldots, \x^{d'-1} P, Q, \x Q, \ldots, \x^{d-1} P$ with respect to the basis $1, \x, \ldots, \x^{d+d'-1}$.  

More generally, if $D \geq 1$, $1 \leq i \leq D$, $d_i, d'_i \geq 1$, and $P, Q$ are two polynomials in $R[\x_1,\ldots,\x_D]$ with $\deg_{\x_i}(P) \leq d_D$ and $\deg_{\x_i}(Q) \leq d'_i$, then we define the resultant $\Res_{d_D,d'_D,\x_D}(P,Q) \in 
R[\x_1,\ldots,\x_{i-1}, \x_{i+1},\ldots,\x_D]$ by viewing $P$ and $Q$ as one-dimensional polynomials of $\x_i$ over the ring
$R[\x_1,\ldots,\x_{i-1}, \x_{i+1},\ldots,\x_D]$ and using the one-dimensional resultant defined earlier.
\end{definition}

\begin{example} If $d=d'=1$, then the resultant of $a + b\x$ and $c+d\x$ is $ad-bc$, and the resultant of
$a(\x_1) + b(\x_1) \x_2$ and $c(\x_1) + d(\x_1) \x_2$ in the $\x_2$ variable is $a(\x_1) d(\x_1) - b(\x_1) c(\x_1)$.
\end{example}

Let $P,Q \in R[\x]$ have degrees $d, d'$ respectively for some $d,d' \geq 1$, where $R$ is a unique factorization domain.
By the determinant into a matrix and its adjugate, we obtain an identity
\begin{equation}\label{reso}
 \Res_{d,d'}(P,Q) = AP + BQ
\end{equation}
for some polynomials $A,B \in R[\x]$ of degree at most $d'-1$ and $d-1$ respectively.  Thus if $P,Q$ are irreducible and coprime, then the resultant cannot vanish by unique factorization.  The same extends to higher dimensions:

\begin{lemma}\label{coprime-resultant}  Let $R$ be a unique factorization domain, let $1 \leq i \leq D$, and suppose that $P, Q \in \R[\x_1,\ldots,\x_D]$ are such that $\deg_{\x_i}(P) = d_i \geq 1$ and $\deg_{\x_i}(Q) = d'_i \geq 1$.  If $P$ and $Q$ are irreducible and coprime,
then $\Res_{d,d',\x_i}(P,Q) \neq 0$.
\end{lemma}

\begin{proof} View $P,Q$ as one-dimensional polynomials over the unique factorization domain $R[\x_1,\ldots,\x_{i-1}, \x_{i+1},\ldots,\x_D]$ and apply the preceding argument.
\end{proof}

\begin{lemma}[Generic points of multiple polynomials]\label{coprime-descent}  Let $P_1, \ldots, P_J \in R[\x_1,\ldots,\x_D]$ have degrees at most $d$ over a unique factorization domain $R$, and suppose that all the $P_1,\ldots,P_J$ are non-zero and jointly coprime.  Then $P_1(\cdot,x_D)$, $\ldots$, $P_J(\cdot,x_D)$ are non-zero and jointly coprime for all but
at most $O_{J,d}(1)$ values of $x_D \in R$.
\end{lemma}

\begin{proof} 
By alignedting each of the $P_j$ into factors we may assume that all the $P_j$ are irreducible.  By eliminating any two polynomials which are scalar multiples of each other we may then assume that the $P_j$ are pairwise coprime.  The claim is vacuous for $k < 2$, so it will suffice to verify the claim for $k=2$.

Suppose first that $P_1$ is constant in $\x_D$.  Then $P_1(\cdot,x_D) = P_1$ is irreducible, and the only way it can fail to be coprime to $P_2(\cdot,x_D)$ is if $P_2(\cdot,x_D)$ is a multiple of $P_1$.  But we know that $P_2$ itself is not a multiple of $P_1$; viewing $P_2$ modulo $P_1$ as a polynomial of degree at most $d$ in $\x_D$ over the ring 
$R[\x_1,\ldots,\x_{D-1}]/(P_1)$ we see from Lemma \ref{factor-one} that the number of exceptional $x_D$ is at most $d$.

A similar argument works if $P_2$ is constant in $\x_D$.  So now we may assume that $\deg_{\x_D}(P_1) = d_1$ and $\deg_{\x_D}(P_2) = d_2$ for some $d_1,d_2 \geq 1$, which allows us to compute the resultant $\Res_{d_1,d_2,\x_D}(P_1,P_2) \in F[\x_1,\ldots,\x_{D-1}]$.  By Lemma \ref{coprime-resultant}, this resultant is non-zero; also by definition we see that the resultant has degree $O_{d}(1)$.

From \eqref{reso} we see that if $P_1(\cdot,x_D)$ and $P_2(\cdot,x_D)$ have any common factor in $R[\x_1,\ldots,\x_{D-1}]$ (which we may assume to be irreducible), then this factor must also divide $\Res_{d_1,d_2,\x_D}(P_1,P_2)$.  From degree considerations we see that there are at most $O_d(1)$ such factors.  Let $Q \in R[\x_1,\ldots,\x_{D-1}]$ be one such possible factor. Since $P_1$ and $P_2$ are coprime, we know that $Q$ cannot divide both $P_1$ and $P_2$; say it does not divide $P_2$.  Then by viewing $P_2$ modulo $Q$ as a polynomial of $y_D$ over $R[x_1,\ldots,x_{D-1}]/Q$ as before we see that there are at most $d$ values of $x_D$ for which $Q$ divides $P_2(\cdot,x_D)$.  Putting this all together we obtain the claim.
\end{proof}

This gives us a variant of Lemma \ref{babynull}.

\begin{lemma}[Second baby Nullstellensatz]\label{babynull-2}
Let $P_1,\ldots,P_J \in R[\x_1,\ldots,\x_D]$ be polynomials of $D$ variables of degree at most $d$ over a unique factorization domain $R$.  Let $A_1,\ldots,A_D$ be finite subsets of $F$ with $|A_1|,\ldots,|A_D| \geq M$ for some $M > 0$. If all the $P_1,\ldots,P_J$ are non-zero and jointly coprime, then
$$ \E_{ x_1 \in A_1, \ldots, x_D \in A_D} \prod_{j \in [J]} 1_{P_j(x_1,\ldots,x_D) \hbox{non-zero, jointly coprime}}   = 1 - O_{D,d,J}(\frac{1}{M})$$
and
$$ \E_{ x_1 \in A_1, \ldots, x_D \in A_D } \prod_{j \in [J]} 1_{P_j(x_1,\ldots,x_D) = 0}  \ll_{D,d,J} \frac{1}{M^2}.$$
\end{lemma}

\begin{remark} One could obtain sharper results by using Bezout's lemma, but the result here will suffice for our applications.
\end{remark}

\begin{proof} The first claim follows by repeating the proof of Lemma \ref{babynull} (replacing Corollary \ref{factor-two} by
Lemma \ref{coprime-descent} and we leave it to the reader.  To prove the second claim, we again induct on $D$.  The base case $D=0$ is again trivial, so assume $D \geq 1$ and the claim has already been proven for $D-1$.  

By Lemma \ref{coprime-descent}, we know that for all but $O_{J,d}(1)$ values of $x_D$, that the polynomials $P_1(\cdot,x_D), \ldots, P_J(\cdot,x_D)$ are all non-zero and coprime.  This case will contribute $O_{D,d,J}(\frac{1}{M^2})$ by the induction hypothesis.  Now consider one of the $O_{J,d}(1)$ exceptional values of $x_D$.  For each such $x_D$, at least one of the polynomials $P_j(\cdot,x_D)$ has to be non-zero, otherwise $\x_D-x_D$ would be a common factor of $P_1,\ldots,P_J$, a contradiction.  Applying Lemma \ref{babynull} we see that the contribution of each such $x_D$ is thus $O_{D,d,J}(\frac{1}{M^2})$, and the claim follows.
\end{proof}

We now specialize the above discussion to compute the local factors $c_p$ and $\overline{c_p}$
defined in Definition \ref{local-def}.  We first observe the following easy upper bounds:

\begin{lemma}[Crude local bound]\label{crude-cor}  Let $P_1,\ldots,P_J \in \Z[\x_1,\ldots,\x_D]$ have degree at most $d$,
and let $p$ be a prime.
\begin{itemize}
\item[(i)] If all the $P_1,\ldots,P_J$ vanish identically modulo $p$, then $c_p(P_1,\ldots,P_J) = 1$.
\item[(ii)] If at least one of $P_1,\ldots,P_J$ vanish identically modulo $p$, then $\overline{c_p}(P_1,\ldots,P_J) = 0$.
\item[(iii)] If at least one of $P_1,\ldots,P_J$ is a non-zero constant modulo $p$, then $c_p(P_1,\ldots,P_J) = 0$.
\item[(iv)] If at least one of $P_1,\ldots,P_J$ is non-constant modulo $p$, then $c_p(P_1,\ldots,P_J) \ll_{d,D} 1/p$.
\item[(v)] If the $P_1,\ldots,P_J$ are jointly coprime modulo $p$, then $c_p(P_1,\ldots,P_J) \ll_{d,D} 1/p^2$.
\end{itemize}
\end{lemma}

\begin{proof} (i), (ii), (iii) are trivial, while (iv) and (v) follow from 
Lemma \ref{babynull} and Lemma \ref{babynull-2} respectively (setting $A_1=\ldots=A_D=R=F_p$).
\end{proof}

Now we can refine the bound for a single polynomial $P$ in the case when $P$ is linear in one variable, with linear and constant coefficients coprime.

\begin{lemma}[Linear case]\label{linear-lemma}  Let $P \in \Z[\x_1,\ldots,\x_D]$ have degree at most $d$, and let $p$ be a prime.  Suppose that $P\ \mod\ p$ is linear in the $\x_i$ variable for some $1 \leq i \leq d$, thus we have 
$$ P(\x_1,\ldots,\x_D) = P_1(\x_1,\ldots,\x_{i-1}, \x_{i+1},\ldots,\x_{D}) \x_i + P_0(\x_1,\ldots,\x_{i-1},\x_{i+1},\ldots, \x_{D})\ \mod\ p$$
for some polynomials $P_0,P_1 \in F_p[\x_1,\ldots,\x_{i-1},\x_{i+1},\ldots,\x_D]$.  Suppose also that the linear coefficient $P_1$ is non-zero and coprime to the constant coefficient $P_0$.  Then $c_p(P) = 1/p + O_{d,D}(1/p^2)$.
\end{lemma}

\begin{proof}
Let us aligned $F_p^{D-1} = A \cup B \cup C$, where $A$ is the subset of $F_p^{D-1}$ where $P_1 \neq 0$, $B$ is the subset of $F_p^{D-1}$ where $P_1=0$ and $P_2 \neq 0$, and $C$ is the subset of $F_p^{D-1}$ where $P_1=P_2=0$.  Then an elementary counting argument shows that
$$ c_p(P) = \frac{|A| + |C| p}{p^D} = \frac{1}{p} - \frac{|B|+|C|}{p^D} + \frac{|C|}{p^{D-1}}.$$
Since $P_1$ is not zero, we see from Lemma \ref{crude-cor}(iv) that $|B|+|C| \ll_d p^{D-2}$.  Since $P_1,P_2$ are coprime modulo $p$, we see from Lemma \ref{crude-cor}(v) that $|C| \ll_d p^{D-2}$.  The claim follows.
\end{proof}

We can now quickly prove Lemma \ref{local-cor}.

\begin{proof}[Proof of Lemma \ref{local-cor}]  The claims (a), (b), (d) follow from Lemma \ref{crude-cor} and Definition \ref{bad-def}, while (c) follows from Lemma \ref{linear-lemma} and Definition \ref{bad-def}.  The claim (e) is trivial, and the claim (f) follows from (a), (b) and \eqref{exclude}.
\end{proof}

\section{The distribution of primes}\label{prime-sec}

In this section we recall some classical results about the distribution of primes.

For $\Re(s) > 1$, define the \emph{Riemann zeta function}
$$ \zeta(s) := \sum_{n=1}^\infty \frac{1}{n^s}.$$
Our argument will be elementary enough that we will not need the meromorphic continuation of $\zeta$ to
the region $\Re(s) \leq 1$. From the unique factorization of the natural numbers, we have the \emph{Euler product formula}
\begin{equation}\label{euler-formula}
 \zeta(s) = \prod_p (1 - \frac{1}{p^s})^{-1}.
\end{equation}
We also have the bounds
\begin{equation}\label{zeta-bounds}
\zeta(s) = \frac{1}{s-1} + O(\log(2 + |\Im(s)|) \hbox{ and } \frac{1}{\zeta(s)} = O(\log(2 + |\Im(s)|))
\end{equation}
whenever $1 < \Re(s) < 10$ (see e.g. \cite[Chapter 3]{titchmarsh}).  

From the prime number theorem 
$$ \sum_{p < x} 1 = (1 + o(1)) \frac{x}{\log x} \hbox{ as } x \to \infty$$
(which, incidentally, can be deduced readily from \eqref{zeta-bounds}),
and summation by parts, we easily obtain the estimates
\begin{align}
\sum_{p < x} \log p &= x + o(x) \hbox{ as } x \to \infty \label{logsum}\\
\sum_{p < x} \frac{1}{p} &= \log \log(10+x) + O(1) \hbox{ for } x > 0 \label{log-chebyshev}\\
\sum_{p < x} \frac{\log^K p}{p} &\ll_K \log^K(10 + x)  \hbox{ for } K > 0, x > 0 \label{plog-chebyshev}\\
|\sum_{p > x} \frac{\log^K p}{p^s}| &\ll_{K,s} \frac{\log^{K-1}(x)}{x^{\Re s-1}}  \hbox{ for } K \geq 0, x > 1, \Re(s) > 1.\label{psum-chebyshev}
\end{align}
In a similar spirit we have whenever $1 < \Re(s) < 2$ and $x > 2$
\begin{align*}
|\sum_{p > x} \log( 1 - \frac{1}{p^s} )| &\ll \sum_{p > x} \frac{1}{p^{\Re s}}  \\
&\ll \sum_{n=0}^\infty \sum_{2^n x < p \leq 2^{n+1} x} \frac{1}{2^{n\Re s} x^{\Re s}} \\
&\ll \sum_{n=0}^\infty \frac{ 1 }{ 2^{n(\Re s - 1)} x^{\Re s - 1} \log x} \\
&\ll \frac{1}{(\Re s - 1) x^{\Re s-1} \log x}.
\end{align*}
In particular, when $\Re(s) = 1+1/\log R$ and $x = R^{\log R}$ we have
$$ \sum_{p > R^{\log R}} \log( 1 - \frac{1}{p^s} ) = o(1)$$
and hence from \eqref{euler-formula}
\begin{equation}\label{euler-formula-2}
\prod_{p \leq R^{\log R}} (1 - \frac{1}{p^s})^{-1} = (1 + o(1)) \zeta(s) \hbox{ whenever } \Re(s) = 1 + 1/\log R.
\end{equation}

We will frequently encounter expressions of the form $\exp( K \sum_{p \in \primes} \frac{1}{p} )$, where $\primes$ ranges over some
set of primes (typically finite).  Such sums can eventually be somewhat large, thanks to \eqref{log-chebyshev}.  Fortunately 
the very slow nature of the divergence of $\sum_p \frac{1}{p}$ lets us estimate this exponential by a slowly divergent sum over primes, conceding only a few logarithms.

\begin{lemma}[Exponentials can be replaced by logarithms]\label{explog}  Let $\primes$ be any set of primes, and let $K \geq 1$.  Then
$$ \exp( K \sum_{p \in \primes} \frac{1}{p} ) \leq 1 + O_K( \sum_{p \in \primes} \frac{\log^{K} p}{p} )$$
or equivalently
$$ \EXP( K \sum_{p \in \primes} \frac{1}{p} ) \ll_K \sum_{p \in \primes} \frac{\log^{K} p}{p} .$$
\end{lemma}

\begin{remark} Note that the sum is only over primes in $\primes$, rather than products of primes in $\primes$, which would have been
the case if we had written $\exp( K \sum_{p \in \primes} \frac{1}{p} ) = \prod_{p \in \primes} \exp( K/p ) = \prod_{p \in \primes} (1 + O_K(1/p) )$.  The fact that we keep the sum over primes is useful for applications, as it allows us to work over fields $F_p$ rather than mere rings $\Z_N$ when performing certain local counting estimates.  This lets us avoid certain technical issues involving zero divisors which would otherwise complicate the argument.  The additional logarithmic powers of $p$ are sometimes dangerous, but in several cases we will be able to acquire an additional factor of $\frac{1}{p}$ from an averaging argument, which will make the summation on the right-hand side safely convergent regardless of how many logarithms are present, thanks to \eqref{psum-chebyshev}.
\end{remark}

\begin{proof} Let us fix $K$ and suppress the dependence of the $O()$ notation on $K$.  By a limiting argument we may take $\prod$ to be finite. We expand the left-hand side as a power series
$$ 1 + \sum_{n=1}^\infty \frac{K^n}{n!} \sum_{p_1,\ldots,p_n \in \primes} \frac{1}{p_1 \ldots p_n}.$$
By paying a factor of $n$ we may assume that $p_n$ is greater than or equal to the other primes, thus bounding the previous expression by
$$ 1 + \sum_{n=1}^\infty \frac{K^n}{(n-1)!} \sum_{p_n \in \primes} \sum_{p_1,\ldots,p_{n-1} \in \primes; p_1,\ldots,p_{n-1} \leq p_n} \frac{1}{p_1 \ldots p_n}.$$
We rewrite $p_n$ as $p$ and rearrange this as
$$ 1 + \sum_{p \in \primes} \frac{1}{p} \sum_{n=1}^\infty \frac{K^n}{(n-1)!} (\sum_{p' \in \primes: p' \leq p} \frac{1}{p'})^{n-1}.$$
From \eqref{log-chebyshev} we have
$$\sum_{p' \in \primes: p' \leq p} \frac{1}{p'} \leq \log \log(10+p) + O(1)$$
and so we can bound the previous expression by
$$ 1 + \sum_{p \in \primes} \frac{1}{p} \sum_{n=1}^\infty \frac{K ( K \log \log(10+p) + O(1) )^{n-1}}{(n-1)!}.$$
Summing the power series we obtain the result.
\end{proof}

Finally, we record a very simple lemma, using the quantities $w$ and $W$ defined in Section \ref{notation}.

\begin{lemma}[Divisor bound]\label{divbound}  Let $\primes$ be any collection of primes such that $\prod_{p \in \primes} p \leq M W^M$ for some $M > 0$.  Then
$$ \sum_{p \geq w: p \in \primes} \frac{1}{p} = o_M(1).$$
\end{lemma}

\begin{proof}  We trivially bound $\frac{1}{p}$ by $\frac{1}{w}$, and observe that the number of primes in $\primes$ larger than $w$ is at most $\log(M W^M)/\log(w) = o_M(\log W)$.  But from \eqref{logsum} we have $\log W \ll w$.  The claim follows.
\end{proof}

\end{document}